\documentclass{report}

\usepackage{amsmath}
\usepackage{amssymb}
\usepackage{amsthm}
\usepackage{tikz}
\usepackage{tikz-cd}
\usepackage{xcolor}
\usepackage{graphicx}
\usepackage{caption}
\usepackage{subcaption}
\usepackage{mathrsfs}

\newtheorem{theorem}{Theorem}[chapter]

\newtheorem{definition}[theorem]{Definition}
\newtheorem{example}[theorem]{Example}

\newtheorem{proposition}[theorem]{Proposition}

\title{A Geometric Approach to Structure-Preserving Integrators for Mechanical Systems}
\author{Viyom Vivek \thanks{Centre for Systems and Control, IIT Bombay, India}, David Martin de Diego \thanks{Instituto de Ciencias Matematicas, Madrid, Spain} and Ravi N. Banavar \thanks{Centre for Systems and Control, IIT Bombay, India}}

\usepackage{hyperref}

\begin{document}

\maketitle

\begin{abstract}
We develop a geometric framework for the numerical integration of mechanical systems evolving on manifolds. After briefly reviewing classical numerical methods and highlighting their geometric limitations, we include a short interlude reviewing the differential geometric notions required in the sequel. We then introduce retraction maps as intrinsic generalizations of the Riemannian exponential, which induce discretization maps tailored to manifold-valued dynamics. Adopting the Tulczyjew unified viewpoint, mechanical systems are formulated as Lagrangian submanifolds, providing a natural and coordinate-free foundation for the construction of structure-preserving integrators for both Hamiltonian and Lagrangian systems. The framework is specialized to Lie groups, where parallelizability allows for the global trivialization of tangent and cotangent bundles and the systematic derivation of integrators for Euler–Poincaré and Lie–Poisson equations. The effectiveness of the proposed approach is illustrated through the rigid body and heavy top, and is further extended to the construction of a geometric integrator for an underactuated mechanical system—a quadrotor—demonstrating the applicability of the framework beyond fully symmetric systems and toward problems relevant in robotics and control.
\end{abstract}

\tableofcontents

\chapter{A Brief Introduction to Numerical Integrators}
\label{ch1} 
Our exposition starts with the definition of numerical integrators, emphasizing their utility and their foundational role in the discretization of continuous dynamical models. We then provide a concise survey of classical methods, illustrating their limitations through representative examples from mechanics, and conclude this introductory discussion with a brief literature survey that situates the present work within the broader development of geometric and structure-preserving numerical methods.
\section{What are numerical integrators and why do we need them?}
\label{sec1.1}
Physical laws are fundamentally concerned with describing the temporal evolution of a system’s state, whether it be the temperature in thermodynamics, the phase-space variables in classical mechanics, or the wavefunction in quantum mechanics. Mathematically, these laws are formulated in the language of differential equations. In thermodynamics, the temporal evolution of temperature, pressure, and entropy is governed by a system of first-order partial differential equations (PDEs) encoded in the laws of thermodynamics. Likewise, the evolution of positions  in classical mechanics is described by a system of second-order ordinary differential equations (SODEs), as expressed by Newton’s laws (or, equivalently, by Hamilton’s principle, discussed later) of first order if we adopt the Hamiltonian formulation in terms of positions and momenta. In quantum mechanics, the wavefunction evolves according to the Schr{\"o}dinger equation, which is a second-order partial differential equation. In what follows, we confine our discussion to ordinary differential equations associated with classical mechanical systems.

Except for a limited class of simplified models, obtaining analytical solutions to such differential equations is extremely challenging and, in many cases, impossible; consequently, one must resort to approximate solution techniques. This is typically achieved by discretizing the differential equation into a corresponding difference equation, which can then be solved iteratively using computational methods. The above procedure is termed numerical integration, with the associated difference equation being called a numerical integrator.

Let us look at the following dynamical system:
\begin{equation}
    \dot{x}(t) = f(x(t)) \quad \text{given} \quad x(t_0) = x_0 \label{1.1}
\end{equation}
where $x : [t_0, t_f] \to \mathbb{R}^n$ denotes the solution trajectory and $f : \mathbb{R}^n \to \mathbb{R}^n$ is a sufficiently smooth vector-valued function governing the dynamics. We initiate the discretization procedure by partitioning the time interval $[t_0, t_f] \subset \mathbb{R}$ into $N \in \mathbb{Z}$ subintervals of equal length $h := (t_f - t_0) / N$, commonly referred to as the step size. The discrete, or approximate, solution trajectory of \eqref{1.1} can be obtained by iteratively solving the following difference equation:
\begin{equation}
    F(x_k, x_{k+1}) = 0 \quad \text{given} \quad x_0 \label{1.2}
\end{equation}
where $x_k := x(kh)$ for $k \in \{0, \dots, N\}$ and $F : \mathbb{R}^n \times \mathbb{R}^n \to \mathbb{R}^n$ is a function constructed to define the numerical integrator. This function $F$ also depends on the parameter $h$, the step size. The performance of a numerical integrator is governed by the choice of the function $F$ along with the step size $h$ and is measured by the discrepancy between the continuous solution of \eqref{1.1} and the discrete solution of \eqref{1.2}.
\section{A quick overview of classical methods and their limitations}
\label{sec1.2}
We now present several classical methods for the numerical integration of dynamical systems, illustrated through examples from mechanics, and highlight their inherent limitations.
\subsection{Classical Harmonic Oscillator}
\label{harmonic oscillator}
\begin{figure}[h]
\centering
\begin{subfigure}{.5\textwidth}
  \centering
  \includegraphics[width=1\linewidth]{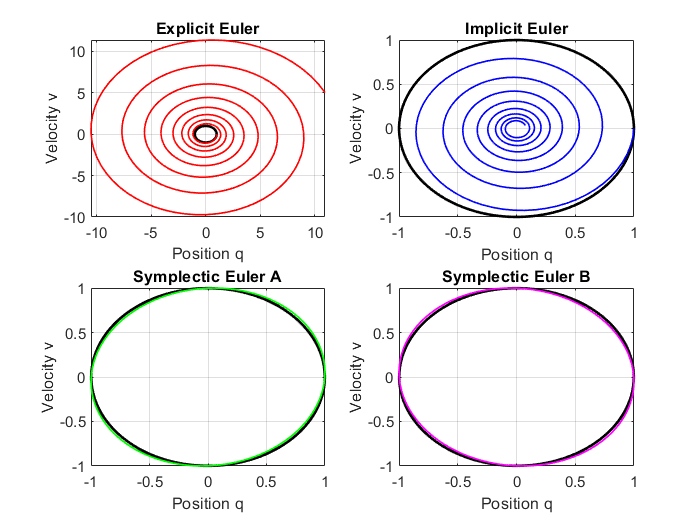}
  \caption{}
\end{subfigure}%
\begin{subfigure}{.5\textwidth}
  \centering
  \includegraphics[width=1\linewidth]{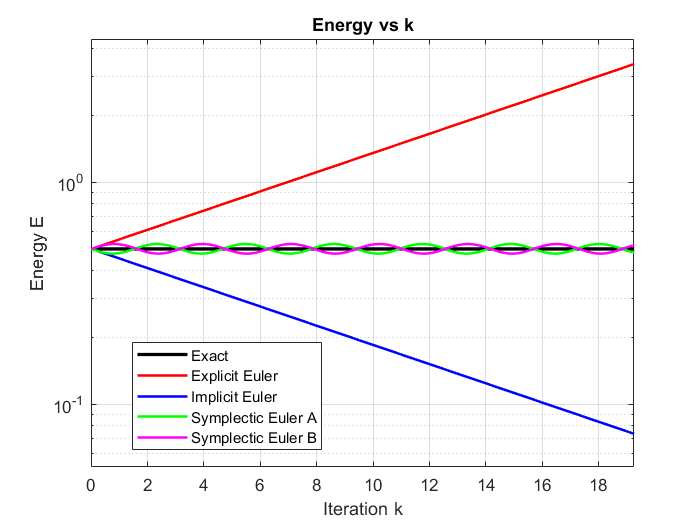}
  \caption{}
\end{subfigure}
\caption{(a) Comparison of various numerical integrators for the classical harmonic oscillator with parameter values (in standard units) $k/m=1$, $h=0.1$, and initial condition $x_0=(1,0)$ for 50 iterations. The exact solution, which traces an elliptical trajectory, is shown in black. (b) Energy versus iteration for various numerical integrators applied to the classical harmonic oscillator. The energy is shown on a logarithmic scale.}
\label{fig1.1}
\end{figure}

\begin{center}
    \begin{tikzpicture}
    \draw[ultra thick] (1,1)--(1,2);
    \foreach \i in {0,0.8,1.6,2.4,3.2,4}
    {\draw (1+\i,1.5)--(1.2+\i,1.8)--(1.6+\i,1.2)--(1.8+\i,1.5);}
    \draw[fill] (6,1.5) circle (7pt);
    \draw[dashed] (4,2.2)--(4,0.8);
    \draw[->] (4,0.9)--(6,0.9);
    \node at (5,0.6) {$q$};
    \node at (6,2) {$m$};
    \node at (3.5,2.2) {$k$};
    \end{tikzpicture}
\end{center}
Consider a point mass $m$ attached to a spring with stiffness $k$, as illustrated in the figure above. The equations of motion are governed by Newton’s second law:
\begin{equation}
m \ddot{q} = -kq \label{1.3}
\end{equation}
where $q \in \mathbb{R}$ denotes the displacement from equilibrium, and the restoring force on the right-hand side is given by Hooke’s law. We first rewrite \eqref{1.3} as a first-order dynamical system:
\begin{align}
x := \begin{pmatrix}
q \\
v
\end{pmatrix} \implies
\dot{x} =
\begin{pmatrix}
    \dot{q} \\
    \dot{v}
\end{pmatrix} = \begin{pmatrix}
    v \\
    -kq/m
\end{pmatrix} =: f(x) \label{1.4}
\end{align}
where $x \in \mathbb{R}^2$ represents the state of the system in the state space $\mathbb{R}^2$, comprising both the position $q$ and the velocity $v := \dot{q}$. Since the vector-valued function $f$ in \eqref{1.4} is linear, the system admits a closed-form analytical solution, given by:
\begin{align}
x(t) = \begin{pmatrix}
\cos{\sqrt{k/m} \, t} & \sqrt{m/k} \,  \sin{\sqrt{k/m} \,t} \\
-\sqrt{k/m} \, \sin{\sqrt{k/m} \, t} & \cos{\sqrt{k/m} \, t}
\end{pmatrix} x(0). \label{1.5}
\end{align}
Consequently, the state trajectory traces elliptical curves in the state space. A fundamental physical invariant associated with mechanical systems is the total energy:
\begin{align}
    E := \frac{1}{2} m v^2 + \frac{1}{2} k q^2 \label{1.6}
\end{align}
which consists of both kinetic and potential energy contributions. The conservation of this quantity follows from a direct verification:
\begin{align*}
    \dot{E} = m v \dot{v} + k q \dot{q} =-kvq + kqv = 0 
\end{align*}
where we have substituted the expressions for $\dot{q}$ and $\dot{v}$ from \eqref{1.4}. Let us now introduce some classical methods for the numerical integration of dynamical systems.
\paragraph{Explicit Euler method} 
This method represents the earliest approach to the numerical integration of dynamical systems, with its origins dating back to the mid-eighteenth century. The derivative in \eqref{1.1} is approximated by a forward finite difference according to
\begin{align*}
\frac{x_{k+1}-x_k}{h} = f(x_k) \implies x_{k+1} = x_k + h f(x_k).
\end{align*}
The method is termed explicit since $x_{k+1}$ can be expressed directly in terms of $x_k$. Moreover, it is a single-step method, as only the current state is required to compute the next one, and it is first order, with the error proportional to the step size $h$.
\paragraph{Implicit Euler method}
In this method, the derivative in \eqref{1.1} is approximated by a backward finite difference, given by
\begin{align*}
    \frac{x_k-x_{k-1}}{h} = f(x_k) \implies x_{k+1} = x_k + h f(x_{k+1}).
\end{align*}
The method is termed implicit since $x_{k+1}$ cannot be expressed explicitly in terms of $x_k$. Consequently, an implicit algebraic equation must be solved for $x_{k+1}$ at each iteration, typically using a Newton–Raphson–type method. As with the explicit Euler scheme, this is a single-step, first-order method. 
\vspace{10pt} \\
We first construct numerical integrators for \eqref{1.4} using the explicit and implicit Euler methods defined above: 
\begin{align}
    &x_{k+1} = x_k + h f(x_k) \implies x_{k+1} 
    = \begin{pmatrix}
        1 & h \\
        -hk/m & 1
    \end{pmatrix}x_k \quad \text{and} \label{explicit euler harmonic oscillator} \\
    &x_{k+1} = x_k + h f(x_{k+1}) \implies x_{k+1} 
     = \frac{1}{1+h^2k/m} \begin{pmatrix}
        1 & -h \\
        hk/m & 1
    \end{pmatrix}x_k, \label{implicit euler harmonic oscillator}
\end{align}
respectively. We know present the first example of a structure-preserving numerical integrator.
\paragraph{Symplectic Euler method} 
\label{para symplectic Euler}
This method constitutes one of the first structure-preserving numerical integration scheme to have been developed. It applies to systems whose state variables can be decomposed into position and velocity (or momentum) like components, namely,
\begin{align*}
    \dot{x} = f(x) \implies \begin{pmatrix}
        \dot{q} \\
        \dot{v}
    \end{pmatrix}=\begin{pmatrix}
        f_1(q,v) \\
        f_2(q,v)
    \end{pmatrix}.
\end{align*}
The method is commonly referred to as the semi-implicit or semi-explicit Euler scheme, as it combines an explicit Euler update for the position-like variables with an implicit Euler update for the velocity-like variables (symplectic Euler A): 
\begin{align*}
    \frac{q_{k+1}-q_k}{h} &= f_1(q_k,v_{k+1}) \implies q_{k+1} = q_k + h f_1(q_k,v_{k+1}) \quad \text{and} \\
    \frac{v_{k+1}-v_k}{h} &= f_2(q_k,v_{k+1}) \implies v_{k+1} = v_k + h f_2(q_k,v_{k+1})
\end{align*}
or, alternatively, an implicit Euler update for the position-like variables together with an explicit Euler update for the velocity-like variables (symplectic Euler B): 
\begin{align*}
    \frac{q_{k+1}-q_k}{h} &= f_1(q_{k+1},v_k) \implies q_{k+1} = q_k + h f_1(q_{k+1},v_k) \quad \text{and} \\
    \frac{v_{k+1}-v_k}{h} &= f_2(q_{k+1},v_k) \implies v_{k+1} = v_k + h f_2(q_{k+1},v_k).
\end{align*}
As with the standard Euler schemes, this method is a single-step, first-order integrator.
\vspace{10pt} \\
We next assess the performance of the integrators defined in \eqref{explicit euler harmonic oscillator} and \eqref{implicit euler harmonic oscillator} against those derived from the symplectic Euler A and B schemes:
\begin{align}
    &\begin{pmatrix}
        q_{k+1} \\
        v_{k+1}
    \end{pmatrix} = \begin{pmatrix}
        q_k + h v_{k+1} \\
        v_k - hk q_k/m
    \end{pmatrix} \implies x_{k+1} = \begin{pmatrix}
        1-h^2k/m & h \\
        -hk/m & 1
    \end{pmatrix}x_k \quad \text{and} \label{symplectic euler a harmonic oscillator} \\
    &\begin{pmatrix}
        q_{k+1} \\
        v_{k+1}
    \end{pmatrix} = \begin{pmatrix}
        q_k + h v_k \\
        v_k - hk q_{k+1}/m
    \end{pmatrix} \implies x_{k+1} = \begin{pmatrix}
        1 & h \\
        -hk/m & 1-h^2k/m
    \end{pmatrix}x_k, \label{symplectic euler b harmonic oscillator}
\end{align}
respectively. Consider the plots shown in figure \ref{fig1.1}. The explicit Euler method artificially injects energy into the system, causing the numerical trajectory to spiral outward. This behavior is reflected in the determinant of the update matrix in \eqref{explicit euler harmonic oscillator}, which satisfies $1+h^2\lambda^2 > 1$. In contrast, the implicit Euler method dissipates energy, leading to an inward-spiraling trajectory; correspondingly, the determinant of the matrix in \eqref{implicit euler harmonic oscillator} is $1/(1+h^2\lambda^2) < 1$. The symplectic Euler methods, on the other hand, are able to reproduce the qualitative features of the true dynamics. This is further evidenced by the fact that the update matrices in \eqref{symplectic euler a harmonic oscillator} and \eqref{symplectic euler b harmonic oscillator} have unit determinant. Note that the symplectic Euler methods do not preserve the energy of the system as shown in figure \ref{fig1.1}; however, the energy remains bounded, resulting in periodic behavior. The quantity that is truly preserved by these methods is the symplectic structure, which we define subsequently.
\vspace{-10pt}
\paragraph{Linear Symplectic Structure}
\label{para linear symplectic structure}
Let us reconsider the equations of motion of the spring–mass system from the perspective of Hamiltonian mechanics. We define the momentum as $p := m \dot{q}$ and introduce the energy function, which in position–momentum coordinates, commonly referred to as phase space in physics, is called the Hamiltonian of the system:
\begin{align*}
    H := \frac{1}{2m} p^2 + \frac{1}{2} k q^2.
\end{align*}
The dynamics are then governed by
\begin{align*}
    z := \begin{pmatrix}
        q \\
        p
    \end{pmatrix} &\implies \dot{z} = \begin{pmatrix}
        \dot{q} \\
        \dot{p}
    \end{pmatrix} = \begin{pmatrix}
        p/m \\
        -kq
    \end{pmatrix} = \begin{pmatrix}
        \partial H / \partial p \\
        -\partial H / \partial q
    \end{pmatrix} = \mathbb{J} \nabla H(z) 
\end{align*}
where $\mathbb{J} := \left(\begin{smallmatrix}
    0 & 1\\
    -1 & 0
\end{smallmatrix}\right)$ and $\nabla := \frac{\partial}{\partial z} = \left(\begin{smallmatrix}
    \partial/\partial q \\
    \partial/\partial p
\end{smallmatrix} \right)$ denotes the gradient operator. The operator $\mathbb{J} \nabla := \left(\begin{smallmatrix}
    \partial/\partial p \\
    -\partial/\partial q
\end{smallmatrix}\right)$, sometimes referred to as the symplectic gradient, underlies the conservative nature of the dynamics, as illustrated in the figure below.
\begin{center}
    \begin{tikzpicture}
\draw[<->] (1,3)--(7,3);
\draw[<->] (4,1)--(4,5);
\node at (7.3,3) {$q$};
\node at (4,5.3) {$p$};
\draw[red,thick] (4,3) ellipse (2 and 1);
\node at (6,2) {\footnotesize $\textcolor{red}{H=c}$};
\node at (4.8,3.665) {\footnotesize $z$};
\draw[-stealth,ultra thick,red] (5,3.865)--(5.2,4.5);
\draw[-stealth,ultra thick,red] (3,3.865)--(2.8,4.5);
\draw[-stealth,ultra thick,red] (5,2.135)--(5.2,1.5);
\draw[-stealth,ultra thick,red] (3,2.135)--(2.8,1.5);
\draw[-stealth,ultra thick,red] (4,4)--(4,4.7);
\draw[-stealth,ultra thick,red] (4,2)--(4,1.3);
\draw[-stealth,ultra thick,red] (2,3)--(1.3,3);
\draw[-stealth,ultra thick,red] (6,3)--(6.7,3);
\node at (5.4,4.7) {\footnotesize $\textcolor{red}{\nabla H (z)}$};
\draw[-stealth,ultra thick,orange] (5,3.865)--(5.635,3.665);
\draw[-stealth,ultra thick,orange] (3,3.865)--(3.635,4.065);
\draw[-stealth,ultra thick,orange] (5,2.135)--(4.365,1.925);
\draw[-stealth,ultra thick,orange] (3,2.135)--(2.265,2.335);
\draw[-stealth,ultra thick,orange] (4,4)--(4.7,4);
\draw[-stealth,ultra thick,orange] (4,4)--(4.7,4);
\draw[-stealth,ultra thick,orange] (6,3)--(6,2.3);
\draw[-stealth,ultra thick,orange] (4,2)--(3.3,2);
\draw[-stealth,ultra thick,orange] (2,3)--(2,3.7);
\node at (6.3,3.865) {\footnotesize $\textcolor{orange}{\mathbb{J} \nabla H (z)}$};
\draw[fill] (5,3.865) circle (1.5pt);
\end{tikzpicture}
\end{center}
On a constant energy level set (red ellipse), the gradient vector (red arrows) is normal to the level set at every point. In contrast, the symplectic gradient (orange arrows) is tangent to the level set everywhere, since the operator $\mathbb{J}$ acts as a $90^\circ$ rotation, analogous to multiplication by the imaginary unit $i$ in the complex plane. This analogy is not coincidental and reflects the deep connections between symplectic and complex structures which led to the discovery of K{\"a}hler structures. Consequently, the phase trajectory remains confined to a level set of the energy.

Through the matrix $\mathbb{J}$, the phase space
$\mathbb{R}^2$ is equipped with a nondegenerate skew-symmetric bilinear form, referred to as the symplectic structure. We now turn to a precise formulation of this notion following \cite{da2001lectures, berndt2001introduction}.
\begin{definition}
    Let $V$ be a finite-dimensional real vector space. A \emph{symplectic form} on $V$ is a bilinear map
    \begin{align*}
        \omega : V \times V \to \mathbb{R}
    \end{align*}
    that is
    \begin{itemize}
        \item \emph{skew-symmetric}, i.e. $\omega(u,v) = -\omega(v,u)$ for all $u,v \in V$, and
        \item \emph{non-degenerate}, i.e. if $\omega(u,v) = 0$ for all $v \in V$, then $u = 0$.
    \end{itemize}
    The pair $(V, \omega)$ is called a \emph{symplectic vector space}.
\end{definition}
Notice that the non-degeneracy of $\omega$ forces the vector space $V$ to be even-dimensional. The canonical example of a symplectic vector space is $(\mathbb{R}^{2n}, \omega_{std})$ where 
\begin{align}
    \omega_{std}(u,v) := u^T\mathbb{J}v \quad \text{where} \quad \mathbb{J} := \begin{pmatrix}
        0_n & I_n \\
        -I_n & 0_n
    \end{pmatrix}. \label{1.11}
\end{align}
This example is called canonical because of a fundamental result in symplectic geometry, the linear Darboux theorem, which states that for any $2n-$dimensional symplectic vector space $(V,\omega)$, there exists a basis in which the matrix representation of $\omega$ coincides with \eqref{1.11}. Next, we define a mapping that preserves this structure.
\begin{definition}
    Let $(V_1,\omega_1)$ and $(V_2,\omega_2)$ be symplectic vector spaces. A linear map
    \begin{align*}
        T : V_1 \to V_2
    \end{align*}
    is called \emph{symplectic} if it preserves the symplectic forms, i.e.,
    \begin{align*}
        \omega_2(Tu,Tv) = \omega_1(u,v) \quad \text{for all} \quad u,v \in V_1.
    \end{align*}
    Equivalently, this condition can be written as 
    \begin{align*}
        T^*\omega_2 = \omega_1
    \end{align*}
    where $T^*$ denotes the pullback of bilinear forms by $T$.
\end{definition}
In particular, the set of all symplectic automorphisms of a symplectic vector space $(V,\omega)$ forms a group under composition, called the symplectic group, and is denoted by  $Sp(V,\omega)$. In the case of the standard symplectic vector space $(\mathbb{R}^{2n},\omega_{std})$, this group is written as
\begin{align*}
    Sp(2n,\mathbb{R}) := \{A \in \mathbb{R}^{2n \times 2n} \, | \, A^T\mathbb{J}A = \mathbb{J} \}
\end{align*}
where $\mathbb{J}$ is the standard symplectic matrix defined in \eqref{1.11}. As the notion of a Lagrangian submanifold will play a central role in our later constructions, we first introduce its linear analogue, namely that of a Lagrangian subspace.
\begin{definition}
    Let $(V, \omega)$ be a symplectic vector space and let $L \subseteq V$ be a linear subspace. The \emph{symplectic orthogonal complement} of $L$ is defined by
    \begin{align*}
        L^\omega := \{ u \in V \,|\, \omega(u,v) = 0 \,\, \text{for every} \,\, v \in L \}.
    \end{align*}
\end{definition}
%
Unlike the Euclidean orthogonal complement, the symplectic orthogonal complement need not be complementary to $L$. In particular, the following cases may occur:
\begin{itemize}
    \item $L \subseteq L^\omega$ (\emph{isotropic}),
    \item $L^\omega \subseteq L$ (\emph{coisotropic}),
    \item $L = L^\omega$ (\emph{Lagrangian}),
    \item $L \cap L^\omega = \{0\}$ (\emph{symplectic}).
\end{itemize}
A subspace $L$ is Lagrangian if and only if it is maximally isotropic, i.e., $\text{dim} \, L = \frac{1}{2} \, \text{dim} \, V$. This case is of particular importance and will be the primary focus in later chapters. For example, in the symplectic vector space $(\mathbb{R}^{2n},\omega_{std})$, the subspaces
\begin{align*}
    \{(q,p) \,|\, p=0\} \quad \text{and} \quad \{(q,p) \,|\, q=0\}
\end{align*}
consisting of pure positions and pure momenta, respectively, are Lagrangian subspaces of the phase space.
\vspace{10pt} \\
We can now justify the terminology symplectic Euler methods by verifying that the update matrices in \eqref{symplectic euler a harmonic oscillator} and \eqref{symplectic euler b harmonic oscillator} indeed belong to $Sp(2,\mathbb{R})$.
\paragraph{Poisson Structure}
\label{para Poisson structure}
Symplectic geometry arises as a special case of a more general framework known as Poisson geometry, which provides a natural bridge between classical and quantum mechanics. To motivate this structure, recall the equations of motion of the spring-mass system. Physical observables such as position, momentum, and energy are represented by smooth functions on phase space.

Let $F : \mathbb{R}^2 \to \mathbb{R}$ be an observable. Its time evolution along a Hamiltonian trajectory satisfies
\begin{align*}
    \frac{d}{dt} F(q,p)&=\frac{\partial F}{\partial q}\dot{q}+\frac{\partial F}{\partial p}\dot{p}
    = \frac{\partial F}{\partial q}\frac{\partial  H}{\partial p}-\frac{\partial F}{\partial p}\frac{\partial H}{\partial q} \\
    \implies \dot{F}(z) &= \nabla F(z)^T \mathbb{J} \nabla H(z)
\end{align*}
where $z:=(q,p) \in \mathbb{R}^2$, $H: \mathbb{R}^2 \to \mathbb{R}$ is the Hamiltonian and $\mathbb{J}:=\left(\begin{smallmatrix}
    0 & 1\\
    -1 & 0
\end{smallmatrix}\right)$. The bilinear operation $\nabla(\cdot)^T\mathbb{J}\nabla(\cdot)$ on the space of smooth functions on $\mathbb{R}^2$, denoted by $C^\infty(\mathbb{R}^2)$, is called the Poisson bracket, denoted by $\{\cdot \,,\cdot\}$. It encodes the dynamics through the relation
\begin{align*}
    \dot{F} = \{F,H\},
\end{align*}
thereby expressing time evolution as a derivation on the algebra of observables, as illustrated in the figure below.
\begin{center}
\begin{tikzpicture}
\draw[<->] (1,3)--(7,3);
\draw[<->] (4,1)--(4,5);
\node at (7.3,3) {$q$};
\node at (4,5.3) {$p$};
\draw[red,thick] (4,3) ellipse (2 and 1);
\node at (6,2) {\footnotesize $\textcolor{red}{H=c_1}$};
\node at (4.8,3.665) {\footnotesize $z$};
\draw[-stealth,ultra thick,red] (5,3.865)--(5.2,4.5);
\draw[-stealth,ultra thick,red] (3,3.865)--(2.8,4.5);
\draw[-stealth,ultra thick,red] (5,2.135)--(5.2,1.5);
\draw[-stealth,ultra thick,red] (3,2.135)--(2.8,1.5);
\draw[-stealth,ultra thick,red] (4,4)--(4,4.7);
\draw[-stealth,ultra thick,red] (4,2)--(4,1.3);
\draw[-stealth,ultra thick,red] (2,3)--(1.3,3);
\draw[-stealth,ultra thick,red] (6,3)--(6.7,3);
\node at (5.4,4.6) {\footnotesize $\textcolor{red}{\nabla H (z)}$};
\draw[-stealth,ultra thick,orange] (5,3.865)--(5.635,3.665);
\draw[-stealth,ultra thick,orange] (3,3.865)--(3.635,4.065);
\draw[-stealth,ultra thick,orange] (5,2.135)--(4.365,1.925);
\draw[-stealth,ultra thick,orange] (3,2.135)--(2.265,2.335);
\draw[-stealth,ultra thick,orange] (4,4)--(4.7,4);
\draw[-stealth,ultra thick,orange] (4,4)--(4.7,4);
\draw[-stealth,ultra thick,orange] (6,3)--(6,2.3);
\draw[-stealth,ultra thick,orange] (4,2)--(3.3,2);
\draw[-stealth,ultra thick,orange] (2,3)--(2,3.7);
\node at (6.3,3.7) {\footnotesize $\textcolor{orange}{\mathbb{J} \nabla H (z)}$};
\draw[thick,blue] (1.5,5) to[out=-60,in=-120] (6.5,5);
\draw[-stealth,ultra thick,blue] (5,3.865)--(5.2,3.4);
\draw[-stealth,ultra thick,blue] (3,3.865)--(2.8,3.4);
\draw[-stealth,ultra thick,blue] (4,3.75)--(4,3.2);
\draw[-stealth,ultra thick,blue] (2,4.4)--(1.6,4);
\draw[-stealth,ultra thick,blue] (6,4.4)--(6.4,4);
\node at (5.3,3.2) {\footnotesize $\textcolor{blue}{\nabla F (z)}$};
\node at (6.5,5.3) {\footnotesize $\textcolor{blue}{F=c_2}$};
\draw[fill] (5,3.865) circle (1.5pt);
\end{tikzpicture}
\end{center}
Along a level set of the observable (blue curve), the gradient vector (blue arrows) indicates the direction of maximal increase. Its component along the symplectic gradient of the Hamiltonian (orange arrows) yields the rate of change of the observable along the Hamiltonian flow. We now formalize this notion following \cite{vaisman1994lectures, crainic2021lectures}.
\begin{definition}
    Let $A$ be an associative algebra. A \emph{Poisson bracket} on $A$ is a bilinear map
    \begin{equation*}
        \{\cdot \,, \cdot\} : A \times A \to A
    \end{equation*}
    that satisfies
    \begin{itemize}
        \item \emph{skew-symmetry}, i.e. $\{F,G\}=-\{G,F\}$,
        \item \emph{Jacobi identity}, i.e. $\{F,\{G,H\}\}+\{G,\{H,F\}\}+\{H,\{F,G\}\}=0$, and
        \item \emph{Leibniz identity}, i.e. $\{FG,H\}=F\{G,H\}+\{F,H\}G$
    \end{itemize}
    for all $F,G,H \in A$. The pair $(A.\{\cdot\,,\cdot\})$ is called a \emph{Poisson algebra}.
\end{definition}
In our context, $A$ will usually be the algebra of smooth functions $C^\infty(V)$ on a vector space $V$. The same definition extends verbatim to the case when $V=M$ is a smooth manifold.  Note that, unlike symplectic structures, Poisson structures can be defined in both odd and even dimensions. We now define maps that preserves this structure.
\begin{definition}
    Let $(A_1,\{\cdot \,, \cdot\}_1)$ and $(A_2,\{\cdot \,, \cdot\}_2)$ be Poisson algebras. A linear map
    \begin{equation*}
        T : A_1 \to A_2
    \end{equation*}
    is called a \emph{Poisson map} if it preserves the Poisson bracket, that is,
    \begin{align*}
        \{F,G\}_2 \circ T = \{F \circ T,G \circ T\}_1 \quad \text{for all} \quad F,G \in A_2.
    \end{align*}
\end{definition}
Another important notion associated with Poisson structures is the concept of Casimir elements. These encode constants of motion in mechanics and play a central role in the geometry of Poisson manifolds.
\begin{definition}
    Let $(A,\{\cdot\,,\cdot\})$ be a Poisson algebra. An element $C \in A$ is called a \emph{Casimir element} if 
    $$\{C,F\}=0 \quad \text{for all} \quad F \in A.$$
\end{definition}
\subsection{Kepler Problem}
\label{Kepler problem}
\begin{center}
    \begin{tikzpicture}
        \draw[fill] (1,1) circle (1pt);
        \node at (0.8,0.8) {o};
        \draw[dotted] (1,1)--(6,2);
        \node at (4,1.2) {$r_2$};
        \draw[dotted] (1,1)--(2,4);
        \node at (1.2,3) {$r_1$};
        \draw[dotted] (2,4)--(6,2);
        \node at (4,3.5) {$r$};
        \draw[fill] (4.5,2.765) circle (1pt);
        \draw[dotted] (1,1)--(4.5,2.765);
        \node at (3,2.5) {$R$};
        \draw[red,ultra thick,-stealth] (2,4)--(3,3.5);
        \draw[red,ultra thick,-stealth] (6,2)--(5,2.5);
        \draw[fill] (2,4) circle (3pt);
        \node at (2,4.4) {$m_1$};
        \draw[fill] (6,2) circle (4pt);
        \node at (6.6,2) {$m_2$};
    \end{tikzpicture}
\end{center}
Consider two point masses $m_1$ and $m_2$ moving in three-dimensional space under their mutual gravitational interaction, see \cite{arnold2007mathematical}. Let $r_1 \in \mathbb{R}^3$ and $r_2 \in \mathbb{R}^3$ denote the positions of $m_1$ and $m_2$, respectively, with respect to an inertial reference frame, as illustrated in the figure above. According to Newton’s second law, their equations of motion are given by:
\begin{align}
    m_1 \ddot{r}_1 &= \frac{Gm_1m_1}{\|r_2-r_1\|^3} (r_2-r_1) \quad \text{and} \label{1.12} \\
    m_2 \ddot{r}_2 &= -\frac{Gm_1m_2}{\|r_2-r_1\|^3}(r_2-r_1) \label{1.13}
\end{align}
where the gravitational force on the right-hand side is a central force obeying the inverse square law. We know introduce the position of the centre of mass, $R := \frac{m_1 r_1 + m_2 r_2}{m_1 + m_2}$, and the relative displacement vector between the two masses, $r := r_2 -r_1$. By adding and subtracting equations \eqref{1.12} and \eqref{1.13}, we obtain
\begin{align}
    \ddot{R} &= 0 \quad \text{and} \label{1.14} \\
    \ddot{r} &= -\frac{G(m_1 + m_2)}{\|r\|^3}r. \label{1.15}
\end{align}
Thus, the center of mass undergoes uniform motion with zero acceleration, while the relative motion reduces to that of a single particle moving under an effective central gravitational potential. Since $r_1 = R - \frac{m_2 r}{m_1 + m_2}$ and $r_2 = R + \frac{m_1 r}{m_1 + m_2}$, the full two-body problem reduces to an equivalent one-body problem whose dynamics are governed by \eqref{1.15}. In addition to the total energy,
\begin{align}
    E := \frac{1}{2} m_1 \|\dot{r}_1\|^2 + \frac{1}{2} m_2 \|\dot{r}_2\|^2 - \frac{G m_1 m_2}{\|r_2 - r_1\|} \label{1.16}
\end{align}
the system admits another fundamental invariant, namely the total angular momentum,
\begin{align}
    L := r_1 \times m_1 \dot{r}_1 + r_2 \times  m_2 \dot{r}_2. \label{1.17}
\end{align}
The conservation of both quantities follows directly from the central nature of the gravitational force. In particular, the conservation of $L$ implies that the motion of the two bodies is confined to a fixed plane. Consequently, it suffices to solve \eqref{1.15} in $\mathbb{R}^2$ rather than $\mathbb{R}^3$. The resulting planar problem admits closed-form solutions in terms of conic sections, as first established by Kepler. 

\begin{figure}
\centering
\begin{subfigure}{.5\textwidth}
  \centering
  \includegraphics[width=1\linewidth]{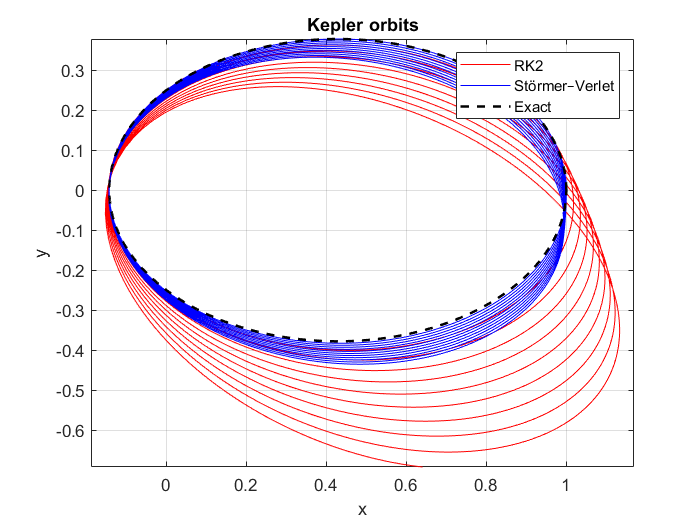}
  \caption{}
\end{subfigure}%
\begin{subfigure}{.5\textwidth}
  \centering
  \includegraphics[width=1\linewidth]{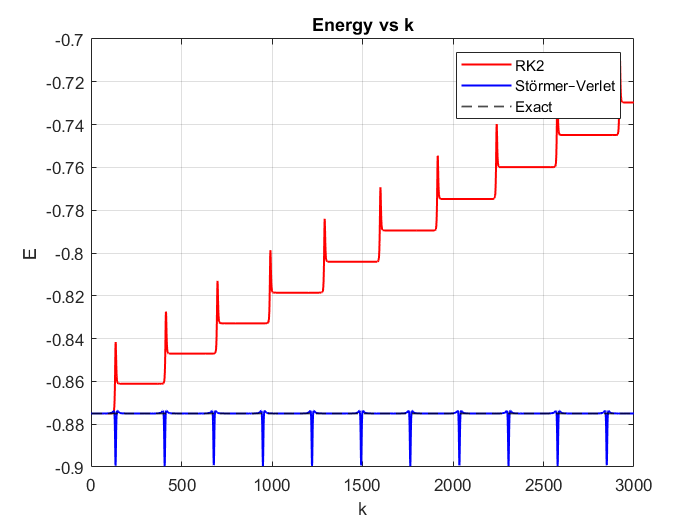}
  \caption{}
\end{subfigure}
\begin{subfigure}{.5\textwidth}
  \centering
  \includegraphics[width=1\linewidth]{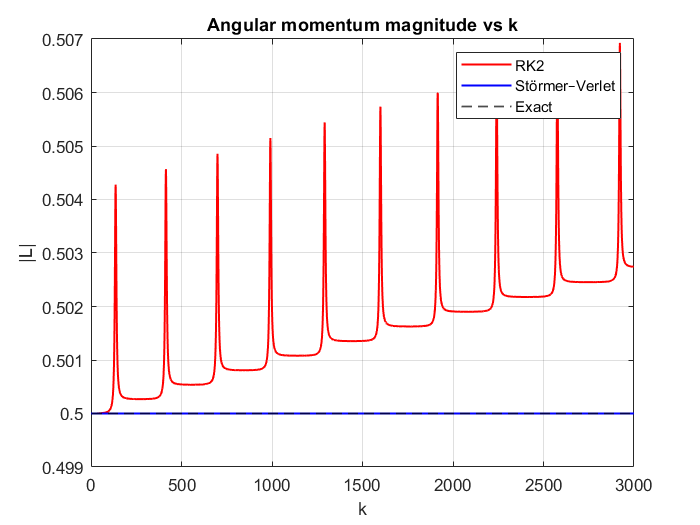}
  \caption{}
\end{subfigure}
\caption{(a) Comparison of Keplerian orbits computed using the second-order Runge–Kutta (RK2) method and the St{\"o}rmer–Verlet (SV) scheme, with parameter values (in standard units) $G(m_1+m_2)=1, h =0.01$, and initial conditions $x_0=(1,0,0,0.5)$ over 3000 iterations. The exact orbit is shown for reference. (b) Energy versus iteration for the RK2 and SV methods. (c) Magnitude of the angular momentum versus iteration for the RK2 and SV methods.}
\label{fig1.2}
\end{figure}

Let us rewrite \eqref{1.15} as a first-order dynamical system in order to make it amenable to numerical integration by introducing the state as
\begin{align}
    x := \begin{pmatrix}
        r \\
        \dot{r}
    \end{pmatrix} \implies
    \dot{x} = \begin{pmatrix}
        \dot{r} \\
        \ddot{r}
    \end{pmatrix} = \begin{pmatrix}
        \dot{r} \\
        -G(m_1+m_2)r/\|r\|^3
    \end{pmatrix} =: f(x).
\end{align}
We now present an overview of the class of numerical integration schemes that are most commonly used in physics and engineering.
\paragraph{Runge-Kutta methods} 
These methods, which date back to the beginning of the twentieth century, together with their numerous variants, constitute the workhorses of physics and engineering. For a system governed by \eqref{1.1}, an $s$-stage Runge-Kutta integrator is defined by
\begin{align*}
    x_{k+1} = x_k + h \sum_{i=1}^{s} b_ik_i \quad \text{where} \quad k_i = f(x_k + h\sum_{j=1}^{s} a_{ij}k_j). 
\end{align*}
The coefficients $a_{ij}$ and $b_i$, for $i,j=1,\dots,s$, uniquely characterize the scheme. If $a_{ij} = 0$ for $j \geq i$, the method is explicit. If, on the other hand, some coefficients $a_{ij} \neq 0$ for $j \geq i$, the method is implicit, in which case a system of nonlinear algebraic equations must be solved at each time step. The order of accuracy of the scheme depends on the choice of these coefficients. For example, the conditions:
\begin{align*}
    &\sum_{i=1}^s b_i =1 \quad \text{(order 1)}, \\
    &\sum_{i=1}^n \left( b_i \sum_{j=1}^s a_{ij} \right) = 1/2 \quad \text{(order 2)}, \\
    &\sum_{i=1}^s \left( b_i \left( \sum_{j=1}^s a_{ij} \right)^2 \right) = 1/3, \quad \sum_{i,j=1}^s \left( b_i a_{ij} \sum_{k=1}^s a_{jk} \right) =1/6 \quad \text{(order 3)}
\end{align*}
must be satisfied for first-, second-, and third-order accuracy, respectively. Analogous, but increasingly intricate, conditions arise for methods of order four and higher, rendering the analysis progressively more cumbersome. All such schemes are single-step methods. Both the explicit ($s=1,\, b_1=1,\, a_{11}=0$) and implicit ($s=1,\, b_1=1, \, a_{11}=1$) Euler methods belong to the Runge–Kutta family and represent its simplest realizations.
\paragraph{Partitioned Runge-Kutta methods}
In some situations, it is advantageous to decompose the state variables into two components,
\begin{align*}
    \dot{x} = f(x) \implies \begin{pmatrix}
        \dot{q} \\
        \dot{p}
    \end{pmatrix} = \begin{pmatrix}
        f_1(q,p) \\
        f_2(q,p)
    \end{pmatrix}
\end{align*}
where one subset of variables, say $q$, is advanced in time using a given Runge–Kutta scheme, while the remaining subset, say $p$, is evolved using a different Runge–Kutta scheme. Such a scheme is called a partitioned Runge-Kutta method and is defined by
\begin{align*}
    q_{k+1} &= q_k + h \sum_{i=1}^s b_i k_i \quad \text{and} \quad p_{k+1} = p_k + h \sum_{i=1}^s \hat{b}_i l_i \quad \text{where} \\
    k_i &= f_1(q_k + h \sum_{j=1}^s a_{ij} k_j,\, p_k + h \sum_{j=1}^s \hat{a}_{ij} l_j) \quad \text{and} \\
    l_i &= f_2(q_k + h \sum_{j=1}^s a_{ij} k_j,\, p_k + h \sum_{j=1}^s \hat{a}_{ij} l_j)
\end{align*}
where the coefficients $a_{ij},\, b_i,\, \hat{a}_{ij},\, \hat{b}_i$, for $i,j=1,\dots,s$, uniquely characterize the scheme. If the Runge-Kutta method with coefficients $(a_{ij},b_i)$ is of order $r$, and the Runge-Kutta method with coefficients $(\hat{a}_{ij},\hat{b}_i)$ is also of order $r$, then the overall order of the partitioned scheme is at most $r$.
\paragraph{Symplectic Partitioned Runge-Kutta methods}
Certain partitioned Runge–Kutta methods preserve symplecticity provided their coefficients satisfy the additional conditions
\begin{align*}
    b_i \hat{a}_{ij} + \hat{b}_j a_{ji} - b_i \hat{b}_j &= 0 \quad \text{for} \quad i,j=1,\dots,s \quad \text{and} \\
    b_i - \hat{b}_i &= 0 \quad \text{for} \quad i=1,\dots,s.
\end{align*}
Partitioned Runge–Kutta schemes satisfying these conditions are referred to as symplectic partitioned Runge–Kutta methods. The symplectic Euler method, corresponding to the choice
\begin{align*}
    s=1, \quad b_1 =1, \quad \hat{b}_1=1,\quad a_{11}=0, \quad \hat{a}_{11}=1,
\end{align*}
provides a canonical example of a first-order symplectic partitioned Runge–Kutta scheme. An example of a second-order symplectic partitioned Runge–Kutta method is obtained by taking
\begin{align*}
    s=2, \quad b = \hat{b} = \begin{pmatrix}
        1/2 \\
        1/2
    \end{pmatrix}, \quad a = \begin{pmatrix}
        0 & 0 \\
        1 & 0
    \end{pmatrix}, \quad \hat{a} = \begin{pmatrix}
        1/2 & 0 \\
        1/2 & 0
    \end{pmatrix}.
\end{align*}
This scheme is commonly know as the St{\"o}rmer-Verlet method.
\vspace{10pt} \\
Consider the plots shown in figure \ref{fig1.2}. The trajectory generated by a standard Runge–Kutta scheme spirals away from the exact elliptical orbit as energy is artificially injected into the system. In contrast, the trajectory produced by a symplectic partitioned Runge–Kutta scheme exhibits a slow drift but remains closed. In this case, the energy remains bounded and oscillates around its exact value, while the angular momentum is preserved exactly.
\subsection{Planar Pendulum}
\label{planar pendulum}
\begin{figure}
\centering
\begin{subfigure}{.5\textwidth}
  \centering
  \includegraphics[width=1\linewidth]{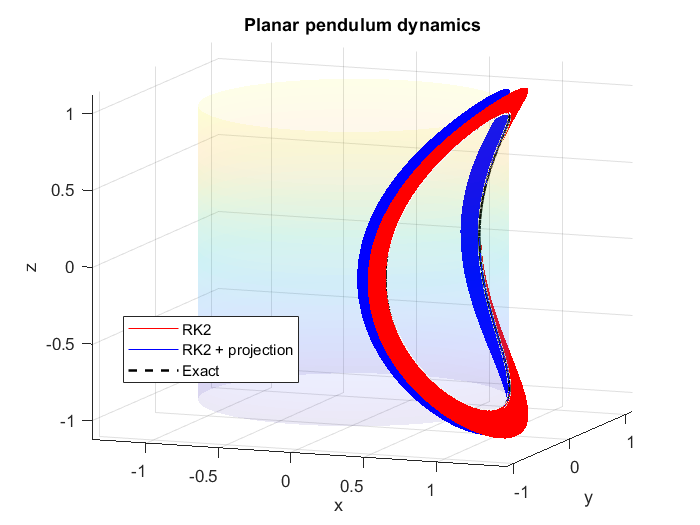}
  \caption{}
\end{subfigure}%
\begin{subfigure}{.5\textwidth}
  \centering
  \includegraphics[width=1\linewidth]{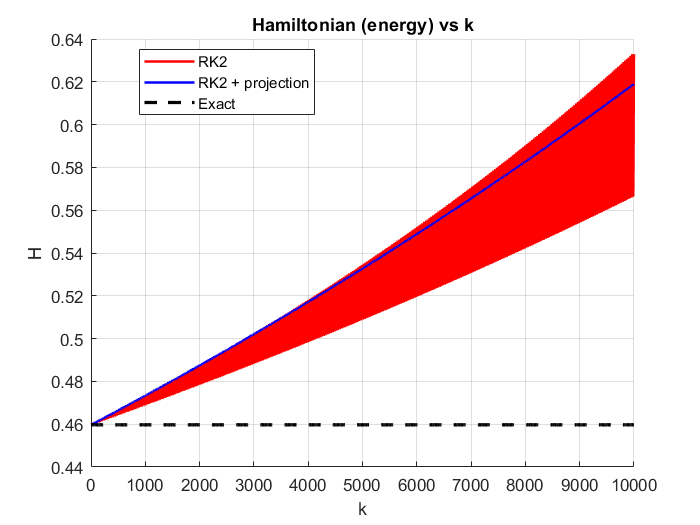}
  \caption{}
\end{subfigure}
\caption{(a) Comparison of projected and unprojected trajectories for the embedded planar pendulum dynamics integrated using the second-order Runge–Kutta (RK2) scheme, with parameter values (in standard units) $ml^2 =1, mgl =1, h= 0.1$, and initial conditions $(\theta_0,p_0)=(1,0)$, over 1000 iterations. The exact trajectory is shown for reference. (b) Hamiltonian (energy) versus iteration for the RK2 method with and without projection.}
\label{fig1.3}
\end{figure}

\begin{center}
    \begin{tikzpicture}
    \draw[] (3,4)--(5,1);
    \draw[fill] (5,1) circle (5pt);
    \draw[dotted] (3,4)--(3,1);
    \draw[-stealth] (3,3) arc (-90:-57:1);
    \node at (3.3,2.7) {$\theta$};
    \node at (4.5,2.5) {$l$};
    \node at (5,0.5) {$m$};
    \draw[-stealth] (2,3)--(2,2);
    \draw[ultra thick] (2.5,4)--(3.5,4);
    \node at (2,1.5) {$g$};
\end{tikzpicture}
\end{center}
Consider a point mass $m$ attached to a rigid, massless rod of length $l$, which is fixed at one end to a stationary support and free to rotate in a vertical plane under the influence of a uniform gravitational field, as illustrated in the figure above. The Hamiltonian (total energy) of this system is given by
\begin{align}
    H := \frac{1}{2ml^2} p^2 - mgl \cos{\theta} \label{1.19}
\end{align}
where $\theta \in \mathbb{S}^1$ (unit circle) is the angular displacement and $p \in \mathbb{R}$ is the angular momentum. The equations of motion can be written in Hamiltonian form as 
\begin{align}
    \begin{pmatrix}
        \dot{\theta} \\
        \dot{p}
    \end{pmatrix} = \begin{pmatrix}
        \partial H / \partial p \\
        -\partial H / \partial \theta
    \end{pmatrix} = \begin{pmatrix}
        p/ml^2 \\
        -mgl \sin{\theta}
    \end{pmatrix}. \label{1.20}
\end{align}
Unlike the previous examples, the dynamics governed by \eqref{1.20} evolve on a nonlinear phase space, namely the cylinder $\mathbb{S}^1 \times \mathbb{R}$. Consequently, the numerical methods introduced so far cannot be applied directly, as they rely on the linear structure of the underlying state space. One approach to circumvent this limitation is outlined below.
\paragraph{Projection-based methods} 
According to Whitney’s embedding theorem, any smooth manifold can be embedded into a sufficiently high-dimensional Euclidean space. Exploiting this fact, we embed the cylindrical phase space into $\mathbb{R}^3$, via the explicit map
\begin{align*}
    \mathbb{S}^1 \times \mathbb{R} \ni (\theta, p) \longmapsto (\cos{\theta}, \sin{\theta}, p) \in \mathbb{R}^3
\end{align*}
where the axis of the cylinder coincides with the z-axis. Using this embedding, the dynamics \eqref{1.20} can be pushed forward from the cylinder to $\mathbb{R}^3$, yielding
\begin{align}
    \begin{pmatrix}
        \dot{x} \\
        \dot{y} \\
        \dot{z}
    \end{pmatrix} = \begin{pmatrix}
        -yz/ml^2 \\
        xz/ml^2 \\
        -mgly
    \end{pmatrix} \label{1.21}
\end{align}
where $(x,y,z) \in \mathbb{R}^3$. If one naively applies a standard Runge–Kutta method to integrate \eqref{1.21}, the resulting discrete trajectory generally drifts away from the embedded cylinder
\begin{align*}
    \{ (x,y,z) \in \mathbb{R}^3 \,|\, x^2 +y^2 =1 \}.
\end{align*}
Consequently, it becomes necessary to project the numerical solution back onto the cylinder after each time step, as illustrated in figure \ref{fig1.3}.
Even though this projection based method succeeds in confining the discrete trajectory to the cylinder, it fails to respect the qualitative features of the dynamics. This deficiency is manifested in the unbounded growth of the system’s energy, as illustrated in figure \ref{fig1.3}. Just as the conservative nature of dynamics on a linear state space is rooted in a linear symplectic structure, the dynamics on a nonlinear state space are underpinned by a nonlinear symplectic structure, which we now define.
\paragraph{Nonlinear Symplectic Structure}
\label{para nonlinear symplectic structure}
Let us examine the planar pendulum from the perspective of differential geometry, see chapter \ref{ch2}. The cylindrical phase space of the pendulum is locally modeled on the symplectic plane $(\mathbb{R}^2, \omega_{std})$. Consequently, a global symplectic structure on the cylinder $\mathbb{S}^1 \times \mathbb{R}$ can be viewed as a smoothly varying family of local linear symplectic structures that are compatible on overlaps. We now make these notions mathematically precise following \cite{weinstein1971symplectic, weinstein1977lectures, da2001lectures, mcduff2017introduction}.
\begin{definition}
A \emph{symplectic manifold} is a pair $(M,\omega)$ consisting of a smooth manifold $M$ and a differential two-form $\omega \in \Omega^2(M)$ such that
\begin{itemize}
    \item $\omega$ is \emph{closed}, i.e. $d\omega = 0$, and
    \item $\omega$ is \emph{non-degenerate}, i.e. the bundle map $\omega^\flat : TM \to T^*M$,   
    is a vector bundle isomorphism.
\end{itemize}
\end{definition}
Here $\tau_M : TM \to M$ and $\pi_M : T^*M \to M$ denote the tangent and cotangent bundles of $M$, respectively, along with their bundle projections. The non-degeneracy of $\omega$ implies that at each point $x \in M$, the pair $(T_xM, \omega_x)$ is a symplectic vector space. In particular, every symplectic manifold is necessarily even-dimensional. The closedness condition ensures that $\omega$ induces a Poisson structure on $C^\infty(M)$, the algebra of smooth functions on $M$. Among the most important and ubiquitous examples of symplectic manifolds are cotangent bundles. 
\begin{definition}
    Let $Q$ be a smooth manifold with cotangent bundle $T^*Q$. The cotangent bundle carries a canonical differential one-form $\theta_Q \in \Omega^1(T^*Q)$, called the \emph{tautological} or \emph{Liouville one-form}, defined by 
    \begin{align*}
        \langle \theta_Q , Z \rangle := \langle \pi_{T^*Q} (\theta_Q), T \pi_Q (Z) \rangle 
    \end{align*}
    for all vector fields $Z \in \mathfrak{X}(T^*Q)$.
\end{definition}
\begin{center}
    \begin{tikzcd}[column sep=large, row sep=large]
        & TT^*Q \arrow[r,"T\pi_Q"] \arrow[d,shift right,"\tau_{T^*Q}",swap] & TQ \arrow[d,"\tau_Q"] \\
        T^*T^*Q \arrow[r, shift right,"\pi_{T^*Q}", swap] & T^*Q \arrow[r,"\pi_Q",swap] \arrow[u,shift right,"Z",swap] \arrow[l, shift right,"\theta_Q", swap] & Q
    \end{tikzcd}
\end{center}
Here $\langle \cdot \,, \cdot \rangle$ denotes the natural pairing between vectors and covectors. The maps, sections, and bundle projections appearing in the above definition are summarized in the commutative diagram displayed above. Every cotangent bundle $T^*Q$ is naturally equipped with a canonical symplectic two-form defined by
\begin{align*}
    \omega_Q := -d \theta_Q.
\end{align*}
Let $(q^i)$ be local coordinates on $Q$, and let $(q^i,p_i)$ denote the induced local coordinates on $T^*Q$. In these coordinates, the Liouville one-form and the canonical symplectic two-form on $T^*Q$ are $\theta_Q := p_i dq^i$ and $\omega_Q := dq^i \wedge dp_i$, respectively, where the Einstein summation convention is understood. According to Darboux's theorem, for any symplectic manifold $(M,\omega)$ and any point $x \in M$, there exist local coordinates $(q^i,p_i)$, called Darboux coordinates, defined in a neighbourhood of $x$ such that the symplectic form $\omega$ takes the canonical form $\omega = dq^i \wedge dp_i$ in these coordinates. The cylindrical phase space of the planar pendulum can be identified with the cotangent bundle  $T^*\mathbb{S}^1 \cong \mathbb{S}^1 \times \mathbb{R}$, equipped with its canonical symplectic form $\omega_{\mathbb{S}^1} = d\theta \wedge dp$. We are now ready to define Hamiltonian dynamics on symplectic manifolds.
\begin{definition}
    Let $(M, \omega)$ be a symplectic manifold and $H \in C^\infty(M)$ be a smooth function, called the \emph{Hamiltonian}. The \emph{Hamiltonian vector field} associated with $H$ is the unique vector field $X_H \in \mathfrak{X}(M)$ defined by
    \begin{equation*}
        X_H \,\lrcorner \, \omega = dH \quad \text{or equivalently,} \quad \omega(X_H, \cdot) = dH(\cdot). 
    \end{equation*}
\end{definition}
%
In Darboux coordinates $(q^i,p_i)$, the Hamiltonian vector field associated with $H$ takes the explicit form
\begin{align*}
    X_H = \frac{\partial H}{\partial p_i} \frac{\partial}{\partial q^i} - \frac{\partial H}{\partial q^i} \frac{\partial}{\partial p_i}.
\end{align*}
The integral curves of $X_H$ are governed by the celebrated Hamilton’s equations,
\begin{align*}
    \dot{q}^i = \frac{\partial H}{\partial p_i} \quad \text{and} \quad \dot{p}_i = - \frac{\partial H}{\partial q^i}.
\end{align*}
The planar pendulum dynamics given in \eqref{1.20} arise precisely in this manner from the Hamiltonian vector field associated with $H$ given in \eqref{1.19}, expressed in Darboux coordinates $(\theta,p)$ on the cotangent bundle $T^*\mathbb{S}^1$. We now define the canonical Poisson structure associated with a symplectic manifold.
\begin{definition}
    Let $(M,\omega)$ be a symplectic manifold and let $F,G \in C^\infty(M)$ be smooth functions. The bilinear operation
    $$\{F,G\} := \omega(X_F,X_G)$$
    where $X_F$ and $X_G$ are Hamiltonian vector fields associated with $F$ and $G$, respectively, defines the \emph{canonical Poisson bracket} on $(M,\omega)$.
\end{definition}
In Darboux coordinates $(q^i,p_i)$, the Poisson bracket takes the explicit form
$$\{F,G\} = \frac{\partial F}{\partial q^i} \frac{\partial G}{\partial p_i}-\frac{\partial F}{\partial p_i} \frac{\partial G}{\partial q^i}.$$
We now define the nonlinear analogue of a symplectic linear map.
\begin{definition}
    Let $(M_1, \omega_1)$ and $(M_2, \omega_2)$ be two symplectic manifolds. A diffeomorphism 
    \begin{align*}
        f : M_1 \to M_2
    \end{align*}
    is called a \emph{symplectomorphism} if it preserves the symplectic form, i.e.,
    \begin{align*}
        f^*\omega_2 = \omega_1.
    \end{align*}
    Equivalently, for all $x \in M_1$ and all $u,v \in T_x M_1$,
    \begin{align*}
        \omega_1(u,v) = \omega_2(T_xf(u),T_xf(v)).
    \end{align*}
\end{definition}
A canonical and ubiquitous source of symplectomorphisms is provided by cotangent lifts of diffeomorphisms, which preserve the canonical symplectic structures on cotangent bundles. We now define the notion of a Lagrangian submanifold, a central geometric object that will underpin our construction of symplectic integrators in subsequent sections.
\begin{definition}
    Let $(M,\omega)$ be a symplectic manifold and let $L \subset M$ be an embedded submanifold. Then $L$ is called a \emph{Lagrangian submanifold} if, for every $x \in L$,
    \begin{align*}
        (T_xL)^\omega = T_xL \quad \text{and} \quad \text{dim}\,L = \tfrac{1}{2} \,\text{dim}\,M.
    \end{align*}
\end{definition}
For example, the image of a \emph{closed} differential one-form $\alpha \in \Omega^1(Q)$ defines a Lagrangian subamnifold of the cotangent bundle $(T^*Q,\omega_Q)$. In particular, if $\alpha$ is \emph{exact}, i.e. $\alpha = df$ for some $f \in C^\infty(Q)$, then the image of $df$ is a Lagrangian submanifold. In this case, the function $f$ is called a \emph{generating function} for the Lagrangian subamnifold. One of the most important properties of Lagrangian submanifolds is that they are invariant under symplectic maps: the image of a Lagrangian submanifold under a symplectomorphism is again a Lagrangian submanifold.
\vspace{10pt}\\
We now return to the planar pendulum problem. The numerical integration of its dynamics presents two distinct challenges. The primary challenge is to ensure that the discrete trajectory remains confined to the nonlinear phase space, while the secondary challenge is to respect the intrinsic geometric structure of that phase space. Projection-based schemes can successfully address the first issue by enforcing the constraint extrinsically; however, they fail to preserve the underlying geometric structure of the dynamics. Even if one were to construct a symplectic embedding of $(T^*\mathbb{S}^1, \omega_{\mathbb{S}^1})$ into a linear symplectic space such as $(\mathbb{R}^4, \omega_{std})$ instead of $\mathbb{R}^3$, a non-trivial task in its own right, and then apply a symplectic partitioned Runge–Kutta scheme in the ambient space, the subsequent projection back onto the embedded manifold would generally fail to be a symplectomorphism. This limitation stems from their fundamentally extrinsic viewpoint. These observations motivate the construction of numerical integrators that are formulated intrinsically on the phase space itself. Before doing so, we provide a brief survey of the relevant literature, including modern developments in the field.
\section{Literature Survey}
\label{sec lit survey}
The development of numerical methods that respect the intrinsic geometry of dynamical systems has been a central theme in applied mathematics for several decades. Early work in geometric numerical integration \cite{hairer2013geometric, feng2010symplectic, sanz1993symplectic} demonstrated that standard integrators often fail to preserve qualitative features of Hamiltonian and mechanical systems over long time intervals. This insight motivated the systematic development of structure-preserving methods \cite{leimkuhler2004simulating, mclachlan2006geometric, quispel2008new, celledoni2009energy}, designed to conserve symplectic structure, momentum maps, and Lie–Poisson geometry.

A foundational tool in this direction is the notion of retraction maps \cite{absil2007optimization}, which provide computationally efficient approximations to the exponential map and enable intrinsically defined updates on manifolds as we shall see in chapter \ref{ch3}. Combined with discrete variational mechanics \cite{marsden2001discrete}, retractions led to variational integrators that are symplectic and satisfy discrete Noether theorems, ensuring momentum conservation in the presence of symmetry \cite{lee2007lie, lee2007lie_2, ma2010lie, leok2011discrete, leok2012general, bou2009hamilton, bou2018geometric}. On Lie groups, parallelizability further simplifies these constructions while preserving geometric fidelity. This line of work was unified and clarified in \cite{barbero2023retraction}, which identified retractions and discretization maps as the fundamental building blocks underlying a broad class of geometric integrators. Extensions to reduced Hamiltonian systems yielded Lie–Poisson integrators \cite{marsden1999discrete}, enabling faithful simulation of systems such as rigid bodies, heavy tops, and fluids while preserving Casimirs and coadjoint orbits.

Recent research has significantly expanded this classical framework. Stochastic variational integrators on Lie groups \cite{wu_gaybalmaz_2023_stochastic_liegroup_gsi} preserve symplectic and Lie–Poisson structure under stochastic forcing, including systems with advected quantities. Time-adaptive variational integrators \cite{duruisseaux2021adaptive, duruisseaux_leok_2023_time_adaptive_variational} extend geometric methods to multi-scale dynamics and optimization. Structure-preserving discretizations have also been embedded into learning architectures, most notably in Lie Forced Variational Integrator Networks \cite{duruisseaux_duong_leok_atanasov_2023_liefvin}, which improve accuracy and generalization by enforcing geometric constraints. Further extensions address nonsmooth and hybrid systems \cite{leok_tran_2025_collision_integrators}, higher-order and alternative variational principles \cite{sharma_patil_woolsey_2022_hermite_variational, tran_leok_2023_typeII_liegroup}, large-scale multibody dynamics \cite{brudigam_kremheller_popp_2023_variational_graph_multibody, holzinger_hofstetter_2024_liegroup_multibody}, and general Poisson systems \cite{eidnes_2022_order_theory_discrete_gradient, cosserat2024numerical}.

Against this backdrop, the present work emphasizes retraction and discretization maps as a unifying conceptual core for constructing structure-preserving integrators on Lie groups and their associated bundles. This viewpoint recovers classical variational and Lie–Poisson schemes while providing a natural entry point for recent advances in stochastic dynamics, adaptivity, and learning-based methods. Overall, the literature reflects a mature yet rapidly expanding field, where geometric structure remains the guiding principle across increasingly broad classes of dynamical systems.

\chapter{A Short Interlude on Differential Geometry}
\label{ch2}
We begin with a brief introduction to smooth manifolds and the elements of differential geometry required for our exposition. This overview is not intended to be exhaustive and primarily serves to establish notation and conventions. For a more comprehensive treatment, the reader is referred to \cite{lee2012introduction, michortopics, abraham1993manifolds, bott2013differential, mishchenko1988course, burns1985modern, sternberg1964lectures}.
\section{Smooth Manifolds}
\label{sec manifolds}
We begin by making precise what we mean by a nonlinear space, or more formally, a manifold. Intuitively, a smooth manifold may be viewed as a collection of open subsets of $\mathbb{R}^n$ that are smoothly ``stitched together'' so that the tools of calculus remain applicable on the resulting space. We now formalize this notion.
\begin{center}
\begin{tikzpicture}
\shadedraw[white] (2,1) to[out=0,in=180] (5,2) to[out=0,in=0] (3,4) to[out=180,in=90] (1,3) to[out=270,in=180] (2,1);
\node at (3.5,1) {$M$};
\fill[red,fill opacity=0.2] (2.5,2.5) ellipse (20 pt and 10 pt); 
\node at (2,2) {\footnotesize $\textcolor{red}{U_i}$};
\fill[blue,fill opacity=0.2] (3,2.75) ellipse (20 pt and 10 pt); 
\node at (3,3.3) {\footnotesize $\textcolor{blue}{U_j}$};
\fill (2.7,2.65) circle (0.8 pt); 
\node at (2.85,2.6) {\tiny $x$};
\node at (1.85,3.05) {\footnotesize $\textcolor{violet}{U_i \cap U_j}$}; 
\draw[-stealth] (6.5,0.5)--(8.5,0.5); 
\draw[-stealth] (6.5,0.5)--(6.5,2.5); 
\node at (9,0.5) {$\mathbb{R}^n$};
\fill[red,fill opacity=0.2] (7.5,1.5) circle (20 pt);
\fill[blue,fill opacity=0.2] (7.8,0.86) arc (-65:65:20pt) (7.8,2.14) arc (115:245:20 pt);
\node at (7.2,2.4) {\footnotesize $\textcolor{red}{\varphi_i(U_i)}$};
\fill (7.8,1.5) circle (0.8 pt); 
\node at (7.8,1.3) {\tiny $\varphi_i(x)$};
\node at (9.1,1.5) {\footnotesize $\textcolor{violet}{\varphi_i(U_i \cap U_j)}$};
\draw[-stealth] (6.5,3)--(8.5,3); 
\draw[-stealth] (6.5,3)--(6.5,5); 
\node at (9,3) {$\mathbb{R}^n$};
\fill[blue,fill opacity=0.2] (7.5,4) circle (20 pt);
\fill[red,fill opacity=0.2] (7.8,3.36) arc (-65:65:20 pt) (7.8,4.64) arc (115:245:20 pt);
\node at (7.2,4.9) {\footnotesize $\textcolor{blue}{\varphi_j(U_j)}$};
\fill (7.8,4) circle (0.8 pt);
\node at (7.8,3.8) {\tiny $\varphi_j(x)$};
\node at (9.1,4) {\footnotesize $\textcolor{violet}{\varphi_j(U_i \cap U_j)}$};
\draw[red,->] (3,2.2)--(6.7,1.5);
\node at (5,1.6) {\footnotesize $\textcolor{red}{\varphi_i}$};
\draw[blue,->] (3.75,2.8)--(6.7,4);
\node at (5,3.6) {\footnotesize $\textcolor{blue}{\varphi_j}$};
\draw[violet,->] (7.9,2.2)--(7.9,3.3);
\node at (8.6,2.5) {\footnotesize $\textcolor{violet}{\varphi_j \circ \varphi_i^{-1}}$};
\end{tikzpicture}
\end{center}
Consider a topological space $M$ as shown above. Let $x \in M$ and $U \subset M$ be an open neighborhood of $x$. A continuous map $\varphi : U \to \mathbb{R}^n$ that is a homeomorphism onto its image is called a (local) chart around $x \in M$. Two charts $(U_i,\varphi_i)$ and $(U_j,\varphi_j)$ are said to be compatible if either $U_i \cap U_j = \emptyset$, or, whenever $U_i \cap U_j \neq \emptyset$, the transition map $\varphi_j \circ \varphi_i^{-1} : \varphi_i(U_i \cap U_j) \to \varphi_j(U_i \cap U_j)$ is a homeomorphism between open subsets of $\mathbb{R}^n$. A collection of pairwise compatible charts whose domain cover the entire space, $\bigcup_{i\in I} U_i = M$, is called an atlas $\mathcal{A}:=\{(U_i,\varphi_i) \, | \,\, i \in I\}$ on $M$. Two atlases $\mathcal{A}_\alpha$ and $\mathcal{A}_\beta$ are said to be compatible if their union $\mathcal{A}_\alpha \cup \mathcal{A}_\beta$ is itself an atlas. The collection of all atlases compatible with a given atlas defines a unique maximal atlas $\mathcal{A}_{max}$ on $M$.  
\begin{definition} \label{def manifold}
    A \emph{smooth manifold} is a topological space $M$ equipped with a maximal atlas $\mathcal{A}_{max}$ whose transition maps are diffeomorphisms. 
\end{definition}
Here, $n$ denotes the (topological) dimension of the manifold $M$. Standard examples of smooth manifolds are Euclidean spaces $\mathbb{R}^n$, the $n$-sphere $\mathbb{S}^n$, and the $n$-torus $\mathbb{T}^n := \mathbb{S}^1 \times \cdots \times \mathbb{S}^1$ ($n$-copies). The configuration spaces of mechanical systems provide a rich class of further examples of smooth manifolds. For instance, the configuration space of a planar pendulum is the circle $\mathbb{S}^1$, that of a planar double pendulum is the torus $\mathbb{T}^2$, and that of a planar pendulum mounted on a linear cart is the product manifold $\mathbb{R} \times \mathbb{S}^1$.
\section{Tangent and Cotangent Spaces}
\label{sec tangent cotangent spaces}
We now turn to one of the central ideas of differential geometry: the linearization of a manifold at a point, which leads to the notion of the tangent space. There are several equivalent ways to define the tangent space of a manifold; here we adopt a definition motivated by physical intuition, viewing tangent vectors as velocities of smooth curves passing through a point as shown below.
\begin{center}
\begin{tikzpicture}
\shadedraw[white] (3,1) to[out=0,in=180] (6,2) to[out=0,in=0] (4,4) to[out=180,in=90] (2,3) to[out=270,in=180] (3,1);
\node at (4.5,1) {$M$};
\draw[|-|,red] (1,1.7)--(1,3.3); 
\fill[red] (1,2.5) circle (0.8 pt);
\node at (0.8,2.5) {\footnotesize $\textcolor{red}{0}$};
\draw[red] (2.6,2.2) to[out=50,in=170] (5,3.2); 
\fill[red] (4,3.17) circle (0.8 pt);
\node at (4.3,3) {\tiny $\textcolor{red}{x=c(0)}$};
\draw[blue,dotted] (4,3.17)--(4,4.17);
\fill[blue] (4,4.17) circle (0.8 pt);
\draw[blue,fill=blue,fill opacity=0.2] (3.5,3.5)--(5.7,4)--(4.5,4.8)--(2.3,4.3)--cycle; 
\draw[thick,-stealth,blue] (4,4.17)--(4.5,4.33);
\node at (5.1,4.5) {\tiny $\textcolor{blue}{v=c'(0)=[c]_x}$};
\node at (3,4.8) {$\textcolor{blue}{T_xM}$};
\fill[orange,fill opacity=0.2] (4,3) ellipse (25 pt and 15 pt);
\node at (4,2.3) {\footnotesize $\textcolor{orange}{U}$};
\draw[-stealth] (7.5,1.5)--(10,1.5); 
\draw[-stealth] (7.5,1.5)--(7.5,4); 
\node at (10.5,1.5) {$\mathbb{R}^n$};
\node at (8.2,2) {\footnotesize $\textcolor{orange}{\varphi(U)}$}; 
\draw[red] (8.3,2.5) to[out=40,in=190] (9.8,3.3);
\node at (9.7,2.8) {\tiny $\textcolor{red}{\varphi(x)=(\varphi \circ c)(0)}$};
\fill[orange,fill opacity=0.2] (9,3) circle (25 pt);
\draw[blue,-stealth] (9,3.02)--(10.8,3.02);
\draw[blue,-stealth] (9,3.02)--(9,4.82);
\draw[blue,thick,-stealth] (9,3.02)--(9.5,3.35);
\node at (9.7,3.6) {\tiny $\textcolor{blue}{(\varphi \circ c)'(0)}$};
\node at (9.4,5.1) {$\textcolor{blue}{T_{\varphi(x)}\mathbb{R}^n \cong \mathbb{R}^n}$};
\fill[red] (9,3.02) circle (0.8 pt);
\draw[->,red] (1.2,2.3)--(2.5,2.3);
\node at (1.7,2.5) {\footnotesize $\textcolor{red}{c}$};
\draw[orange,->] (5,2.9)--(8,2.9);
\node at (6.5,3.1) {\footnotesize $\textcolor{orange}{\varphi}$};
\end{tikzpicture}
\end{center}
Let $c : \,]-\epsilon,\epsilon[\, \to M$ be a smooth curve on a manifold $M$ with $c(0)=x$. From an extrinsic viewpoint - thinking of $M$ as a hypersurface embedded in a sufficiently large Euclidean space - our physical intuition tells us that the velocity $c'(0)$ lies in a hyperplane tangent to the hypersurface at $x$. This velocity is what we naturally interpret as a tangent vector, and the collection of all such velocities forms the tangent space at $x$, denoted $T_xM$.

In mechanics, however, it is essential to work intrinsically, without reference to any embedding, since the dynamics depend on configurations, velocities, and momenta defined solely on the manifold itself. To formalize the notion of velocity intrinsically, two curves passing through $x$ are said to be equivalent if they have the same local velocity $(\varphi \circ c)'(0)$ in some (and hence any) chart $(U,\varphi)$ containing $x$. An equivalence class $[c]_x$ under this relation is called a tangent vector at $x$.
Thus, intrinsically, a tangent vector represents all curves passing through $x$ with the same instantaneous velocity. The set of all such equivalence classes defines the tangent space $T_xM$, which is independent of the choice of chart and naturally carries a linear structure. As in elementary calculus, this space provides a first-order, linear approximation of the nonlinear manifold near $x$.
\begin{center}
\begin{tikzpicture}
\shadedraw[white] (3,1) to[out=0,in=180] (6,2) to[out=0,in=0] (4,4) to[out=180,in=90] (2,3) to[out=270,in=180] (3,1);
\node at (5,1.5) {$M$};
\fill[] (4,3.17) circle (0.8 pt);
\node at (4,2.9) {$x$};
\draw[dotted] (4,3.17)--(4,5.17);
\fill[red] (4,4.17) circle (0.8 pt);
\draw[red,fill=red,fill opacity=0.2] (3.5,3.5)--(5.7,4)--(4.5,4.8)--(2.3,4.3)--cycle; 
\draw[thick,-stealth,red] (4,4.17)--(4.5,4.33);
\node at (4.6,4.1) {$\textcolor{red}{v}$};
\node at (2,4) {$\textcolor{red}{T_xM}$};
\fill[blue] (4,5.17) circle (0.8 pt);
\draw[blue,fill=blue,fill opacity=0.2] (3.5,4.5)--(5.7,5)--(4.5,5.8)--(2.3,5.3)--cycle;
\draw[thick,-stealth,blue] (4,5.17)--(4.5,5.33);
\node at (4.6,5.1) {$\textcolor{blue}{p}$};
\node at (2,5) {$\textcolor{blue}{T_x^*M}$};
\end{tikzpicture}    
\end{center}
Since the tangent space $T_xM$ is a vector space, it admits a dual space, called the cotangent space $T_x^*M$ as shown above. From a physical perspective, tangent vectors correspond to velocities, while cotangent vectors naturally represent momenta, acting as linear functionals on velocities. The natural pairing $\langle p,v \rangle$ between a tangent vector $v$ and a cotangent vector $p$ can be interpreted as an energy-like quantity. At each point $x \in M$, both the tangent and cotangent spaces have the same dimension as the underlying manifold.
\section{Tangent and Cotangent Bundles}
\label{sec tangent cotangent bundles}
\begin{center}
    \begin{tikzpicture}
        \draw (3,1)--(3,3);
        \node at (3,0.7) {$T_xM$};
        \draw (2,1)--(2,3);
        \draw (6,1)--(6,3);
        \draw[red,thick,-stealth] (3,2)--(3,2.4);
        \node at (3.2,2.4) {$\textcolor{red}{v}$};
        \draw (4,1)--(4,3);
        \draw (5,1)--(5,3);
        \draw[dashed] (1.5,2)--(6.5,2);
        \draw[fill] (3,2) circle (1pt);
        \node at (3.2,1.8) {$x$};
        \node at (1,2) {$M$};
        \draw (7,1) to[out=10,in=190] (7.2,2) to[out=170,in=-10] (7,3);
        \node at (7.7,2) {$TM$};
    \end{tikzpicture}
\end{center}
By collecting the tangent spaces at all points of a smooth manifold as illustrated above, we obtain what is known as the tangent bundle.
\begin{definition}\label{def tangent bundle}
    Let $M$ be a smooth manifold. The disjoint union of its tangent spaces, 
    \begin{align*}
        TM := \bigsqcup_{x \in M} T_xM
    \end{align*}
    is called the \emph{tangent bundle} of $M$.
\end{definition}
The tangent bundle is a canonical example of a vector bundle, with the bundle projection $\tau_M : TM \to M$ which maps a tangent vector to its base point by forgetting the velocity component. Moreover, $TM$ is itself a smooth manifold of twice the dimension of $M$. From a physical perspective, the tangent bundle represents the space of all configurations and velocities and serves as the natural geometric setting for Lagrangian mechanics, as will be discussed later. 
\begin{center}
    \begin{tikzpicture}
        \draw (3,1)--(3,3);
        \node at (3,0.7) {$T_x^*M$};
        \draw (2,1)--(2,3);
        \draw (6,1)--(6,3);
        \draw[blue,thick,-stealth] (3,2)--(3,2.4);
        \node at (3.2,2.4) {$\textcolor{blue}{p}$};
        \draw (4,1)--(4,3);
        \draw (5,1)--(5,3);
        \draw[dashed] (1.5,2)--(6.5,2);
        \draw[fill] (3,2) circle (1pt);
        \node at (3.2,1.8) {$x$};
        \node at (1,2) {$M$};
        \draw (7,1) to[out=10,in=190] (7.2,2) to[out=170,in=-10] (7,3);
        \node at (7.7,2) {$T^*M$};
    \end{tikzpicture}
\end{center}
Similarly, by collecting the cotangent spaces at all points of the manifold as illustrated above, we obtain the cotangent bundle.
\begin{definition}\label{def cotangent bundle}
    Let $M$ be a smooth manifold. The disjoint union of its cotangent spaces, 
    \begin{align*}
        T^*M := \bigsqcup_{x \in M} T_x^*M
    \end{align*}
    is called the \emph{cotangent bundle} of $M$.
\end{definition}
The cotangent bundle is also a vector bundle, with the bundle projection $\pi_M : T^*M \to M$ which maps a cotangent vector to its base point by forgetting the momentum component. As in the tangent bundle case, $T^*M$ is a smooth manifold of twice the dimesnion of $M$. Physically, it represents the space of all configurations and momenta and provides the natural phase space for Hamiltonian mechanics as expounded in \ref{para nonlinear symplectic structure}.

Any smooth map $f : M \to N$ between smooth manifolds naturally induces a map between their tangent bundles, called the tangent lift, $Tf : TM \to TN$, as illustrated in the figure below.
\begin{center}
\begin{tikzpicture}
\shadedraw[white] (3,1) to[out=0,in=180] (6,2) to[out=0,in=0] (4,4) to[out=180,in=90] (2,3) to[out=270,in=180] (3,1);
\node at (5,1.3) {$M$};
\draw[|-|,red] (1,1.7)--(1,3.3); 
\fill[red] (1,2.5) circle (0.8 pt);
\node at (0.8,2.5) {\footnotesize $\textcolor{red}{0}$};
\draw[red] (2.6,2.2) to[out=50,in=170] (5,3.2); 
\fill[red] (4,3.17) circle (0.8 pt);
\node at (4.3,3) {\tiny $\textcolor{red}{x=c(0)}$};
\draw[blue,dotted] (4,3.17)--(4,4.17);
\fill[blue] (4,4.17) circle (0.8 pt);
\draw[blue,fill=blue,fill opacity=0.2] (3.5,3.5)--(5.7,4)--(4.5,4.8)--(2.3,4.3)--cycle; 
\draw[thick,-stealth,blue] (4,4.17)--(4.5,4.33);
\node at (4.8,4.5) {\tiny $\textcolor{blue}{v=c'(0)}$};
\node at (3,4.8) {$\textcolor{blue}{T_xM}$};
\draw[->,red] (1.2,2.3)--(2.5,2.3);
\node at (1.7,2.5) {\footnotesize $\textcolor{red}{c}$};
\shadedraw[white] (8,1) to[out=0,in=180] (11,2) to[out=0,in=0] (9,4) to[out=180,in=90] (7,3) to[out=270,in=180] (8,1);
\node at (10,1.3) {$N$};
\draw[red] (7.6,2.2) to[out=50,in=170] (10,3.2); 
\fill[red] (9,3.17) circle (0.8 pt);
\node at (9.7,3) {\tiny $\textcolor{red}{f(x)=(f \circ c)(0)}$};
\draw[blue,dotted] (9,3.17)--(9,4.17);
\fill[blue] (9,4.17) circle (0.8 pt);
\draw[blue,fill=blue,fill opacity=0.2] (8.5,3.5)--(10.7,4)--(9.5,4.8)--(7.3,4.3)--cycle; 
\draw[thick,-stealth,blue] (9,4.17)--(9.5,4.33);
\node at (10.3,4.5) {\tiny $\textcolor{blue}{T_xf(v)=(f \circ c)'(0)}$};
\node at (8,4.8) {$\textcolor{blue}{T_{f(x)}N}$};
\draw[->] (5.8,3.2) to[out=10,in=170] (6.9,3.2);
\node at (6.4,3.6) {$f$};
\draw[->,blue] (5.6,4.2) to[out=10,in=170] (7.2,4.2);
\node at (6.4,4.6) {$\textcolor{blue}{T_xf}$};
\end{tikzpicture}
\end{center}
Concretely, if $c: \,]-\epsilon,\epsilon[\, \to M$ is a smooth curve with $c(0)=x$, then the composition $f \circ c$ is a smooth curve in $N$ passing through $f(x)$. the velocity of this curve at $t=0$ defines a tangent vector at $f(x)$, and this construction yields a linear map $T_xf : T_xM \to T_{f(x)}N$, sending the tangent vector $v=c'(0)$ to $(f\circ c)'(0)$. This map is the intrinsic, coordinate-free generalization of the Jacobian matrix from multivariate calculus.

Since $T_xf$ is a linear map between vector spaces, it admits a natural dual (or transpose) map between the corresponding cotangent spaces, $T_x^*f : T_{f(x)}^*N \to T_x^*M$, defined via the natural pairing
\begin{equation*}
    \langle T_x^*f(\alpha), v \rangle := \langle \alpha, T_xf(v) \rangle
\end{equation*}
for all $v \in T_xM$ and $\alpha \in T_{f(x)}^*N$.
In physical terms, while $T_xf$ describes how velocities are pushed forward under the map $f$, the dual map $T_x^*f$ pulls back covectors, such as momenta or forces. By assembling these maps pointwise, the smooth map $f : M \to N$ can be lifted to a smooth map between cotangent bundles, $T^*f : T^*N \to T^*M$, called the cotangent lift, which plays a fundamental role in mechanics, symplectic geometry, and variational formulations of dynamics.
\section{Vector Fields, Flows and Integral Curves}
\label{sec vector fields}
We now arrive at one of the most fundamental notions from the perspective of dynamical systems, namely that of a vector field. Intuitively, a vector field assigns a tangent vector to each point of a manifold in a smooth manner. This idea can be made precise as follows.
\begin{definition} \label{def vector field}
    Let $M$ be a smooth manifold. A smooth map $X : M \to TM$ satisfying $\tau_M \circ X = id_M$ is called a smooth vector field on $M$.
\end{definition}
In the language of vector bundles, a vector field is simply a smooth section of the tangent bundle as highlighted in the figure below. 
\begin{center}
    \begin{tikzpicture}
        \draw (3,1)--(3,3);
        \node at (3,0.7) {$T_xM$};
        \draw (2,1)--(2,3);
        \draw[red,thick,-stealth] (2,2)--(2,2.45);
        \draw (6,1)--(6,3);
        \draw[red,thick,-stealth] (6,2)--(6,1.45);
        \draw[red,thick,-stealth] (3,2)--(3,2.55);
        \node at (3.5,2.7) {$\textcolor{red}{X(x)}$};
        \draw (4,1)--(4,3);
        \draw[red,thick,-stealth] (4,2)--(4,2.3);
        \draw[red] (1.5,2.2) to[out=30,in=150] (4.5,2) to[out=-30,in=190] (6.5,1.5) ;
        \draw (5,1)--(5,3);
        \draw[red,thick,-stealth] (5,2)--(5,1.7);
        \node at (1.5,2.6) {$\textcolor{red}{X}$};
        \draw[dashed] (1.5,2)--(6.5,2);
        \draw[fill] (3,2) circle (1pt);
        \node at (3.2,1.8) {$x$};
        \node at (1,2) {$M$};
        \draw (7,1) to[out=10,in=190] (7.2,2) to[out=170,in=-10] (7,3);
        \node at (7.7,2) {$TM$};
    \end{tikzpicture}
\end{center}
The space of all vector fields on $M$ is denoted by $\mathfrak{X}(M) := \Gamma^\infty(TM)$. A closely related notion to that of a vector field is the concept of a flow as shown below.
\begin{center}
\begin{tikzpicture}
\shadedraw[white] (3,1) to[out=0,in=180] (6,2) to[out=0,in=0] (4,4) to[out=180,in=90] (2,3) to[out=270,in=180] (3,1);
\node at (5,1.5) {$M$};
\draw[orange,dashed,->] (3,1.5) to[out=40,in=180] (5.5,2.5);
\draw[orange,dashed,->] (2.5,2.2) to[out=42.5,in=177.5] (5.3,3.2);
\draw[orange,dashed,->] (2.3,3) to[out=45,in=175] (4.6,3.7);
\draw[red,thick,-stealth] (3.2,1.65)--(3.5,1.95);
\draw[red,thick,-stealth] (4.1,2.22)--(4.4,2.4);
\draw[red,thick,-stealth] (5,2.47)--(5.3,2.51);
\draw[red,thick,-stealth] (2.6,2.3)--(3,2.65);
\draw[red,thick,-stealth] (3.3,2.8)--(3.7,3.05);
\draw[red,thick,-stealth] (4,3.08)--(4.3,3.2);
\draw[red,thick,-stealth] (4.7,3.2)--(5,3.25);
\draw[red,thick,-stealth] (2.8,3.4)--(3.2,3.65);
\draw[red,thick,-stealth] (3.8,3.7)--(4.2,3.8);
\fill[orange] (2.6,2.3) circle (1pt);
\node at (2.6,2.1) {\tiny $\textcolor{orange}{x}$};
\fill[orange] (3.3,2.8) circle (1pt);
\node at (3.6,2.6) {\tiny $\textcolor{orange}{\phi_t^X(x)}$};
\node at (2.8,2.8) {\tiny $\textcolor{red}{X(x)}$};
\end{tikzpicture}
\end{center}
\begin{definition} \label{def flow}
    Let $M$ be a smooth manifold and let $X \in \mathfrak{X}(M)$. A smooth map $\phi^X : \,]-\epsilon, \epsilon[\, \times M \to M$ satisfying 
    \begin{align*}
        &\frac{d}{dt} \phi^X_t(x) = X(\phi^X_t(x)), \quad \text{for all} \quad x \in M, \\
        &\phi^X_0 = id_M \quad \text{and} \quad \phi^X_s \circ \phi^X_t = \phi^X_{s+t}
    \end{align*}
    for all $s,t,s+t \in \,]-\epsilon,\epsilon[\,$ is called the \emph{local flow} of the vector field $X$.
\end{definition}
Such a flow always exists for some $\epsilon > 0$ by standard existence and uniqueness results for ordinary differential equations. If the flow exists for all $t \in \mathbb{R}$, it is called global, and the corresponding vector field is said to be complete. In this case, the flow defines a smooth $\mathbb{R}$-action on $M$, with the vector field $X$ serving as its infinitesimal generator. Intuitively, the flow describes the time evolution of a dynamical system governed by the vector field. Fixing an initial condition yields a distinguished curve on the manifold, called an integral curve.
\begin{definition} \label{def integral curve}
    Let $M$ be a smooth manifold and let $X \in \mathfrak{X}(M)$. For a given point $x_0 \in M$, the curve $\phi^X_t(x_0)$ is called the \emph{integral curve} of $X$ through $x_0$.
\end{definition}
Thus, solving an initial value problem for a dynamical system on a manifold amounts to finding an integral curve of the associated vector field. The flow further provides a natural notion of directional differentiation on manifolds, known as the Lie derivative
\section{Covector Fields, Differential Forms and Exterior Derivatives}
\label{sec covector fields}
Just as a vector field assigns a tangent vector to each point of a manifold, a covector field assigns a cotangent vector to each point. 
\begin{definition}
    Let $M$ be a smooth manifold. A smooth map $\alpha : M \to T^*M$ satisfying $\pi_M \circ \alpha = id_M$ is called a smooth covector field on $M$.
\end{definition}
Hence, in the language of vector bundles, a covector field is a smooth section of the cotangent bundle as highlighted in the figure below.
\begin{center}
    \begin{tikzpicture}
        \draw (3,1)--(3,3);
        \node at (3,0.7) {$T_x^*M$};
        \draw (2,1)--(2,3);
        \draw[blue,-stealth,thick] (2,2)--(2,2.45);
        \draw (6,1)--(6,3);
        \draw[blue,thick,-stealth] (6,2)--(6,1.45);
        \draw[blue,thick,-stealth] (3,2)--(3,2.55);
        \node at (3.5,2.7) {$\textcolor{blue}{\alpha(x)}$};
        \draw (4,1)--(4,3);
        \draw[blue,thick,-stealth] (4,2)--(4,2.3);
        \draw[blue] (1.5,2.2) to[out=30,in=150] (4.5,2) to[out=-30,in=190] (6.5,1.5) ;
        \draw (5,1)--(5,3);
        \draw[blue,thick,-stealth] (5,2)--(5,1.7);
        \node at (1.5,2.6) {$\textcolor{blue}{\alpha}$};
        \draw[dashed] (1.5,2)--(6.5,2);
        \draw[fill] (3,2) circle (1pt);
        \node at (3.2,1.8) {$x$};
        \node at (1,2) {$M$};
        \draw (7,1) to[out=10,in=190] (7.2,2) to[out=170,in=-10] (7,3);
        \node at (7.7,2) {$T^*M$};
    \end{tikzpicture}
\end{center}
Covector fields are more commonly referred to as differential one-forms denoted by $\Omega^1(M) := \Gamma^\infty(T^*M)$. More generally, differential forms play a central role in geometry, analysis, and physics, particularly in the context of integration on manifolds. 

Differential one-forms naturally pair with tangent vectors and measure their oriented magnitude in a given direction, serving as the integrands of line integrals. Likewise, differential two-forms take pairs of vectors and return the oriented area of the parallelogram they span, encoding flux and surface integrals. This construction extends to differential $k$-forms, denoted by $\Omega^k(M)$, which act on $k$ vectors and represent oriented $k$-dimensional volume elements for all $k \leq \text{dim}(M)$. Since oriented volume elements change sign under the interchange of arguments, differential forms are naturally alternating multilinear objects. Consequently, they are equipped with the structure of an exterior algebra and can be described intrinsically as smooth sections of the exterior powers of the cotangent bundle, $\Omega^k(M) := \Gamma^\infty(\bigwedge^k T^*M)$. 

Recall gradient, curl and divergence operators from multivariate calculus. On a smooth manifold, these familiar operators are unified by a single geometric construction: the exterior derivative, denoted
\begin{equation*}
    d : \Omega^k(M) \to \Omega^{k+1}(M).
\end{equation*}
The exterior derivative is a natural differential operator that satisfies the expected properties of linearity and a Leibniz rule with respect to the wedge product, and crucially obeys $d \circ d = 0$. In Euclidean space, the gradient, curl, and divergence arise as special cases of $d$ acting on 0-, 1-, and 2- forms, respectively.

The operator $d$ gives rise to two important classes of differential forms. A $k$-form $\omega \in \Omega^k(M)$ is called closed if $d\omega = 0$, and exact if there exists a $(k-1)$-form $\alpha$ such that $\omega = d \alpha$. These notions play a central role in geometry and physics, encoding conservation laws, potentials, and topological information about the underlying manifold.
\section{Lie Groups and their Actions}
\label{sec Lie groups}
Finally, we conclude this preliminary discussion with a brief overview of Lie groups and their actions. A comprehensive reference for Lie groups is \cite{helgason1979differential}.
\begin{definition}
    A smooth manifold $G$ endowed with a group structure for which the group multiplication $G \times G \to G$ and inversion $G \to G$ maps are smooth is called a Lie group.
\end{definition}
Lie groups thus unify algebraic and geometric structures: they are manifolds whose points represent configurations, while the group operation encodes symmetries and composition laws. Typical examples of abelian (commutative) Lie groups include $\mathbb{R}^n$ and $\mathbb{T}^n$. In contrast, matrix Lie groups such as $\mathbb{GL}(n,\mathbb{R}), \mathbb{SL}(n,\mathbb{R}), \mathbb{O}(n,\mathbb{R})$ etc. provide fundamental examples of nonabelian Lie groups and arise naturally as configuration spaces of rigid bodies and mechanical systems. 

A distinguished role is played by the tangent space at the identity element $e \in G$, called the Lie algbera of $G$ and denoted by $\mathfrak{g}$. This vector space captures the infinitesimal structure of the Lie group: its elements may be interpreted as infinitesimal generators of one-parameter subgroups, corresponding physically to velocities or angular velocities. The algebraic structure on $\mathfrak{g}$, encoded by the Lie bracket, reflects the local noncommutativity of the group.
\begin{center}
\begin{tikzpicture}
\shadedraw[white] (3,1) to[out=0,in=180] (6,2) to[out=0,in=0] (4,4) to[out=180,in=90] (2,3) to[out=270,in=180] (3,1);
\node at (5,1.5) {$G$};
\fill[] (4,3.17) circle (0.8 pt);
\node at (4,2.9) {$e$};
\draw[dotted] (4,3.17)--(4,5.17);
\fill[red] (4,4.17) circle (0.8 pt);
\draw[red,fill=red,fill opacity=0.2] (3.5,3.5)--(5.7,4)--(4.5,4.8)--(2.3,4.3)--cycle; 
\node at (6.7,4) {$\textcolor{red}{T_eG = \mathfrak{g}}$};
\fill[blue] (4,5.17) circle (0.8 pt);
\draw[blue,fill=blue,fill opacity=0.2] (3.5,4.5)--(5.7,5)--(4.5,5.8)--(2.3,5.3)--cycle;
\node at (6.7,5) {$\textcolor{blue}{T_e^*G=\mathfrak{g}^*}$};
\end{tikzpicture}
\end{center}
The dual space $\mathfrak{g}^*$ is equally important, particularly in mechanics. It naturally carries a Poisson structure, making it the canonical phase space for reduced Hamiltonian dynamics and Lie–Poisson systems.
\begin{definition}
    Let $\mathfrak{g}$ be a lie algebra and let $\mathfrak{g}^*$ denote its dual. For $F,G \in C^\infty(\mathfrak{g}^*)$ and $\mu \in \mathfrak{g}^*$, define
    $$\{F,G\}_{\pm}(\mu) := \pm \langle \mu, [dF_\mu,dG_\mu] \rangle.$$ 
    Here the differentials are identified via the natural isomorphisms $dF_\mu$, $dG_\mu \in T_\mu^*\mathfrak{g}^* \cong (\mathfrak{g}^*)^* \cong \mathfrak{g}$. The bilinear operation $\{\cdot\,,\cdot\}_\pm$ defines a \emph{Lie-Poisson bracket} on $\mathfrak{g}^*$.
\end{definition}
This duality between $\mathfrak{g}$ and $\mathfrak{g}^*$ underlies much of geometric mechanics and will be central to the developments that follow.

The reason groups naturally encode symmetries of a space is that they can act on it; that is, elements of the group move points of the space in a manner compatible with the group structure. We now make this precise.
\begin{definition} \label{def group action}
    Let $M$ be a smooth manifold and $G$ a Lie group. A smooth map
    \begin{equation*}
        \Phi : G \times M \to M
    \end{equation*}
    is called a (left) action of $G$ on $M$ if it satisfies
    \begin{itemize}
        \item $\Phi(e,x)=x$, equivalently $\Phi_e = id_M$ and
        \item $\Phi(g,\Phi(h,x)) = \Phi(gh,x)$, equivalently $\Phi_g \circ \Phi_h = \Phi_{gh}$
    \end{itemize}
    for all $g,h \in G$ and $x \in M$.
\end{definition}
The set of all points obtained by acting on $x \in M$ with every element of the group is called the orbit of $x$ and is denoted by
\begin{equation*}
    G \cdot x := \{\Phi_g(x) \, | \, g \in G\}.
\end{equation*}
Just as Lie algebras encode the infinitesimal structure of Lie groups, the infinitesimal effect of a group action on a manifold is captured by its infinitesimal generators, defined via one–parameter subgroups.
\begin{definition}
    Let $\Phi : G \times M \to M$ be a smooth group action. For each $\xi \in \mathfrak{g}$, the assignment
    \begin{equation*}
        \xi_M(x) := \left. \frac{d}{dt} \right|_{t=0} \Phi_{\exp{(t\xi)}}(x) 
    \end{equation*}
    for all $\xi \in \mathfrak{g}$ defines a vector field on $M$, called the infinitesimal generator associated with $\xi$.
\end{definition}
\begin{center}
\begin{tikzpicture}
\shadedraw[white] (3,1) to[out=0,in=180] (6,2) to[out=0,in=0] (4,4) to[out=180,in=90] (2,3) to[out=270,in=180] (3,1);
\node at (5,1.5) {$G$};
\draw[orange] (2.6,2.2) to[out=50,in=170] (5,3.2); 
\fill[orange] (4,3.17) circle (0.8 pt);
\node at (4.1,3) {\footnotesize $\textcolor{orange}{e}$};
\node at (2.8,2) {\footnotesize $\textcolor{orange}{\exp{(t\xi)}}$};
\draw[red,dotted] (4,3.17)--(4,4.17);
\fill[red] (4,4.17) circle (0.8 pt);
\draw[red,fill=red,fill opacity=0.2] (3.5,3.5)--(5.7,4)--(4.5,4.8)--(2.3,4.3)--cycle; 
\draw[thick,-stealth,red] (4,4.17)--(4.5,4.33);
\node at (4.7,4.2) {\footnotesize $\textcolor{red}{\xi}$};
\node at (3,4.8) {$\textcolor{red}{\mathfrak{g}}$};
\shadedraw[white] (8,1) to[out=0,in=180] (11,2) to[out=0,in=0] (9,4) to[out=180,in=90] (7,3) to[out=270,in=180] (8,1);
\node at (10,1.5) {$M$};
\draw[orange,dashed,->] (8,1.5) to[out=40,in=180] (10.5,2.5);
\draw[orange,dashed,->] (7.5,2.2) to[out=42.5,in=177.5] (10.3,3.2);
\draw[orange,dashed,->] (7.3,3) to[out=45,in=175] (9.6,3.7);
\draw[red,thick,-stealth] (8.2,1.65)--(8.5,1.95);
\draw[red,thick,-stealth] (9.1,2.22)--(9.4,2.4);
\draw[red,thick,-stealth] (10,2.47)--(10.3,2.51);
\draw[red,thick,-stealth] (7.6,2.3)--(8,2.65);
\draw[red,thick,-stealth] (8.3,2.8)--(8.7,3.05);
\draw[red,thick,-stealth] (9,3.08)--(9.3,3.2);
\draw[red,thick,-stealth] (9.7,3.2)--(10,3.25);
\draw[red,thick,-stealth] (7.8,3.4)--(8.2,3.65);
\draw[red,thick,-stealth] (8.8,3.7)--(9.2,3.8);
\fill[orange] (7.6,2.3) circle (1pt);
\node at (7.6,2.1) {\tiny $\textcolor{orange}{x}$};
\fill[orange] (8.3,2.8) circle (1pt);
\node at (9,2.6) {\tiny $\textcolor{orange}{\Phi_{\exp{(t\xi)}}(x)}$};
\node at (7.8,2.8) {\tiny $\textcolor{red}{\xi_M(x)}$};
\node at (10.3,3) {\tiny $\textcolor{orange}{G \cdot x}$};
\draw[<-] (6.8,2.9) arc (30:330:8 pt);
\node at (6.5,3.4) {$\Phi$};
\end{tikzpicture}
\end{center}
By construction, the vector field $\xi_M$ is everywhere tangent to the group orbits $G \cdot x$ as shown in the figure above. 

Since Lie groups are themselves smooth manifolds, they admit natural actions on themselves. We now list several commonly encountered examples of such actions. 
\paragraph{Left Translation, Right Translation and Inner Automorphism}
The map $G \ni h \longmapsto L_g(h) := gh \in G$ is called left translation by $g$, while $G \ni h \longmapsto R_g(h) := hg \in G$ is called the right translation by $g$. Combining left translation by $g$ with right translation by $g^{-1}$ yields the map $G \ni h \longmapsto I_g(h) := ghg^{-1} \in G$ called an inner automorphism or conjugation by $g$.
\paragraph{Adjoint and Coadjoint Action}
The differential of the inner automorphism at the identity, $T_eI_g : \mathfrak{g} \to \mathfrak{g}$ defines a natural action of $G$ on its Lie algebra called the adjoint action, denoted $Ad_g := T_eI_g$ for all $g \in G$. Via the natural dual pairing, this induces the dual map $T_e^*I_g : \mathfrak{g}^* \to \mathfrak{g}^*$, which is the
corresponding action of $G$ on the dual of its Lie algebra, called the coadjoint action and denoted $Ad^*_g := T_e^*I_g$. The coadjoint action plays a fundamental role in mechanics, since its orbits carry a canonical symplectic structure. The infinitesimal generators of the adjoint and coadjoint actions are denoted by $ad_\xi$ and $ad^*_\xi$, respectively, for each $\xi \in \mathfrak{g}$ as shown below.
\begin{center}
\begin{tikzpicture}
\shadedraw[white] (3,1) to[out=0,in=180] (6,2) to[out=0,in=0] (4,4) to[out=180,in=90] (2,3) to[out=270,in=180] (3,1);
\node at (5,1.5) {$G$};
\draw[orange] (3,2) to[out=68,in=190] (5,3.5);
\node at (3,1.8) {\tiny $\textcolor{orange}{\exp{(t\xi)}}$};
\fill[] (4,3.17) circle (0.8 pt);
\node at (4,3) {\tiny $e$};
\draw[dotted] (4,3.17)--(4,5.17);
\fill[red] (4,4.17) circle (0.8 pt);
\draw[red,fill=red,fill opacity=0.2] (3.5,3.5)--(5.7,4)--(4.5,4.8)--(2.3,4.3)--cycle; 
\draw[thick,-stealth,red] (4,4.17)--(4.5,4.5);
\node at (4.6,4.5) {\tiny $\textcolor{red}{\xi}$};
\draw[thick,-stealth,red] (4,4.17)--(4.6,4);
\node at (4.7,3.9) {\tiny $\textcolor{red}{\eta}$};
\draw[orange] (4,3.8) to[out=10,in=240] (5,4.4);
\draw[-stealth,thick,red] (4.6,4)--(5,4.25);
\node at (5.4,4.15) {\tiny $\textcolor{red}{ad_\xi (\eta)}$};
\node at (3,3.5) {$\textcolor{red}{\mathfrak{g}}$};
\node at (5.75,4.5) {\tiny $\textcolor{orange}{Ad_{\exp{(t\xi)}}(\eta)}$};
\fill[blue] (4,5.17) circle (0.8 pt);
\draw[blue,fill=blue,fill opacity=0.2] (3.5,4.5)--(5.7,5)--(4.5,5.8)--(2.3,5.3)--cycle;
\draw[thick,-stealth,blue] (4,5.17)--(4.6,5);
\node at (4.7,4.9) {\tiny $\textcolor{blue}{\mu}$};
\draw[orange] (4,4.8) to[out=10,in=240] (5,5.4);
\draw[-stealth,thick,blue] (4.6,5)--(5,5.25);
\node at (5.4,5.15) {\tiny $\textcolor{blue}{ad^*_\xi (\mu)}$};
\node at (5.75,5.5) {\tiny $\textcolor{orange}{Ad^*_{\exp{(t\xi)}}(\mu)}$};
\node at (3,5.8) {$\textcolor{blue}{\mathfrak{g}^*}$};
\draw[<-] (2.2,4.4) arc (30:330:5 pt);
\node at (1.5,4.4) {\footnotesize $Ad$};
\draw[<-] (2.2,5.4) arc (30:330:5 pt);
\node at (1.5,5.4) {\footnotesize $Ad^*$};
\end{tikzpicture}
\end{center}

\chapter{Retraction and Discretization Maps}
\label{ch3}
In this chapter, we introduce the mathematical machinery required for the intrinsic construction of numerical integrators for dynamical systems evolving on manifolds. The framework presented here is based on the approaches developed in \cite{barbero2023retraction}, \cite{BARBEROLINAN2024155}, \cite{10.1007/978-3-031-10047-5_65}, and \cite{10383230}. 
\section{Retraction Maps}
\label{sec retraction maps}
Intuitively, a retraction map takes a tangent vector at a point on a manifold and maps it to a nearby point on the manifold in the direction of that vector. This is accomplished by generating a curve that emanates from the given point with the prescribed initial velocity, as illustrated in the figure \ref{fig retraction}, and evaluating this curve at unit time, see \cite{absil2007optimization}. The curve is required to satisfy certain natural compatibility conditions, which are formalized in the following definition.
\begin{figure}
\centering
    \begin{tikzpicture}[scale = 1]
        \draw[gray,fill=gray,fill opacity=0.1] (1,3) to[out=15,in=105] (3,1) to[out=-5,in=-175] (7,1) to[out=105,in=15] (5,3) to[out=-175,in=-5] (1,3);
        \node at (2,1.5) {$M$};
        \draw[fill] (4,2.5) circle (1pt);
        \draw (4,2.5) to[out=-25,in=115] (5,1.5);
        \draw[fill] (4.7,2) circle (1pt);
        \node at (3.7,2.5) {$x$};
        \node at (5.7,1.5) {$\mathcal{R}(x,tv)$};
        \node at (5.4,2.2) {$\mathcal{R}(x,v)$};
        \draw[dotted] (4,2.5)--(4,4);
        \draw[fill=gray,fill opacity=0.4] (3.5,3.5)--(6,3.5)--(5,4.5)--(2.5,4.5)--cycle;
        \draw[fill] (4,4) circle (1pt);
        \draw[thick,-stealth] (4,4)--(4.8,3.7);
        \node at (3.7,4) {$0$};
        \node at (5,3.7) {$v$};
        \node at (2.5,3.5) {$T_xM$};
    \end{tikzpicture}
    \caption{Retraction map on $M$}
    \label{fig retraction}
\end{figure}
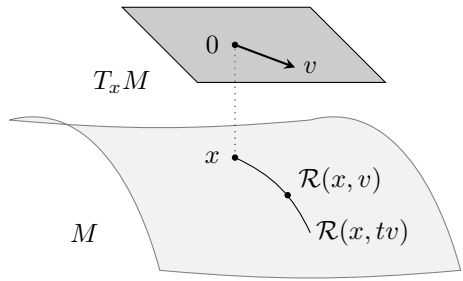
\begin{definition} \label{def retraction map}
      Let $M$ be a smooth manifold and $TM$ its tangent bundle. A smooth map 
\begin{align*}
    \mathcal{R} : TM &\longrightarrow M \\
    (x,v) &\longmapsto \mathcal{R}(x,v)
\end{align*}
satisfying
\begin{itemize}
    \item (R1) $\mathcal{R}(x,0)=x$ and
    \item (R2) $\left.\frac{d}{dt}\right|_{t=0}\mathcal{R}(x,tv)=v$ 
\end{itemize}
for every $(x,v) \in TM$ is called a \emph{retraction map}.
\end{definition}
Condition (R1) simply states that, in the absence of any velocity, the point on the manifold remains fixed. Condition (R2) ensures that, as the velocity is scaled, the point moves away from the base point along the manifold in the corresponding direction, with the correct first-order behavior. We now turn to some illustrative examples of retraction maps.
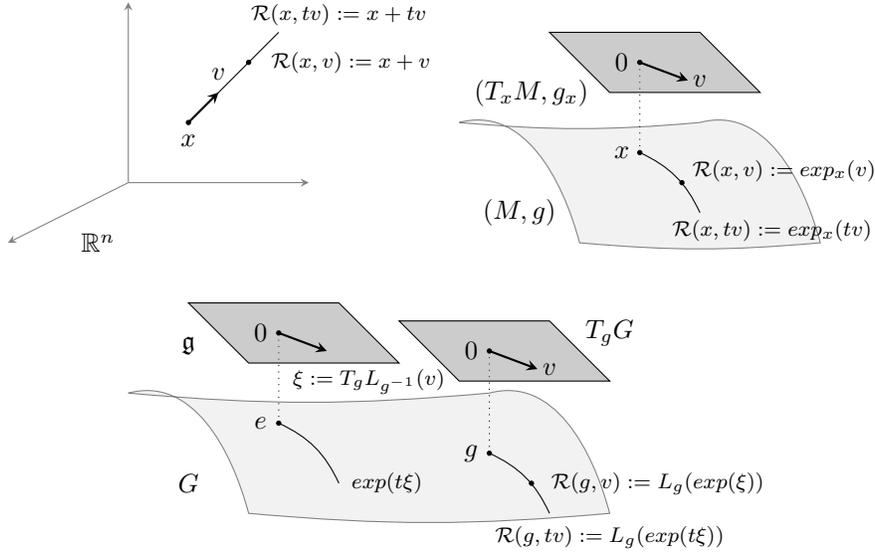
\begin{figure}
    \centering
    \begin{tikzpicture}[scale = 0.8]
        \draw[gray,-stealth] (3,2)--(6,2);
        \draw[gray,-stealth] (3,2)--(3,5);
        \draw[gray,-stealth] (3,2)--(1,1);
        \node at (2.5,1) {$\mathbb{R}^n$};
        \draw[fill] (4,3) circle (1pt);
        \draw (4,3)--(5.5,4.5);
        \draw[thick,-stealth] (4,3)--(4.5,3.5);
        \draw[fill] (5,4) circle (1pt);
        \node at (4,2.7) {$x$};
        \node at (4.5,3.8) {$v$};
        \node at (6.5,4.8) {\footnotesize $\mathcal{R}(x,tv):=x+tv$};
        \node at (6.7,4) {\footnotesize $\mathcal{R}(x,v):=x+v$};
        \draw[gray,fill=gray,fill opacity=0.1] (8.5,3) to[out=15,in=105] (10.5,1) to[out=-5,in=-175] (14.5,1) to[out=105,in=15] (12.5,3) to[out=-175,in=-5] (8.5,3);
        \node at (9.5,1.5) {$(M,g)$};
        \draw[fill] (11.5,2.5) circle (1pt);
        \draw (11.5,2.5) to[out=-25,in=115] (12.5,1.5);
        \draw[fill] (12.2,2) circle (1pt);
        \node at (11.2,2.5) {$x$};
        \node at (13.7,1.2) {\footnotesize $\mathcal{R}(x,tv):=exp_x(tv)$};
        \node at (13.9,2.2) {\footnotesize $\mathcal{R}(x,v):=exp_x(v)$};
        \draw[dotted] (11.5,2.5)--(11.5,4);
        \draw[fill=gray,fill opacity=0.4] (11,3.5)--(13.5,3.5)--(12.5,4.5)--(10,4.5)--cycle;
        \draw[fill] (11.5,4) circle (1pt);
        \draw[thick,-stealth] (11.5,4)--(12.3,3.7);
        \node at (11.2,4) {$0$};
        \node at (12.5,3.7) {$v$};
        \node at (9.7,3.5) {$(T_xM,g_x)$};
        \draw[gray,fill=gray,fill opacity=0.1] (3,-1.5) to[out=15,in=105] (5,-3.5) to[out=-5,in=-175] (11,-3.5) to[out=105,in=15] (9,-1.5) to[out=-175,in=-5] (3,-1.5);
        \node at (4,-3) {$G$};
        \draw[fill] (5.5,-2) circle (1pt);
        \draw (5.5,-2) to[out=-25,in=115] (6.5,-3);
        \node at (5.2,-2) {$e$};
        \node at (7.3,-3) {\footnotesize $exp(t\xi)$};
        \draw[fill] (9,-2.5) circle (1pt);
        \draw (9,-2.5) to[out=-25,in=115] (10,-3.5);
        \draw[fill] (9.7,-3) circle (1pt);
        \node at (8.7,-2.5) {$g$};
        \node at (11,-3.9) {\footnotesize $\mathcal{R}(g,tv):=L_g(exp(t\xi))$};
        \node at (11.8,-3) {\footnotesize $\mathcal{R}(g,v):=L_g(exp(\xi))$};
        \draw[dotted] (5.5,-2)--(5.5,-0.5);
        \draw[fill=gray,fill opacity=0.4] (5,-1)--(7.5,-1)--(6.5,0)--(4,0)--cycle;
        \draw[fill] (5.5,-0.5) circle (1pt);
        \draw[thick,-stealth] (5.5,-0.5)--(6.3,-0.8);
        \node at (5.2,-0.5) {$0$};
        \node at (7,-1.3) {\footnotesize $\xi:=T_gL_{g^{-1}}(v)$};
        \node at (4,-0.7) {$\mathfrak{g}$};
        \draw[dotted] (9,-2.5)--(9,-0.8);
        \draw[fill=gray,fill opacity=0.4] (8.5,-1.3)--(11,-1.3)--(10,-0.3)--(7.5,-0.3)--cycle;
        \draw[fill] (9,-0.8) circle (1pt);
        \draw[thick,-stealth] (9,-0.8)--(9.8,-1.1);
        \node at (8.7,-0.8) {$0$};
        \node at (10,-1.1) {$v$};
        \node at (11,-0.5) {$T_gG$};
    \end{tikzpicture}
    \caption{Retraction maps on $\mathbb{R}^n$ (top left), on a Riemannian manifold $(M,g)$ (top right) and on a Lie group $G$ (bottom)}
    \label{fig retraction examples}
\end{figure}
\begin{example} \label{eg retraction R^n}
    On $\mathbb{R}^n$, a retraction map can be defined in the simplest possible way: given a point $x \in \mathbb{R}^n$ and a tangent vector $v \in T_x \mathbb{R}^n \cong \mathbb{R}^n$, the map sends $(x,v)$ to the point $x+v \in \mathbb{R}^n$, as illustrated in figure \ref{fig retraction examples}. 
\end{example}
\begin{example} \label{eg retraction (M,g)}
     Let $(M,g)$ be a Riemannian manifold. A natural retraction map can be defined using the Riemannian exponential map; see, for example, \cite{do1992riemannian}, \cite{postnikov2013geometry}, or \cite{lee2019introduction}. Recall that every Riemannian metric determines a unique affine connection, the Levi–Civita connection, which is metric-compatible and torsion-free. This connection allows one to define geodesics (or autoparallel curves) on the manifold given an initial point and velocity. The exponential map then sends a pair $(x,v) \in TM$ to the point reached at unit time along the geodesic $\exp_x{(tv)}$ starting at $x$ with initial velocity $v$, as illustrated in figure \ref{fig retraction examples}.
\end{example}
\noindent Note that the straight line in the previous example is precisely a geodesic with respect to the standard Euclidean metric on $\mathbb{R}^n$.
\begin{example} \label{eg retraction G}
    Let $G$ be a Lie group. A natural retraction map on $G$ can be defined using the group exponential map (see \cite{helgason1979differential}, \cite{olver1993applications}, or \cite{lee2012introduction}). Recall that the group exponential map is constructed from the integral curves of left-invariant vector fields on $G$. Given $v \in T_gG$, we first left-translate it to the identity to obtain $\xi := T_gL_{g^{-1}}(v) \in T_eG \cong \mathfrak{g}$. We then form the associated one-parameter subgroup $\exp{(t\xi)}$ and translate it back to $g$ via left multiplication. The retraction map thus sends $(g,v) \in TG$ to the point one unit of time along the curve $L_g(\exp{(t\xi)})$, as illustrated in figure \ref{fig retraction examples}.    
\end{example}
\noindent Note that when $G$ is equipped with a bi-invariant Riemannian metric, the group exponential map coincides with the Riemannian exponential map. Also note that one may equally use right translation instead of left translation to define a retraction map. In essence, it suffices to have a mapping from the Lie algebra to the Lie group, which leads naturally to the notion of trivialized retraction maps.
\subsection{Trivialized Retraction Maps}
\label{subsec triv retraction maps}
We now specialize retraction maps to the setting of Lie groups, where they play a central role in the construction of numerical integrators for Lie-group–valued dynamics. We begin by discussing the trivialization of the tangent and cotangent bundles of a Lie group, see \cite{vivek2025lcss} .
\paragraph{Tangent bundle of a Lie group}
\label{triv TG}
The tangent bundle $\tau_G : TG \to G$ of a Lie group $G$ carries a natural Lie group structure, see \cite{michortopics}, with group multiplication and inversion given by
\begin{align*}
    u_g \, v_h &:= (T_gR_h(u) + T_hL_g(v))_{gh} \quad \text{and} \\
    (u_g)^{-1} &:= (-T_eL_{g^{-1}} \circ T_gR_{g^{-1}}(u))_{g^{-1}}.
\end{align*}
Here concatenation denotes the group operation on $TG$, with $0_e$ serving as the identity element. An element of $TG$ is denoted as $u_g$ where $\tau_G(u_g) := g \in G$ is the base point and $u \in \tau^{-1}_G(g) = T_gG$ is the tangent vector in the fiber over $g$ as illustrated in the figure below.
\begin{center}
    \begin{tikzpicture}
        \draw (3,1)--(3,3);
        \node at (3,0.7) {\scriptsize $T_gG$};
        \draw[red,thick,-stealth] (3,2)--(3,2.4);
        \node at (3.2,2.4) {\tiny $\textcolor{red}{u}$};
        \draw (4,1)--(4,3);
        \node at (4,0.7) {\scriptsize $T_hG$};
        \draw[blue,thick,-stealth] (4,2)--(4,2.8);
        \node at (4.2,2.8) {\tiny $\textcolor{blue}{v}$};
        \draw (5,1)--(5,3);
        \node at (5,0.7) {\scriptsize $T_{gh}G$};
        \draw[blue,thick,-stealth] (5,2)--(5,2.8);
        \node at (5.7,2.8) {\tiny $\textcolor{blue}{T_hL_g(v)}$};
        \draw[red,thick,-stealth] (5,2)--(5,2.4);
        \node at (5.7,2.4) {\tiny $\textcolor{red}{T_gR_h(u)}$};
        \draw[dashed] (1.5,2)--(6.5,2);
        \draw[fill] (3,2) circle (1pt);
        \node at (3.2,1.8) {\tiny $g$};
        \draw[fill] (4,2) circle (1pt);
        \node at (4.2,1.8) {\tiny $h$};
        \draw[fill] (5,2) circle (1pt);
        \node at (5.3,1.8) {\tiny $gh$};
        \node at (1,2) {$G$};
        \draw (7,1) to[out=10,in=190] (7.2,2) to[out=170,in=-10] (7,3);
        \node at (7.7,2) {$TG$};
    \end{tikzpicture}
\end{center}
Using left translations, the tangent bundle admits a natural trivialization.
\begin{definition} \label{left-trivialization of TG}
    The map $tr^L_{TG} : TG \to G \times \mathfrak{g}$ defined by
    \begin{align*}
        tr^L_{TG}(u_g) := (g,T_gL_{g^{-1}}(u))
    \end{align*}
    for all $u_g \in TG$ is called the \emph{left-trivialization} of $TG$.
\end{definition}
Analogously, one may define a \emph{right-trivialization} map $tr^R_{TG} : TG \to G \times \mathfrak{g}$ using right translations.
Via the left trivialization map, the Lie group structure on $TG$ can be transported to $G \times \mathfrak{g}$, endowing it with the group multiplication and inversion given by
\begin{align*}
    (g,\xi) (h,\eta) &:= (gh,Ad_{h^{-1}}\xi + \eta) \quad \text{and} \\
        (g,\xi)^{-1} &:= (g^{-1},-Ad_g\xi),
\end{align*}
respectively, for $g,h \in G$ and $\xi,\eta \in \mathfrak{g}$. This construction equips $G \times \mathfrak{g}$ with a semidirect product group structure, and we therefore identify $TG \cong G \ltimes \mathfrak{g}$.
\paragraph{Cotangent bundle of a Lie group}
\label{triv T^*G}
The cotangent bundle $\pi_G : T^*G \to G$ of a Lie group $G$ also carries a natural Lie group structure with group multiplication and inversion given by
\begin{align*}
    \alpha_g \, \beta_h &:= (T^*_{gh}R_{h^{-1}}(\alpha) + T^*_{gh}L_{g^{-1}}(\beta))_{gh} \quad \text{and} \\
    (\alpha_g)^{-1} &:= (-T_{g^{-1}}^*L_g \circ T_e^*R_g(\alpha))_{g^{-1}}.
\end{align*}
Here again concatenation denotes the group operation on $T^*G$, with $0_e$ serving as the identity element. An element of $T^*G$ is denoted as $\alpha_g$ where $\pi_G(\alpha_g) := g \in G$ is the base point and $\alpha \in \pi^{-1}_G(g) = T_g^*G$ is the cotangent vector in the fiber over $g$ as illustrated in the figure below.
\begin{center}
    \begin{tikzpicture}
        \draw (3,1)--(3,3);
        \node at (3,0.7) {\scriptsize $T_g^*G$};
        \draw[red,thick,-stealth] (3,2)--(3,2.4);
        \node at (3.2,2.4) {\tiny $\textcolor{red}{\alpha}$};
        \draw (4,1)--(4,3);
        \node at (4,0.7) {\scriptsize $T_h^*G$};
        \draw[blue,thick,-stealth] (4,2)--(4,2.8);
        \node at (4.2,2.8) {\tiny $\textcolor{blue}{\beta}$};
        \draw (5,1)--(5,3);
        \node at (5,0.7) {\scriptsize $T_{gh}^*G$};
        \draw[blue,thick,-stealth] (5,2)--(5,2.8);
        \node at (5.9,2.8) {\tiny $\textcolor{blue}{T_{gh}^*R_{h^{-1}}(\alpha)}$};
        \draw[red,thick,-stealth] (5,2)--(5,2.4);
        \node at (5.9,2.4) {\tiny $\textcolor{red}{T_{gh}^*L_{g^{-1}}(\beta)}$};
        \draw[dashed] (1.5,2)--(6.5,2);
        \draw[fill] (3,2) circle (1pt);
        \node at (3.2,1.8) {\tiny $g$};
        \draw[fill] (4,2) circle (1pt);
        \node at (4.2,1.8) {\tiny $h$};
        \draw[fill] (5,2) circle (1pt);
        \node at (5.3,1.8) {\tiny $gh$};
        \node at (1,2) {$G$};
        \draw (7,1) to[out=10,in=190] (7.2,2) to[out=170,in=-10] (7,3);
        \node at (7.7,2) {$T^*G$};
    \end{tikzpicture}
\end{center}
Using left translations, the cotangent bundle also admits a natural trivialization.
\begin{definition} \label{left-trivialization of T*G}
     The map $tr^L_{T^*G} : T^*G \to G \times \mathfrak{g}^*$ defined by
    \begin{align*}
        tr^L_{T^*G}(\alpha_G) := (g,T_e^*L_g(\alpha))
    \end{align*}
    for all $\alpha_g \in T^*G$ is called the \emph{left-trivialization} of $T^*G$.
\end{definition}
Analogously, one may define a \emph{right-trivialization} map $tr^R_{T^*G} : T^*G \to G \times \mathfrak{g}$ using right translations. Via the left trivialization map, the Lie group structure on $T^*G$ can also be transported to $G \times \mathfrak{g}^*$, endowing it with the group multiplication and inversion given by
\begin{align*}
        (g,\mu)(h,\nu) &:= (gh, Ad^*_h \mu + \nu) \quad \text{and} \\
        (g,\mu)^{-1} & := (g^{-1}, -Ad^*_{g^{-1}}\mu),
\end{align*}
respectively, for $g,h \in G$ and $\mu,\nu \in \mathfrak{g}^*$. This construction endows $G \times \mathfrak{g}^*$ with a semidirect product Lie group structure, yielding the identification $T^*G \cong G \ltimes \mathfrak{g}^*$. Recall that the cotangent bundle of any manifold carries a canonical symplectic structure. Via the left–trivialization map, this canonical symplectic form on $T^*G$ can be pushed forward to $G \ltimes \mathfrak{g}^*$, see \cite{ortega2003momentum}. Thus, $T^*G$ and $G \ltimes \mathfrak{g}^*$ are naturally identified not only as vector bundles and Lie groups, but also as symplectic manifolds.
\vspace{10pt}\\
We now define retraction maps adapted to the left trivialization of the tangent bundle $TG$ of a Lie group $G$.
\begin{definition} \label{def triv retraction map}
    A smooth map $\mathcal{R}^L : G \ltimes \mathfrak{g} \to G$ satisfying
    \begin{itemize}
        \item (TR1) $\mathcal{R}^L(g,0) = g$ and
        \item (TR2) $\left. \frac{d}{dt} \right|_{t=0} \mathcal{R}^L(g,t\xi) = T_eL_g(\xi)$
    \end{itemize}
    for all $(g,\xi) \in G \ltimes \mathfrak{g}$ is called a \emph{left-trivialized retraction} map on $G$.
\end{definition}
A right-trivialized retraction map $\mathcal{R}^R : G \ltimes \mathfrak{g} \to G$ is defined analogously with condition (TR2) replaced by $\left. \frac{d}{dt} \right|_{t=0} \mathcal{R}^R(g,t\xi) = T_eR_g(\xi)$. Equivalently, a left-trivialized retraction can be constructed from any local diffeomorphism $\tau : \mathfrak{g} \to G$
satisfying $\tau(0) = e$ and $\left. \frac{d}{dt} \right|_{t=0} \tau(t\xi) = \xi$, as stated in the following proposition.
\begin{figure}
    \centering
    \begin{tikzpicture}[scale = 1]
        \draw[gray,fill=gray,fill opacity=0.1] (1,3) to[out=15,in=105] (3,1) to[out=-5,in=-175] (9,1) to[out=105,in=15] (7,3) to[out=-175,in=-5] (1,3);
        \node at (2,1.5) {$G$};
        \draw[fill] (3.5,2.5) circle (1pt);
        \draw (3.5,2.5) to[out=-25,in=115] (4.5,1.5);
        \draw[fill] (4.2,2) circle (1pt);
        \node at (3.2,2.5) {$e$};
        \node at (4.5,1.2) {$\tau(t\xi)$};
        \node at (4.7,2.1) {$\tau(\xi)$};
        \draw[fill] (7,2) circle (1pt);
        \draw (7,2) to[out=-25,in=115] (8,1);
        \draw[fill] (7.7,1.5) circle (1pt);
        \node at (6.7,2) {$g$};
        \node at (9,0.6) {\footnotesize $\mathcal{R}^L(g,t\xi):=L_g(\tau(t\xi))$};
        \node at (9.5,1.5) {\footnotesize $\mathcal{R}^L(g,\xi):=L_g(\tau(\xi))$};
        \draw[dotted] (3.5,2.5)--(3.5,4);
        \draw[fill=red,fill opacity=0.1,red] (3,3.5)--(5.5,3.5)--(4.5,4.5)--(2,4.5)--cycle;
        \draw[red,fill] (3.5,4) circle (1pt);
        \draw[red,thick,-stealth] (3.5,4)--(4.3,3.7); 
        \node at (3.2,4) {$\textcolor{red}{0}$};
        \node at (4,4.1) {$\textcolor{red}{\xi}$};
        \node at (2,3.8) {$\textcolor{red}{\mathfrak{g}}$};
        \draw[dotted] (7,2)--(7,3.7);
        \draw[gray,fill=gray,fill opacity=0.3] (6.5,3.2)--(9,3.2)--(8,4.2)--(5.5,4.2)--cycle;
        \draw[fill] (7,3.7) circle (1pt);
        \draw[thick,-stealth] (7,3.7)--(7.8,3.4);
        \node at (6.7,3.7) {$0$};
        \node at (7.8,3.8) {$T_eL_g(\xi)$};
        \node at (9,4) {$T_gG$};
\end{tikzpicture}
    \caption{Local diffeomorphism $\tau$ and left-trivialized retraction map $\mathcal{R}^L$ on a Lie group $G$}
    \label{fig trivialized retraction}
\end{figure}
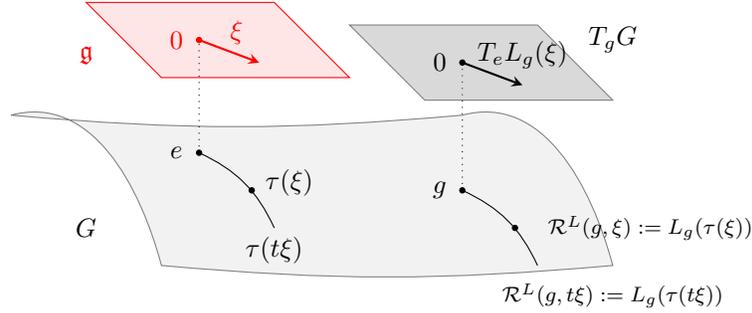
\begin{proposition} \label{prop 2.1}
    A map $\mathcal{R}^L : G \ltimes \mathfrak{g} \to G$ defined as 
    \begin{align*}
        \mathcal{R}^L(g,\xi) := L_g(\tau(\xi)) = g\tau(\xi)
    \end{align*}
    for all $(g,\xi) \in G \times \mathfrak{g}$ is a left-trivialized retraction map on $G$ where $\tau$ is as defined above.
\end{proposition}
Similarly, a \emph{right-trivialized} retraction map $\mathcal{R}^R : G \ltimes \mathfrak{g} \to G$ is defined by replacing left translation with right translation as $\mathcal{R}^R(g,\xi) := R_g(\tau(\xi)) = \tau(\xi) g$. This characterization makes explicit that a trivialized retraction requires only a smooth map from the Lie algebra to the Lie group that is first-order compatible at the identity, as illustrated in figure \ref{fig trivialized retraction}. A natural choice for $\tau$ is the group exponential map, however, we retain this more general formulation since alternative choices with advantageous properties exist for certain Lie groups, see \cite{bou2009hamilton}. We take this opportunity to introduce the notion of the logarithmic derivative, see \cite{michortopics}, in particular for maps $\tau : \mathfrak{g} \to G$, which will play a key role in the constructions later on.
\paragraph{Logarithmic derivatives} 
\label{log derivative}
Let $\tau : G \to \mathfrak{g}$ be a local diffeomorphism as discussed above. Consider the affine line $\xi + t \eta$ in $\mathfrak{g}$. Under $\tau$ , this line gets mapped to a curve in $G$. The velocity of this curve at $t=0$ can be translated back to the identity by left or right translation, as illustrated in the figure below. This construction yields the \emph{left} and \emph{right logarithmic derivatives} of $\tau$, which are formalized in the following definition.
\begin{definition} \label{def log derivative}
    A smooth map $d^L_\xi \tau : \mathfrak{g} \to \mathfrak{g}$ defined as
    \begin{align*}
        d^L_\xi \tau (\eta) := T_{\tau(\xi)}L_{\tau(\xi)^{-1}} \circ T_\xi \tau (\eta)
    \end{align*}
    for all $\eta \in \mathfrak{g}$ is called the \emph{left-logarithmic derivative} of $\tau : \mathfrak{g} \to G$ at $\xi \in \mathfrak{g}$.
\end{definition}
\begin{center}
\begin{tikzpicture}[scale = 1]
        \draw[gray,fill=gray,fill opacity=0.1] (1,3) to[out=15,in=105] (3,1) to[out=-5,in=-175] (9,1) to[out=105,in=15] (7,3) to[out=-175,in=-5] (1,3);
        \node at (2,1.5) {$G$};
        \draw[->] (4.7,4) to[out=-75,in=75] (4.7,3);
        \node at (5,3.5) {$\tau$};
        \draw[fill] (4.5,2) circle (1pt);
        \draw (3,2.1) to[out=5,in=135] (6,1.1);
        \node at (3,0.5) {\footnotesize $R_{\tau(\xi)^{-1}}(\tau(\xi+t\eta))$};
        \draw (3.1,1.1) to[out=45,in=180] (5.9,2.2);
        \node at (6,0.5) {\footnotesize $L_{\tau(\xi)^{-1}}(\tau(\xi+t\eta))$};
        \node at (4.5,1.7) {$e$};
        \draw[fill] (7.55,2) circle (1pt);
        \node at (7.2,1.7) {$\tau(\xi)$};
        \draw (6.5,2.7) to[out=-25,in=115] (8.1,1.2);
        \node at (8.4,0.5) {\footnotesize$\tau(\xi+t\eta)$};
        \draw[dotted] (4.5,2)--(4.5,5);
        \draw[fill=red,fill opacity=0.1,red] (3.5,4.3)--(7,4.3)--(6,5.7)--(2.5,5.7)--cycle;
        \draw[red,fill] (4.5,5) circle (1pt);
        \draw[red,thick,-stealth] (4.5,5)--(5.25,5); 
        \node at (5.3,5.2) {\textcolor{red}{$\xi$}};
        \draw[red,thick,-stealth] (5.25,5)--(5.7,4.82);
        \node at (5.8,5) {\textcolor{red}{$\eta$}};
        \draw[red] (4,5.5)--(6.5,4.5);
        \node at (7,4) {$\textcolor{red}{\xi+t\eta}$};
        \node at (4.3,5.1) {$\textcolor{red}{0}$};
        \node at (2.5,4.8) {$\textcolor{red}{\mathfrak{g}}$};
        \draw[thick,orange,-stealth] (4.5,5)--(5.4,4.4);
        \node at (5.5,4) {\textcolor{orange}{$d^L_\xi\tau(\eta)$}};
        \draw[thick,orange,-stealth] (4.5,5)--(3.6,4.4);
        \node at (3.5,4) {\textcolor{orange}{$d^R_\xi\tau(\eta)$}};
\end{tikzpicture}
\end{center}
Similarly, the \emph{right-logarithmic derivative} $d^R_\xi\tau : \mathfrak{g} \to \mathfrak{g}$ is defined as $d^R_\xi \tau (\eta) := T_{\tau(\xi)}R_{\tau(\xi)^{-1}} \circ T_\xi \tau (\eta)$ for all $\eta \in \mathfrak{g}$. The left- and right-logarithmic derivatives are related by the adjoint action:
\begin{align*} \label{left-right log derivative}
    d^R_\xi\tau := Ad_{\tau(\xi)} \circ d^L_\xi
\end{align*}
for all $\xi \in \mathfrak{g}$.
\section{Discretization Maps}
\label{sec discretization maps}
A discretization map generalizes the notion of a retraction map and, as its name suggests, is intended for discretization. Unlike a retraction map, which maps $TM$ to $M$, a discretization map assigns to each $(x,v) \in TM$ an ordered pair of points in $M$, as illustrated in figure \ref{fig discretization on M}. It is required to satisfy certain natural properties, formalized in the following definition.
\begin{figure}
    \centering
    \begin{tikzpicture}[scale = 0.8]
        \draw[gray,fill=gray,fill opacity=0.1] (1,3) to[out=15,in=105] (3,1) to[out=-5,in=-175] (7,1) to[out=105,in=15] (5,3) to[out=-175,in=-5] (1,3);
        \node at (2,1.5) {$M$};
        \draw[fill] (4,2.5) circle (1pt);
        \draw[] (4,2.5) to[out=-30,in=110] (5.2,1);
        \draw[] (4,2.5) to[out=0,in=160] (6.4,2);
        \draw[fill] (4.65,2) circle (1pt);
        \draw[fill] (5,2.42) circle (1pt);
        \node at (3.7,2.5) {$x$};
        \node at (5.5,2.7) {\footnotesize $D^2(x,v)$};
        \node at (3.8,2) {\footnotesize $D^1(x,v)$};
        \node at (5,0.6) {\footnotesize $D^1(x,tv)$};
        \node at (7.5,2) {\footnotesize $D^2(x,tv)$};
        \draw[dotted] (4,2.5)--(4,4);
        \draw[fill=gray,fill opacity=0.4] (3.5,3.5)--(6,3.5)--(5,4.5)--(2.5,4.5)--cycle;
        \draw[fill] (4,4) circle (1pt);
        \draw[thick,-stealth] (4,4)--(4.8,3.7);
        \node at (3.7,4) {$0$};
        \node at (5,3.7) {$v$};
        \node at (2.5,3.5) {$T_xM$};
        \draw[gray,fill=gray,fill opacity=0.1] (8,3) to[out=15,in=105] (10,1) to[out=-5,in=-175] (14,1) to[out=105,in=15] (12,3) to[out=-175,in=-5] (8,3);
        \node at (9,1.5) {$M$};
        \draw[fill] (11,2.6) circle (1pt);
        \draw (10.5,2.8) to[out=-15,in=105] (12,1);
        \draw[fill] (11.6,2) circle (1pt);
        \draw[fill] (11.95,1.2) circle (1pt);
        \node at (9.8,2.5) {\footnotesize $\mathcal{R}(x,-\theta v)$};
        \node at (13.4,1.2) {\footnotesize $\mathcal{R}(x,(1-\theta)v)$};
        \node at (12,2) {$x$};
        \draw[dotted] (11.6,2)--(11.6,4);
        \draw[fill=gray,fill opacity=0.4] (11.1,3.5)--(13.6,3.5)--(12.6,4.5)--(10.1,4.5)--cycle;
        \draw[fill] (11.6,4) circle (1pt);
        \draw[thick,-stealth] (11.6,4)--(12.4,3.7);
        \node at (11.3,4) {$0$};
        \node at (12.6,3.7) {$v$};
        \node at (10.1,3.5) {$T_xM$};
        \node at (12,0.6) {\footnotesize $\mathcal{R}(x,tv)$};
    \end{tikzpicture}
    \caption{Discretization map on $M$ (left) and discretization map induced by a retraction on $M$ (right)}
    \label{fig discretization on M}
\end{figure}
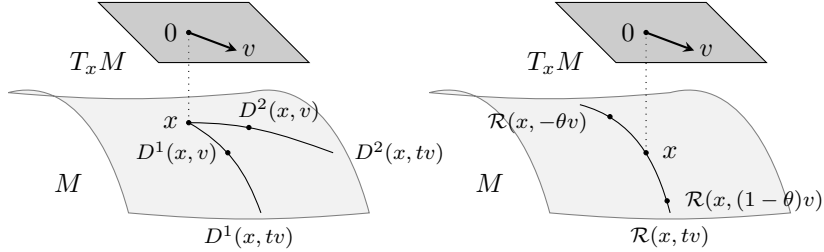
\begin{definition} \label{def discretization map}
    Let $M$ be a smooth manifold and $TM$ its tangent bundle. A smooth map
\begin{align*}
    \mathcal{D}:TM &\longrightarrow M \times M \\
    (x,v) &\longmapsto (D^1(x,v),D^2(x,v))
\end{align*}
satisfying
\begin{itemize}
    \item (D1) $\mathcal{D}(x,0)= (D^1(x,0),D^2(x,0))=(x,x)$ and
    \item (D2) $\left.\frac{d}{dt}\right|_{t=0}D^2(x,tv)-\left.\frac{d}{dt}\right|_{t=0}D^1(x,tv)=v$ 
\end{itemize}
for every $(x,v) \in TM$ is called a \emph{discretization map}.
\end{definition}
The conditions (D1) and (D2) ensure that, when 
$D^1(x,v) = x$, the map $D^2 : TM \to M$ is a retraction. Moreover, a single retraction map naturally induces an entire family of discretization maps, as illustrated in figure \ref{fig discretization on M}.
\begin{proposition}
    Let $\mathcal{R}: TM \to M$ be a retraction map on $M$. A smooth map 
    \begin{align*}
        \mathcal{D}: TM &\longrightarrow M \times M\\
        (x,v) &\longmapsto (\mathcal{R}(x,-\theta v),\mathcal{R}(x,(1-\theta)v))
    \end{align*}
    is a discretization map for every $\theta \in [0,1]$.
\end{proposition}
Since $TM$ and $M \times M$ have the same dimension, it is natural to ask whether a discretization map $\mathcal{D} : TM \to M \times M$ is a diffeomorphism. In general $\mathcal{D}$ need not be globally invertible, however, it is a local diffeomorphism.
\begin{proposition}
    Let $\mathcal{D}: TM \to M \times M$ be a discretization map on $M$ then $\mathcal{D}$ is a local diffeomorphism around the zero section of $TM$.
\end{proposition}
Having established that a retraction map induces a family of discretization maps, each of the retraction maps introduced in the previous section gives rise to a corresponding discretization map, as illustrated in figure \ref{fig discretization examples}.
\begin{figure}
    \centering
    \begin{tikzpicture}[scale = 0.8]
        \draw[gray,-stealth] (3,2)--(6,2);
        \draw[gray,-stealth] (3,2)--(3,5);
        \draw[gray,-stealth] (3,2)--(1,1);
        \node at (2.5,1) {$\mathbb{R}^n$};
        \draw[fill] (4,3) circle (1pt);
        \draw (4,3)--(5.5,4.5);
        \draw (4,3)--(3.2,2.2);
        \draw[thick,-stealth] (4,3)--(4.5,3.5);
        \draw[fill] (5,4) circle (1pt);
        \draw[fill] (3.5,2.5) circle (1pt);
        \node at (4,3.3) {\footnotesize$x$};
        \node at (4.5,3.8) {\footnotesize$v$};
        \node at (5.5,4.8) {\footnotesize $x+tv$};
        \node at (6.2,4) {\footnotesize $x+(1-\theta)v$};
        \node at (4.2,2.5) {\footnotesize $x-\theta v$};
        \draw[gray,fill=gray,fill opacity=0.1] (8.5,3) to[out=15,in=105] (10.5,1) to[out=-5,in=-175] (14.5,1) to[out=105,in=15] (12.5,3) to[out=-175,in=-5] (8.5,3);
        \node at (9.5,1.5) {$(M,g)$};
        \draw[fill] (11.5,2.6) circle (1pt);
        \draw (11,2.8) to[out=-15,in=105] (12.5,1);
        \draw[fill] (12.1,2) circle (1pt);
        \draw[fill] (12.45,1.2) circle (1pt);
        \node at (12.8,2.5) {\footnotesize $\exp_x(-\theta v)$};
        \node at (14,1.2) {\footnotesize $\exp_x((1-\theta)v)$};
        \node at (12.5,2) {$x$};
        \draw[dotted] (12.1,2)--(12.1,4);
        \draw[fill=gray,fill opacity=0.4] (11.6,3.5)--(14.1,3.5)--(13.1,4.5)--(10.6,4.5)--cycle;
        \draw[fill] (12.1,4) circle (1pt);
        \draw[thick,-stealth] (12.1,4)--(12.9,3.7);
        \node at (11.8,4) {$0$};
        \node at (13.1,3.7) {$v$};
        \node at (10.3,3.5) {$(T_xM,g_x)$};
        \node at (12.5,0.6) {\footnotesize $\exp_x(tv)$};
        \draw[gray,fill=gray,fill opacity=0.1] (3,-1.5) to[out=15,in=105] (5,-3.5) to[out=-5,in=-175] (11,-3.5) to[out=105,in=15] (9,-1.5) to[out=-175,in=-5] (3,-1.5);
        \node at (4,-3) {$G$};
        \draw[fill] (5.5,-2) circle (1pt);
        \draw (5.5,-2) to[out=-25,in=115] (6.5,-3);
        \node at (5.2,-2) {$e$};
        \node at (6.5,-3.3) {\footnotesize$\exp(t\xi)$};
        \draw[fill] (9,-2.5) circle (1pt);
        \draw (9,-2.5) to[out=-25,in=115] (10,-3.5);
        \draw (9,-2.5) to[out=155,in=-15] (7.8,-2.1);
        \draw[fill] (9.7,-3) circle (1pt);
        \draw[fill] (8,-2.15) circle (1pt);
        \node at (8.7,-2.7) {$g$};
        \node at (10,-3.9) {\footnotesize $L_g(\exp(t\xi))$};
        \node at (11.5,-3) {\footnotesize $L_g(\exp((1-\theta)\xi))$};
        \node at (9.5,-2) {\footnotesize $L_g(\exp(-\theta\xi))$};
        \draw[dotted] (5.5,-2)--(5.5,-0.5);
        \draw[fill=gray,fill opacity=0.4] (5,-1)--(7.5,-1)--(6.5,0)--(4,0)--cycle;
        \draw[fill] (5.5,-0.5) circle (1pt);
        \draw[thick,-stealth] (5.5,-0.5)--(6.3,-0.8);
        \node at (5.2,-0.5) {$0$};
        \node at (7,-1.3) {\footnotesize $\xi:=T_gL_{g^{-1}}(v)$};
        \node at (4,-0.7) {$\mathfrak{g}$};
        \draw[dotted] (9,-2.5)--(9,-0.8);
        \draw[fill=gray,fill opacity=0.4] (8.5,-1.3)--(11,-1.3)--(10,-0.3)--(7.5,-0.3)--cycle;
        \draw[fill] (9,-0.8) circle (1pt);
        \draw[thick,-stealth] (9,-0.8)--(9.8,-1.1);
        \node at (8.7,-0.8) {$0$};
        \node at (10,-1.1) {$v$};
        \node at (11,-0.5) {$T_gG$};
        \end{tikzpicture}
        \caption{Discretization maps induced by retraction maps on $\mathbb{R}^n$ (top left), on a Riemannian manifold $(M,g)$ (top right) and on a Lie group $G$ (bottom)}
        \label{fig discretization examples}
\end{figure}
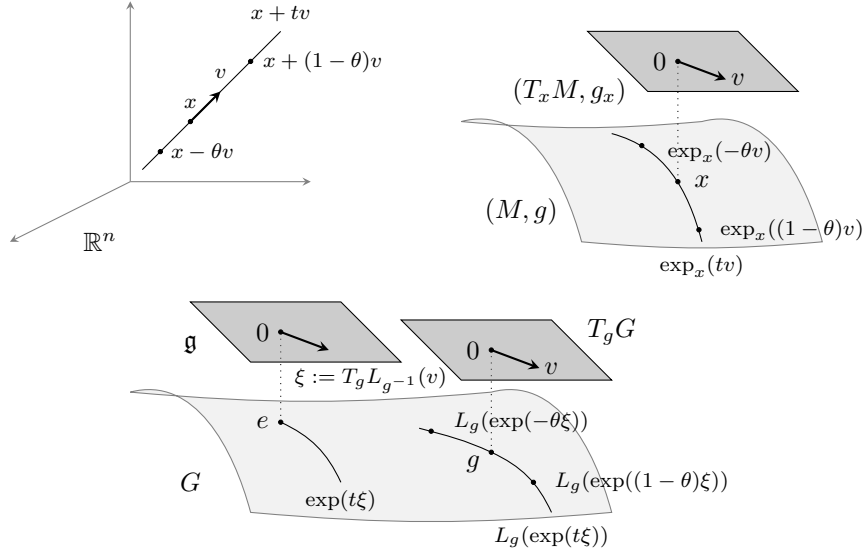
\begin{example} \label{eg disc R^n}
    On $\mathbb{R}^n$, we define a discretization map as 
    \begin{align*}
        \mathcal{D}(x,v):=&(\mathcal{R}(x,-\theta v),\mathcal{R}(x,(1-\theta)v) \\
        =&(x-\theta v,x+(1-\theta)v)
    \end{align*}
    using the retraction $\mathcal{R}(x,v):=x+v$. 
\end{example}
\begin{example} \label{eg disc (M,g)}
    On a Riemannian manifold $(M,g)$, we define a discretization map as
    \begin{align*}
        \mathcal{D}(x,v):=&(\mathcal{R}(x,-\theta v),\mathcal{R}(x,(1-\theta)v))\\
        =&(exp_x(-\theta v),exp_x((1-\theta)v))
    \end{align*}
    using the retraction $\mathcal{R}(x,v):=exp_x(v)$. 
\end{example}
\begin{example} \label{eg disc G}
    On a Lie group $G$, we define a discretization map as
    \begin{align*}
        \mathcal{D}(g,v):=&(\mathcal{R}(g,-\theta v),\mathcal{R}(g,(1-\theta)v))\\
        =&(L_g(exp(-\theta \xi)),L_g(exp((1-\theta)\xi)))
    \end{align*}
    using the retraction $\mathcal{R}(g,v):=L_g(exp(\xi))$. Recall that $\xi:=T_gL_{g^{-1}}(v)$ as shown in figure \ref{fig discretization examples}.
\end{example}
\noindent Just as retraction maps induce discretization maps, trivialized retraction maps naturally induce trivialized discretization maps on Lie groups, which we now proceed to study.
\subsection{Trivialized Discretization Maps}
\label{subsec triv disc maps}
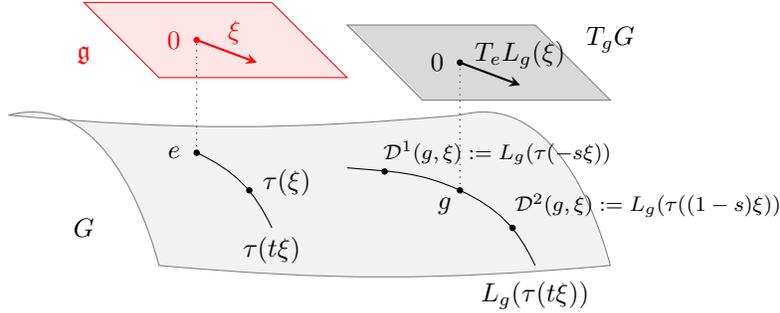
\begin{figure}
    \centering
    \begin{tikzpicture}
        \draw[gray,fill=gray,fill opacity=0.1] (1,3) to[out=15,in=105] (3,1) to[out=-5,in=-175] (9,1) to[out=105,in=15] (7,3) to[out=-175,in=-5] (1,3);
        \node at (2,1.5) {$G$};
        \draw[fill] (3.5,2.5) circle (1pt);
        \draw (3.5,2.5) to[out=-25,in=115] (4.5,1.5);
        \draw[fill] (4.2,2) circle (1pt);
        \node at (3.2,2.5) {$e$};
        \node at (4.5,1.2) {$\tau(t\xi)$};
        \node at (4.7,2.1) {$\tau(\xi)$};
        \draw[fill] (7,2) circle (1pt);
        \draw (7,2) to[out=-25,in=115] (8,1);
        \draw (7,2) to[out=155,in=-5] (5.5,2.3);
        \draw[fill] (7.7,1.5) circle (1pt);
        \draw[fill] (6,2.25) circle (1pt);
        \node at (6.8,1.8) {$g$};
        \node at (8,0.6) {$L_g(\tau(t\xi))$};
        \node at (9.5,1.8) {\footnotesize $\mathcal{D}^2(g,\xi):=L_g(\tau((1-s)\xi))$};
        \node at (7.5,2.5) {\footnotesize $\mathcal{D}^1(g,\xi):=L_g(\tau(-s\xi))$};
        \draw[dotted] (3.5,2.5)--(3.5,4);
        \draw[fill=red,fill opacity=0.1,red] (3,3.5)--(5.5,3.5)--(4.5,4.5)--(2,4.5)--cycle;
        \draw[red,fill] (3.5,4) circle (1pt);
        \draw[red,thick,-stealth] (3.5,4)--(4.3,3.7); 
        \node at (3.2,4) {$\textcolor{red}{0}$};
        \node at (4,4.1) {$\textcolor{red}{\xi}$};
        \node at (2,3.8) {$\textcolor{red}{\mathfrak{g}}$};
        \draw[dotted] (7,2)--(7,3.7);
        \draw[gray,fill=gray,fill opacity=0.3] (6.5,3.2)--(9,3.2)--(8,4.2)--(5.5,4.2)--cycle;
        \draw[fill] (7,3.7) circle (1pt);
        \draw[thick,-stealth] (7,3.7)--(7.8,3.4);
        \node at (6.7,3.7) {$0$};
        \node at (7.8,3.8) {$T_eL_g(\xi)$};
        \node at (9,4) {$T_gG$};
\end{tikzpicture}
    \caption{Left-trivialized discretization map induced by a left-trivialized retraction on $G$}
    \label{fig trivialized discretization}
\end{figure}

We now introduce the notion of trivialized discretization maps on Lie groups.
\begin{definition} \label{def triv disc map}
    A smooth map $\mathcal{D}^L : G \ltimes \mathfrak{g} \to G \times G$ which maps $(g,\xi) \mapsto (\mathcal{D}^1(g,\xi),\mathcal{D}^2(g,\xi))$ and satisfies
    \begin{itemize}
        \item (TD1) $\mathcal{D}^L(g,0)=(\mathcal{D}^1(g,0), \mathcal{D}^2(g,0)) = (g,g)$ and
        \item (TD2) $\left. \frac{d}{dt} \right|_{t=0} \mathcal{D}^2(g,t\xi) - \left. \frac{d}{dt} \right|_{t=0} \mathcal{D}^1(g,t\xi) = T_eL_g(\xi)$
    \end{itemize}
    for all $(g,\xi) \in G \ltimes \mathfrak{g}$ is called a \emph{left-trivialized discretization} map on $G$.
\end{definition}
\noindent A \emph{right-trivialized discretization} map $\mathcal{D}^R : G \ltimes \mathfrak{g} \to G \ltimes G$ is defined analogously with condition (TD2) replaced by $\left. \frac{d}{dt} \right|_{t=0} \mathcal{D}^2(g,t\xi) - \left. \frac{d}{dt} \right|_{t=0} \mathcal{D}^1(g,t\xi) = T_eR_g(\xi)$. Notice that if $\mathcal{D}^1(g,\xi) = g$, then conditions (TD1) and (TD2) together imply that $\mathcal{D}^2 : G \ltimes \mathfrak{g} \to G$ is itself a left-trivialized retraction map on $G$. A left-trivialized discretization is induced by a left-trivialized retraction map as illustrated in figure \ref{fig trivialized discretization}. 
\begin{proposition} \label{prop 2.4}
     A left-trivialized retraction map $\mathcal{R}^L : G \ltimes \mathfrak{g} \to G$ gives rise to a family of left-trivialized discretization maps $\mathcal{D}^L : G \ltimes \mathfrak{g} \to G \times G$ defined as
    \begin{align*}
        \mathcal{D}^L(g,\xi) := (\mathcal{R}^L(g,-s\xi),\mathcal{R}^L(g,(1-s)\xi))
    \end{align*}
    for all $(g,\xi) \in G \ltimes \mathfrak{g}$ and parameterized by $s \in [0,1]$.
\end{proposition}
Just like a discretization map, trivialized discretization maps are not globally invertible, however, they are locally invertible.
\begin{proposition}
    A left-trivialized discretization map $\mathcal{D}^L : G \times \mathfrak{g} \to G \ltimes G$ is locally invertible around $(g,0)$ for all $g \in G$.
\end{proposition}
\section{Numerical Integrators Using Discretization Maps}
\label{sec numerical integrators using disc maps}
Finally, we introduce the central idea for constructing numerical integrators for dynamical systems evolving on a manifold using discretization maps. Let $M$ be a smooth manifold, and let the dynamics be governed by a vector field $X \in \mathfrak{X}(M)$, whose integral curves we wish to approximate numerically. After choosing a discretization map $\mathcal{D} : TM \to M \times M$, we define the numerical integrator implicitly by
\begin{equation}  \label{integrator 1}
    hX(\tau_M(\mathcal{D}^{-1}(x_k,x_{k+1})))=\mathcal{D}^{-1}(x_k,x_{k+1})
\end{equation}
as summarized by the following commutative diagram:
\begin{center}
    \begin{tikzcd}
        M \times M \arrow[r,"\mathcal{D}^{-1}"] \arrow[d,swap,"\mathcal{D}^{-1}"] & TM \arrow[d,"\tau_M"] \\
        TM & M \arrow[l,"X"]
    \end{tikzcd}
\end{center}
Here $\tau_M:TM \to M$ denotes the canonical projection of the tangent bundle, and $h \in \mathbb{R}$ is the step size. The step size $h$ plays a crucial role by scaling the vector field $X$ so that the resulting vector field lies within the local domain of invertibility of
$\mathcal{D}$, ensuring that \eqref{integrator 1} is well defined. Given an initial condition $x_0 \in M$, the implicit difference equation \eqref{integrator 1} can be solved iteratively to generate a discrete trajectory approximating the flow of $X$. We now validate this construction with a simple example. 

Consider the dynamical system given in \eqref{1.1}, whose underlying manifold $M = \mathbb{R}^n$ and the corresponding vector field $X \in \mathfrak{X}(M)$ is given by $X(x) = (x,f(x))$ for every $x \in M$. Let us choose the discretization map introduced in example \ref{eg disc R^n}. Substituting this choice into the general integrator \eqref{integrator 1} yields the implicit one-step scheme
\begin{align*}
    x_{k+1} = x_k + h f((1-\theta) x_k + \theta x_{k+1})
\end{align*}
where $h \in \mathbb{R}$ is the step size and $\theta \in [0,1]$ is a parameter. For $\theta = 0$, this scheme reduces to the explicit Euler method; for $\theta = 1$, it reduces to the implicit Euler method; and for $\theta = 1/2$, it yields the implicit midpoint rule. 
\subsection{Numerical Integrators on Lie Groups using Trivialized Discretization Maps}
\label{subsec integrators on Lie groups}
We now specialize to the case of dynamical systems evolving on Lie groups. The central idea is to use a left- or right-trivialized discretization map, in direct analogy with the discretization map employed in section \ref{sec numerical integrators using disc maps}. In addition, since the underlying space is a Lie group, the vector field describing the dynamics must also be trivialized. We begin by formalizing this notion.
\begin{definition}
    Let $G$ be a Lie group and $X \in \mathfrak{X}(G)$ be a vector field. The \emph{left-trivialization} of $X$ is the map
    \begin{align*}
        \widetilde{X}^L : G \to G \ltimes \mathfrak{g}
    \end{align*}
    obtained by composing $X$ with the left-trivialization of the tangent bundle as shown in the commutative diagram below. The map $\widetilde{X}^L$ is called the \emph{left-trivialized vector field} associated with $X$.
\end{definition}
\begin{center}
    \begin{tikzcd}
        & G \ltimes \mathfrak{g} \\
        G \arrow[r,"X",swap] \arrow[ur,"\widetilde{X}^L"] & TG \arrow[u,"tr^L_{TG}",swap]
    \end{tikzcd}
\end{center}
A \emph{right-trivialization} of a vector field on $G$ is defined analogously using the right-trivialization of $TG$. With these ingredients in place, we can now mirror the general integrator \eqref{integrator 1} in the setting of Lie groups. Given a left-trivialized discretization map $\mathcal{D}^L : G \ltimes\mathfrak{g} \to G \times G$, the numerical integrator is defined implicitly by 
\begin{align} \label{integrator 2 lie groups}
    h \widetilde{X}^L(pr_1((\mathcal{D}^L)^{-1}(g_k,g_{k+1}))) = (\mathcal{D}^L)^{-1}(g_k,g_{k+1})
\end{align}
where $pr_1 : G \ltimes \mathfrak{g} \to G$ denotes the projection onto the first factor and $h \in \mathbb{R}$ is the time step. This construction is summarized by the following commutative diagram:
\begin{center}
    \begin{tikzcd}
        G \times G \arrow[r,"(\mathcal{D}^L)^{-1}"] \arrow[d,"(\mathcal{D}^L)^{-1}",swap] & G \ltimes \mathfrak{g} \arrow[d,"pr_1"] \\
        G \ltimes \mathfrak{g} & G \arrow[l,"\widetilde{X}^L"]
    \end{tikzcd}
\end{center}
As in the general setting, the step size $h$ scales the vector field so that the inverse of the discretization map is evaluated within its domain of local invertibility, ensuring that the scheme is well defined. A noteworthy point is that both choices are independent: the vector field may be trivialized either on the left or on the right, and the discretization map employed may likewise be left- or right-trivialized. These choices can be mixed, leading to distinct but closely related numerical schemes. 
\vspace{10pt}\\
The general integrators presented in \eqref{integrator 1} and their extension to Lie groups in \eqref{integrator 2 lie groups} successfully address the first challenge, namely ensuring that the discrete trajectory remains on the underlying manifold. However, they do not yet preserve any geometric or physical structure. Before modifying these schemes to construct structure-preserving integrators, we therefore review the foundations of classical mechanics, with particular emphasis on its symplectic underpinnings.
\chapter{Lagrangian and Hamiltonian Mechanics}
\label{ch4}
In this chapter, we present an overview of geometric mechanics tailored to the construction of structure-preserving numerical integrators. As illustrated by the preceding examples, mechanical systems—unlike generic dynamical systems—possess intrinsic geometric structures that fundamentally govern their evolution. Classical mechanics admits two complementary formulations: the Lagrangian and the Hamiltonian viewpoints. Our objective is to unify these perspectives by following Tulczyjew’s geometric framework, which will play a central role in the developments that follow. For further in-depth treatments of this rich interplay between geometry and mechanics, we refer the reader to \cite{abraham2008foundations}, \cite{vogtmann1997mathematical}, \cite{marsden2002introduction}, \cite{de2011methods}, and \cite{libermann1987symplectic}.
\section{Lagrangian Mechanics}
\label{Lagrangian Mechanics}
Let $Q$ be a smooth $n$-dimensional manifold and let $TQ$ denote its tangent bundle. Physically, $TQ$ represents the space of configurations and velocities of a mechanical system. Locally, configurations are described by coordinates $(q^i)$, while velocities are denoted by $(\dot{q}^i)$, following the classical notation introduced by Lagrange. 

The dynamics of the system are determined by \emph{Hamilton's principle of stationary action}, which states that the actual trajectory $q(t)$ followed by the system between two fixed times $t_0$ and $t_1$ is a critical point of the action functional
\begin{align*}
    \mathfrak{S}[q] := \int_{t_0}^{t_1} \mathscr{L}(q(t),\dot{q}(t)) \, dt
\end{align*}
where the endpoint $q(t_0)$ and $q(t_1)$ are held fixed. Here, $\mathscr{L} : TQ \to \mathbb{R}$ is a smooth function called the \emph{Lagrangian} of the system.

Applying standard techniques from the calculus of variations, the condition of stationarity $\delta \mathfrak{S} = 0$ yields the \emph{Euler-Lagrange} equations
\begin{align} \label{Euler-Lagrange eq}
     \frac{d}{dt} \left(\frac{\partial \mathscr{L}}{\partial \dot{q}^i} \right)=\frac{\partial \mathscr{L}}{\partial q^i} \quad \text{for all} \quad i \in \{1,\dots,n\}
\end{align}
which govern the equations of motion of the system. For \emph{simple mechanical systems}, the Lagrangian takes the form
\begin{align*}
    \mathscr{L}(q,\dot{q}) = T(q,\dot{q}) - V(q)
\end{align*}
where $T : TQ \to \mathbb{R}$ denotes the kinetic energy and $V : Q \to \mathbb{R}$ denotes the potential energy. In this case, the Euler–Lagrange equations reduce to Newton’s equations of motion, thereby recovering the classical formulation of mechanics from the variational principle.
\section{Hamiltonian Mechanics}
\label{Hamiltonian Mechanics}
Let $Q$ be a smooth $n$-dimensional manifold and let $T^*Q$ denote its cotangent bundle. Physically, $T^*Q$ is interpreted as the space of all positions and momenta of the mechanical system and is commonly referred to as the \emph{phase space}. Locally, positions are denoted by coordinates $(q^i)$, while the conjugate momenta are denoted by $(p_i)$, following the classical convention introduced by Hamilton.

The equations of motion are governed by
\begin{align} \label{Hamilton's eq}
    \dot{q}^i = \frac{\partial \mathscr{H}}{\partial p_i} \quad \text{and} \quad \dot{p}_i = -\frac{\partial \mathscr{H}}{\partial q^i} \quad \text{for all} \quad i \in \{1,\dots,n\}
\end{align}
where $\mathscr{H}: T^*Q \to \mathbb{R}$ is the \emph{Hamiltonian} of the system. These equations arise naturally from the canonical symplectic structure on $T^*Q$ and are known as \emph{Hamilton's equations}. Usually, the Hamiltonian coincides with the total energy of the system, namely the sum of kinetic and potential energies. In this case, Hamilton’s equations are equivalent to Newton’s equations of motion and provide an alternative, geometrically natural formulation of classical mechanics.
\section{Unification of Lagrangian and Hamiltonian Viewpoints}
\label{unification of L and H}
The unification of the Lagrangian and Hamiltonian viewpoints is achieved by considering the double tangent bundle $TT^*Q$. This space admits a natural physical interpretation as the space of positions, momenta, velocities, and forces, thereby encoding all fundamental variables of classical mechanics. To make this precise, one studies the canonical morphisms between the second-order bundles
\begin{align*}
    TTQ, \quad TT^*Q, \quad T^*TQ \quad \text{and} \quad T^*T^*Q,
\end{align*}
whose interrelations form the geometric backbone of Tulczyjew’s framework, see \cite{tulczyjew1977legendre, grabowska2012tulczyjew, zajkac2016tulczyjew}.
\paragraph{$\kappa_Q$}
Unlike the tangent bundle $TQ$, which carries a canonical vector bundle structure $\tau_Q : TQ \to Q$, the double tangent bundle $TTQ$ admits two distinct vector bundle structures, namely $\tau_{TQ} : TTQ \to TQ$ and $T\tau_Q : TTQ \to TQ$. The vector bundle isomorphism relating these two structures is denoted by $\kappa_Q : TTQ \to TTQ$, and is illustrated in the commutative diagram below. 
\begin{center}
    \begin{tikzcd}
        TTQ \arrow[rr,shift left,"\kappa_Q"] \arrow[dr,"\tau_{TQ}",swap] & & TTQ \arrow[dl,"T\tau_Q"] \arrow[ll,shift left,"\kappa_Q"] \\
        & TQ \arrow[d,"\tau_Q"] & \\
        & Q &
    \end{tikzcd}
\end{center}
To give a concrete description of $\kappa_Q$, let $c(s,t)$ be a smooth two-parameter family of curves on $Q$, with $s,t \in ]\varepsilon,\varepsilon[$ for some $\varepsilon > 0$. Fixing $s$ and differentiating with respect to $t$ at $t=0$ yields a curve in $TQ$,
\begin{align*}
    \left(c(s,0), \left. \frac{d}{dt} \right|_{t=0} c(s,t)\right).
\end{align*}
Differentiating this curve with respect to $s$ at $s=0$ produces an element of $TTQ$, locally represented as
\begin{align*}
    (q,\dot{q},\delta q, \delta \dot{q}) = \left(c(0,0),\left. \frac{d}{dt} \right|_{t=0} c(0,t),\left. \frac{d}{ds} \right|_{s=0} c(s,0),\left. \frac{d^2}{ds\,dt} \right|_{s=t=0} c(s,t)\right).
\end{align*}
The map $\kappa_Q$ acts by interchanging the order of differentiation, corresponding to swapping the parameters $s$ and $t$. In local coordinates, this amounts to
\begin{align} 
    \kappa_Q(q,\dot{q},\delta q, \delta \dot{q}) := (q,\delta q,\dot{q}, \delta \dot{q}). \label{kappa_Q local}
\end{align}
Since $\kappa_Q^2 = id_{TTQ}$, the map $\kappa_Q$ is called the \emph{canonical involution} or \emph{canonical flip} on $TTQ$.
\paragraph{$\alpha_Q$}
The map $\alpha_Q : TT^*Q \to T^*TQ$ is a vector bundle isomorphism defined as the dual of the canonical involution $\kappa_Q : TTQ \to TTQ$. This relationship is summarized by the commutative diagram below.  
\begin{center}
    \begin{tikzcd}
        TT^*Q \arrow[dr,"T\pi_Q"] \arrow[dd,"\alpha_Q",swap] & & TTQ \arrow[dl,"\tau_{TQ}",swap] \\
        & TQ &  \\
        T^*TQ \arrow[ur,"\pi_{TQ}",swap] & & TTQ \arrow[ul,"T\tau_Q"] \arrow[uu,"\kappa_Q",swap]        
    \end{tikzcd}
\end{center}
To make this construction precise, recall the natural duality pairings $(T^*TQ, TTQ)$ and $(T^*Q, TQ)$, both denoted by $\langle \cdot \,, \cdot \rangle$. This induces a pairing between $TT^*Q$ and $TTQ$, defined by
\begin{align*}
    \langle\langle \Lambda, V \rangle\rangle := \left. \frac{d}{dt} \right|_{t=0} \langle \lambda(t), v(t) \rangle
\end{align*}
where $\lambda(t)$ and $v(t)$ are curves on $T^*Q$ and $TQ$ with $\dot{\lambda}(0)=\Lambda \in TT^*Q$ and $\dot{v}(0)=V \in TTQ$. The map $\alpha_Q$ is then defined implicitly by 
\begin{align*}
    \langle \alpha_Q(\Lambda),V \rangle := \langle\langle \Lambda,\kappa_Q(V) \rangle\rangle
\end{align*}
for all $\Lambda \in TT^*Q$ and $V \in TTQ$. In local coordinates, an element of $TT^*Q$ is written as $(q,p,\dot{q},\dot{p})$, and the map $\alpha_Q$ takes the explicit form 
\begin{align} 
    \alpha_Q (q,p,\dot{q},\dot{p}) := (q, \dot{q}, \dot{p},p), \label{alpha_Q local}
\end{align}
where the target is identified with coordinates on $T^*TQ$. Apart from this, the canonical symplectic form $\omega_Q$ on $T^*Q$ induces a natural vector bundle isomorphism between $TT^*Q$ and $T^*T^*Q$. This identification allows one to pull back the canonical symplectic form $\omega_{T^*Q}$ on $T^*T^*Q$ to $TT^*Q$, thereby equipping $TT^*Q$ with a symplectic structure known as the \emph{tangent-lifted symplectic} structure. With respect to this symplectic structure, the map $\alpha_Q : TT^*Q \to T^*TQ$ is a symplectomorphism.
\paragraph{$\beta_Q$}
The map $\beta_Q : TT^*Q \to T^*T^*Q$ is a vector bundle isomorphism induced by the canonical symplectic form $\omega_Q$ on $T^*Q$. This is summarized in the commutative diagram below. When $TT^*Q$ is equipped with the tangent-lifted symplectic structure and $T^*T^*Q$ with its canonical symplectic structure, $\beta_Q$ is a symplectomorphism. 
\begin{center}
    \begin{tikzcd}
        T^*T^*Q \arrow[dr,"\pi_{T^*Q}",swap] & & TT^*Q \arrow[dl,"\tau_{T^*Q}"] \arrow[ll,"\beta_Q",swap] \\
        & T^*Q \arrow[d,"\pi_Q"] & \\
        & Q &
    \end{tikzcd}
\end{center}
In local coordinates, if an element of $TT^*Q$ is written as $(q,p,\dot{q},\dot{p})$, then
\begin{align} 
    \beta_Q(q,p,\dot{q},\dot{p}) := (q,p,\dot{p},-\dot{q}), \label{beta_Q local}
\end{align}
where the target is identified with coordinates on $T^*T^*Q$.
\vspace{10pt}\\
We now return to the unification of the Lagrangian and Hamiltonian formalisms.
\subsection{Lagrangian Picture} 
Let $\mathscr{L} :TQ \to \mathbb{R}$ be a Lagrangian. Its differential defines a Lagrangian submanifold $d\mathscr{L} (TQ) \subset T^*TQ$ where $T^*TQ$ carries its canonical symplectic form $\omega_{TQ}$. Using the symplectomorphism $\alpha_Q : TT^*Q \to T^*TQ$, we pull this submanifold back to $TT^*Q$, which is equipped with the tangent-lifted symplectic structure. Thus the Lagrangian induces the Lagrangian submanifold $\alpha_Q^{-1}(d\mathscr{L}(TQ)) \subset TT^*Q$. 
\begin{center}
     \begin{tikzcd}[row sep = large]
         T^*TQ \arrow[rd,"\pi_{TQ}",shift left] & & TT^*Q \arrow[ll,"\alpha_Q",swap] \arrow[dl,"T\pi_{Q}",shift right,swap] \\
        & TQ \arrow[lu,"\textcolor{red}{d\mathscr{L}}",shift left,red]  \arrow[d,"\tau_Q"] \arrow[ld,"\textcolor{red}{\mathscr{L}}",red] & \\
        \mathbb{R} & Q &
    \end{tikzcd}
\end{center}
In local coordinates $(q^i,\dot{q}^i)$ on $TQ$, $d\mathscr{L}(q^i,\dot{q}^i)=(q^i,\dot{q}^i,\tfrac{\partial\mathscr{L}}{\partial q^i},\tfrac{\partial\mathscr{L}}{\partial\dot{q}^i})$ and hence $\alpha_Q^{-1}(d\mathscr{L}(q^i,\dot{q}^i))=(q^i,\frac{\partial\mathscr{L}}{\partial\dot{q}^i},\dot{q}^i,\frac{\partial\mathscr{L}}{\partial q^i})$. Interpreting the second component as a velocity yields the Euler–Lagrange equations 
\begin{align*}
    \frac{d}{dt}\left(\frac{\partial\mathscr{L}}{\partial\dot{q}^i}\right)=\frac{\partial\mathscr{L}}{\partial q^i} \quad \text{for all} \quad i \in \{1,\dots,n\}.
\end{align*}
Hence, the variational principle is equivalent to specifying a Lagrangian subamanifold in the appropriate symplectic manifold.
\subsection{Hamiltonian Picture}
Let $\mathscr{H} : T^*Q \to \mathbb{R}$ be a Hamiltonian. Its differential defines a Lagrangian submanifold $d\mathscr{H}(T^*Q) \subset T^*T^*Q$ where $T^*T^*Q$ carries its canonical symplectic form $\omega_{T^*Q}$. Using the symplectomorphism $\beta_Q : TT^*Q \to T^*T^*Q$, we pull this submanifold back to $TT^*Q$, which is equipped with the tangent-lifted symplectic structure. Thus the Hamiltonian induces the Lagrangian submanifold $\beta_Q^{-1}(d\mathscr{H}(T^*Q)) \subset TT^*Q$. This submanifold coincides with the image of the Hamiltonian vector field $X_\mathscr{H} : T^*Q \to TT^*Q$ defined by the canonical symplectic structure on $T^*Q$. 
\begin{center}
    \begin{tikzcd}[row sep = large]
        T^*T^*Q \arrow[rd,"\pi_{T^*Q}",shift left] & & TT^*Q \arrow[ll,"\beta_Q",swap] \arrow[dl,"\tau_{T^*Q}",shift right,swap] \\
        & T^*Q \arrow[lu,"\textcolor{blue}{d\mathscr{H}}",shift left,blue] \arrow[ru,"\textcolor{blue}{X_\mathscr{H}}",shift right,blue,swap] \arrow[d,"\pi_Q"] \arrow[ld,"\textcolor{blue}{\mathscr{H}}",blue] & \\
        \mathbb{R} & Q &
    \end{tikzcd}
\end{center}
In local coordinates $(q^i,p_i)$ on $T^*Q$, $d\mathscr{H}(q^i,p_i) = (q^i,p_i, \frac{\partial \mathscr{H}}{\partial q^i},\frac{\partial \mathscr{H}}{\partial p_i})$ and hence $\beta_Q^{-1}(d\mathscr{H}(q^i,p_i))=(q^i,p_i;\frac{\partial\mathscr{H}}{\partial p_i},-\frac{\partial\mathscr{H}}{\partial q^i})$. Interpreting the second component as a velocity yields the Hamilton's equations
\begin{align*}
    \dot{q}^i = \frac{\partial\mathscr{H}}{\partial p_i} \quad \text{and} \quad \dot{p}_i = -\frac{\partial\mathscr{H}}{\partial q^i} \quad \text{for all} \quad i \in \{1,\dots,n\}.
\end{align*}
\subsection{Legendre Transform}
In particular, if the Legendre transform $\mathbb{F}\mathscr{L}: TQ \to T^*Q$ is a diffeomorphism, i.e. the Lagrangian $\mathscr{L}$ is hyperregular, then both the Hamiltonian $\mathscr{H} : T^*Q \to \mathbb{R}$ and the Hamiltonian vector field $X_\mathscr{H} : T^*Q \to TT^*Q$ can be pulled back to $TQ$ via $\mathbb{F}\mathscr{L}$, as summarized in the diagram below.
\begin{center}
    \begin{tikzcd}
        TT^*Q \arrow[dd,"\tau_{T^*Q}",shift left] & & TTQ \arrow[ll,"\textcolor{red}{T(\mathbb{F}\mathscr{L})}",red,swap] \arrow[dd,"\tau_{TQ}",shift right,swap] \arrow[rr,"\kappa_Q",shift left] & & TTQ \arrow[ddll,"T\tau_Q"] \arrow[ll,"\kappa_Q",shift left] \\
        & \mathbb{R} & & \\
        T^*Q \arrow[uu,"\textcolor{blue}{X_\mathscr{H}}",shift left,blue] \arrow[ur,"\textcolor{blue}{\mathscr{H}}",blue] \arrow[dr,"\pi_Q",swap] & & TQ \arrow[uu,"\textcolor{orange}{X_\mathscr{E}}",shift right, orange,swap] \arrow[ll,"\textcolor{red}{\mathbb{F}\mathscr{L}}",red,swap] \arrow[ul,"\textcolor{orange}{\mathscr{E}}",orange,swap] \arrow[dl,"\tau_Q"] \arrow[dr,"\textcolor{red}{\mathscr{L}}",red] & & \\
        & Q & & \mathbb{R} &
    \end{tikzcd}
\end{center}
The resulting vector field $X_\mathscr{E}: TQ \to TTQ$ is Hamiltonian with respect to the pulled back symplectic structure and the energy function $\mathscr{E}:=(\mathbb{F}\mathscr{L})^*\mathscr{H}$. Moreover, the image $X_\mathscr{E}(TQ) \subset TTQ$ is invariant under the canonical involution $\kappa_Q: TTQ \to TTQ$, which implies that $X_\mathscr{E}$ is a second-order vector field or a semispray on $Q$, see \cite{de1989methods}. In local coordinates $(q^i,\dot{q}^i)$ on $TQ$, the pulled back Hamiltonian vector field $X_\mathscr{E}$ also known as the Euler-Lagrange vector field, takes the explicit form
\begin{align} \label{Lagrangian semi-spray eq}
    \ddot{q}^i = \left(\frac{\partial^2 \mathscr{L}}{\partial \dot{q}^i \partial \dot{q}^j}\right)^{-1} \left(\frac{\partial \mathscr{L}}{\partial q^j} - \frac{\partial^2 \mathscr{L}}{\partial q^k \partial \dot{q}^j}\dot{q}^k\right) \quad \text{for all} \quad i,j,k \in \{1,\dots,n\}. 
\end{align}

Consequently, the Lagrangian submanifolds of $TT^*Q$ induced by the Lagrangian and Hamiltonian formalisms coincide if and only if the Lagrangian is hyperregular. For regular or degenerate Lagrangians, the Euler–Lagrange and Hamilton equations define distinct Lagrangian submanifolds of $TT^*Q$, and are therefore not equivalent, despite sharing a common geometric origin.
\section{Lagrangian Mechanics on Lie Groups}
Consider a mechanical system evolving on a Lie group $G$. Its Lagrangian $\mathscr{L} : TG \to \mathbb{R}$ can be left-trivialized as $\widetilde{\mathscr{L}}^L \circ tr_{TG}^L := \mathscr{L}$ as illustrated in the diagram below.
\begin{center}
    \begin{tikzcd}
        G \ltimes \mathfrak{g} \arrow[dr,"\widetilde{\mathscr{L}}^L"] & \\
        TG \arrow[u,"tr_{TG}^L"] \arrow[r,"\mathscr{L}",swap] & \mathbb{R}
    \end{tikzcd}
\end{center}
The \emph{left-trivialized Euler-Lagrange} equations are then given by
\begin{align} 
    \frac{d}{dt} \left( \frac{\delta \widetilde{\mathscr{L}}^L}{\delta \xi} \right) &=  ad_\xi^* \left( \frac{\delta \widetilde{\mathscr{L}}^L}{\delta \xi}  \right) + T_e^*L_g \left( \frac{\delta \widetilde{\mathscr{L}}^L}{\delta g} \right) \quad \text{and} \label{left-trivialized Euler-Lagrange eq1} \\
    \dot{g} &= T_eL_g(\xi). \label{left-trivialized Euler-Lagrange eq2}
\end{align}
Here $\tfrac{\delta \widetilde{\mathscr{L}}^L}{\delta \xi} \in \mathfrak{g}^*$ and $\tfrac{\delta \widetilde{\mathscr{L}}^L}{\delta g} \in T_g^*G$ denote the variational derivatives of the left-trivialized Lagrangian $\widetilde{\mathscr{L}}^L : G \ltimes \mathfrak{g} \to \mathbb{R}$. The right-trivialized Euler-Lagrange equations are defined analogously.

In the special case where the Lagrangian $\mathscr{L} : TG \to \mathbb{R}$ is left-invariant, i.e. $\mathscr{L} \circ TL_g = \mathscr{L}$ for all $g \in G$, it reduces to the \emph{left-reduced Lagrangian} $l : \mathfrak{g} \to \mathbb{R}$. The reduced dynamics are governed by
\begin{align} 
    \frac{d}{dt} \left( \frac{\delta l}{\delta \xi} \right) &= ad_\xi^* \left( \frac{\delta l}{\delta \xi} \right) \quad \text{and} \label{left Euler-Poincare eq1} \\
    \dot{g} &= T_eL_g(\xi) \label{left Euler-Poincare eq2}
\end{align}
where $\tfrac{\delta l}{\delta \xi} \in \mathfrak{g}^*$ denotes the variational derivative of $l$, defined by
\begin{align*}
    \left\langle \frac{\delta l}{\delta \xi}, \delta \xi \right\rangle := \left.\frac{d}{ds}\right|_{s=0} l(\xi +  s \delta\xi)
\end{align*}
for all variations $\delta \xi \in \mathfrak{g}$. These equations are known as the \emph{left Euler-Poincar{\'e}} equations. The reconstruction equation $\dot{g} = T_eL_g(\xi)$ recovers the dynamics on $G$ from the reduced variables, as illustrated in the figure below.
\begin{center}
\begin{tikzpicture}[scale = 0.5]
\draw[] (1,1) to[out=30, in=150] (11,1) to[out=75, in= 210] (13,5) to[out=150, in=30] (3,5) to[out=210, in=75] (1,1);
\node at (12,2) {$G$};
\draw[fill] (4,4) circle (2pt);
\node at (4,3.5) { $e$};
\draw[fill=red,red] (4,6.5) circle (2pt); 
\draw[fill=red, fill opacity=0.1, red] (2,4.5)--(6,5.5)--(7,8.5)--(3,7.5)--cycle;
\draw[dotted] (4,4)--(4,6.5);
\node at (2,6.5) { \textcolor{red}{$\mathfrak{g}$}};
\draw[fill=blue,blue] (4,9.5) circle (2pt); 
\draw[fill=blue, fill opacity=0.1, blue] (2,7.5)--(6,8.5)--(7,11.5)--(3,10.5)--cycle;
\node at (2,9.5) { \textcolor{blue}{$\mathfrak{g}^*$}};
\draw[dotted] (4,6.5)--(4,9.5);
\draw[thick] (7,6) to[out=-15 ,in=105] (10,2);
\draw[fill] (9.3,4) circle (2pt);
\node at (8.9,3.7) { $g$};
\draw[fill=gray,gray] (9.3,7.5) circle (2pt);
\draw[fill=gray, fill opacity=0.1, gray] (7.3,6.5)--(11.3,5.5)--(12.8,8)--(8.8,9)--cycle;
\node at (11,9) {\textcolor{gray}{$T_gG$}};
\draw[dotted] (9.3,4)--(9.3,7.5);
%
\draw[blue,thick] (4,10.5) to[out=0, in=120] (5.8,9.3);
\draw[ultra thick,blue,-stealth] (4,9.5)--(5,10.2);
\node at (4.7,9.2) {\textcolor{blue}{$\frac{\delta l}{\delta \xi}$}};
\draw[ultra thick,blue,-stealth] (5,10.2)--(5.8,9.7);
\node at (9,10.4) {\textcolor{blue}{$\frac{d}{dt}\left(\frac{\delta l}{\delta \xi}\right)=ad^*_\xi\left(\frac{\delta l}{\delta \xi}\right)$}};
\draw[red,thick] (4,7.5) to[out=0, in=120] (5.8,6.3);
\draw[ultra thick,red,-stealth] (4,6.5)--(5,7.2);
\node at (4.7,6.3) {\textcolor{red}{$\xi$}};
\draw[gray,ultra thick,-stealth] (9.3,7.5)--(10,6.5);
\node at (12,7.4) {\textcolor{gray}{$\dot{g}=T_eL_g(\xi)$}};
%
\end{tikzpicture}
\end{center}
The right invariant case follows analogously. 
\section{Hamiltonian Mechanics on Lie Groups}
Again consider a mechanical system evolving on a Lie group $G$. Its Hamiltonian $\mathscr{H} : T^*G \to \mathbb{R}$ can be left-trivialized as $\widetilde{\mathscr{H}}^L \circ tr_{T^*G}^L := \mathscr{H}$ as illustrated in the diagram below.
\begin{center}
    \begin{tikzcd}
        G \ltimes \mathfrak{g}^* \arrow[dr,"\widetilde{\mathscr{H}}^L"] & \\
        T^*G \arrow[u,"tr_{T^*G}^L"] \arrow[r,"\mathscr{H}",swap] & \mathbb{R}
    \end{tikzcd}
\end{center}
The \emph{left-trivialized Hamilton's} equations are then given by
\begin{align} 
    \dot{\mu} &=  ad_{\tfrac{\delta \widetilde{\mathscr{H}}^L}{\delta \mu}}^* (\mu) - T_e^*L_g \left( \frac{\delta \widetilde{\mathscr{H}}^L}{\delta g} \right) \quad \text{and} \label{left-trivialized Hamiltons eq1} \\
    \dot{g} &= T_eL_g\left(\frac{\delta \widetilde{\mathscr{H}}^L}{\delta \mu}\right). \label{left-trivialized Hamiltons eq2}
\end{align}
Here $\tfrac{\delta \widetilde{\mathscr{H}}^L}{\delta \mu} \in \mathfrak{g}$ and $\tfrac{\delta \widetilde{\mathscr{H}}^L}{\delta g} \in T_g^*G$ denote the variational derivatives of the left-trivialized Hamiltonian $\widetilde{\mathscr{H}}^L : G \ltimes \mathfrak{g}^* \to \mathbb{R}$. The right-trivialized Hamilton's equations are defined analogously. 

In the special case where the Hamiltonian $\mathscr{H} : T^*G \to \mathbb{R}$ is left-invariant, i.e. $\mathscr{H} \circ T^*L_g = \mathscr{H}$ for all $g \in G$, it reduces to the \emph{left-reduced Hamiltonian} $h : \mathfrak{g}^* \to \mathbb{R}$. The reduced dynamics are governed by
\begin{align} 
    \dot{\mu} &= ad_{\frac{\delta h}{\delta \mu}}^* (\mu) \quad \text{and} \label{Lie-Poisson eq1} \\
    \dot{g} &= T_eL_g \left(\frac{\delta h}{\delta \mu}\right) \label{Lie-Poisson eq2}
\end{align}
where $\tfrac{\delta h}{\delta \mu} \in \mathfrak{g}$ denotes the variational derivative of $\mathscr{h}$, defined by
\begin{align*}
    \left\langle \frac{\delta h}{\delta \mu}, \delta \mu \right\rangle := \left.\frac{d}{ds}\right|_{s=0} h(\mu +  s \delta\mu)
\end{align*}
for all variations $\delta \mu \in \mathfrak{g}^*$. These equations are known as the \emph{left Lie-Poisson} equations. The reconstruction equation $\dot{g} = T_eL_g(\tfrac{\delta h}{\delta \mu})$ recovers the dynamics on $G$ from the reduced variables, as illustrated in the figure below. 
\begin{center}
\begin{tikzpicture}[scale = 0.5]
\draw[] (1,1) to[out=30, in=150] (11,1) to[out=75, in= 210] (13,5) to[out=150, in=30] (3,5) to[out=210, in=75] (1,1);
\node at (12,2) {$G$};
\draw[fill] (4,4) circle (2pt);
\node at (4,3.5) { $e$};
\draw[fill=red,red] (4,6.5) circle (2pt); 
\draw[fill=red, fill opacity=0.1, red] (2,4.5)--(6,5.5)--(7,8.5)--(3,7.5)--cycle;
\draw[dotted] (4,4)--(4,6.5);
\node at (2,6.5) { \textcolor{red}{$\mathfrak{g}$}};
\draw[fill=blue,blue] (4,9.5) circle (2pt); 
\draw[fill=blue, fill opacity=0.1, blue] (2,7.5)--(6,8.5)--(7,11.5)--(3,10.5)--cycle;
\node at (2,9.5) { \textcolor{blue}{$\mathfrak{g}^*$}};
\draw[dotted] (4,6.5)--(4,9.5);
\draw[thick] (7,6) to[out=-15 ,in=105] (10,2);
\draw[fill] (9.3,4) circle (2pt);
\node at (8.9,3.7) { $g$};
\draw[fill=gray,gray] (9.3,7.5) circle (2pt);
\draw[fill=gray, fill opacity=0.1, gray] (7.3,6.5)--(11.3,5.5)--(12.8,8)--(8.8,9)--cycle;
\node at (11,9) {\textcolor{gray}{$T_gG$}};
\draw[dotted] (9.3,4)--(9.3,7.5);
%
\draw[blue,thick] (4,10.5) to[out=0, in=120] (5.8,9.3);
\draw[ultra thick,blue,-stealth] (4,9.5)--(5,10.2);
\node at (4.7,9.4) {\textcolor{blue}{$\mu$}};
\draw[ultra thick,blue,-stealth] (5,10.2)--(5.8,9.7);
\node at (7.5,10.4) {\textcolor{blue}{$\dot{\mu}=ad^*_{\frac{\delta h}{\delta \mu}}(\mu)$}};
\draw[red,thick] (4,7.5) to[out=0, in=120] (5.8,6.3);
\draw[ultra thick,red,-stealth] (4,6.5)--(5,7.2);
\node at (4.7,6.1) {\textcolor{red}{$\frac{\delta h}{\delta \mu}$}};
\draw[ultra thick,-stealth,gray] (9.3,7.5)--(10,6.5);
\node at (12.5,7.3) {\textcolor{gray}{$\dot{g}=T_eL_g \left(\frac{\delta h}{\delta \mu}\right)$}};
%
\end{tikzpicture}
\end{center}
The right invariant case follows analogously. 
\section{Unification of Lagrangian and Hamiltonian Viewpoints \\ on Lie Groups}
The key to unifying the Lagrangian and Hamiltonian viewpoints on Lie groups, along the lines of Section \ref{unification of L and H}, is to trivialize all the bundles involved in a manner that preserves both the Lie group structure and the underlying symplectic structure, see \cite{esen2014tulczyjew} for details. We restrict ourselves to left-trivializations in what follows, as the construction of right-trivializations proceeds in an entirely analogous way..
\paragraph{$TTG$}
We now trivialize the double tangent bundle of a Lie group $G$. Consider the following commutative diagram.
\begin{center}
  \begin{tikzcd}[column sep=large, row sep=large]
    TTG \arrow[r,"Ttr^L_{TG}"] \arrow[d,"\tau_{TG}",swap] & T(G \ltimes \mathfrak{g}) \arrow[r,"tr^L_{T(G \ltimes \mathfrak{g})}"] \arrow[d,"\tau_{G \ltimes \mathfrak{g}}"] & (G \ltimes \mathfrak{g}) \ltimes (\mathfrak{g} \ltimes \mathfrak{g}) \arrow[dl,"pr_{12}"] \\
    TG \arrow[r,"tr^L_{TG}",swap] & G \ltimes \mathfrak{g} & 
\end{tikzcd}  
\end{center}
Hence, $TTG$ can be identified with $(G \ltimes \mathfrak{g}) \ltimes (\mathfrak{g} \ltimes \mathfrak{g})$ via the map $tr^L_{T(G \ltimes \mathfrak{g})} \circ Ttr^L_{TG}$. 
This identification is a Lie group isomorphism. Indeed, both $tr^L_{TG}$ and $tr^L_{T(G \ltimes \mathfrak{g})}$ are Lie group isomorphisms, and the tangent lift of a Lie group isomorphism is again a Lie group isomorphism with respect to the canonical Lie group structure on the tangent bundle.
\paragraph{$TT^*G$}
We now trivialize the tangent bundle of the cotangent bundle of a Lie group $G$. Consider the following commutative diagram.
\begin{center}
  \begin{tikzcd}[column sep=large, row sep=large]
    TT^*G \arrow[r,"Ttr^L_{T^*G}"] \arrow[d,"\tau_{T^*G}",swap] & T(G \ltimes \mathfrak{g}^*) \arrow[r,"tr^L_{T(G \ltimes \mathfrak{g}^*)}"] \arrow[d,"\tau_{G \ltimes \mathfrak{g}^*}"] & (G \ltimes \mathfrak{g}^*) \ltimes (\mathfrak{g} \ltimes \mathfrak{g}^*) \arrow[dl,"pr_{12}"] \\
    T^*G \arrow[r,"tr^L_{T^*G}",swap] & G \ltimes \mathfrak{g}^* & 
\end{tikzcd}  
\end{center}
Hence, $TT^*G$ can be identified with $(G \ltimes \mathfrak{g}^*) \ltimes (\mathfrak{g} \ltimes \mathfrak{g}^*)$ via the map $tr^L_{T(G \ltimes \mathfrak{g}^*)} \circ Ttr^L_{T^*G}$. 
This identification is a Lie group isomorphism. Indeed, both $tr^L_{T^*G}$ and $tr^L_{T(G \ltimes \mathfrak{g}^*)}$ are Lie group isomorphisms, and the tangent lift of a Lie group isomorphism is again a Lie group isomorphism with respect to the canonical Lie group structure on the tangent bundle.

Moreover, this identification is a symplectomorphism. Both $tr^L_{T^*G}$ and $tr^L_{T(G \ltimes \mathfrak{g}^*)}$ are symplectomorphisms, and the tangent lift of a symplectomorphism is again a symplectomorphism with respect to the tangent-lifted symplectic structure on the tangent bundle.
\paragraph{$T^*TG$}
We now trivialize the cotangent bundle of the tangent bundle of a Lie group $G$. Consider the following commutative diagram.
\begin{center}
  \begin{tikzcd}[column sep=large, row sep=large]
    T^*TG \arrow[d,"\pi_{TG}",swap] & T^*(G \ltimes \mathfrak{g}) \arrow[l,"T^*tr^L_{TG}",swap] \arrow[r,"tr^L_{T^*(G \ltimes \mathfrak{g})}"] \arrow[d,"\pi_{G \ltimes \mathfrak{g}}"] & (G \ltimes \mathfrak{g}) \ltimes (\mathfrak{g}^* \ltimes \mathfrak{g}^*) \arrow[dl,"pr_{12}"] \\
    TG \arrow[r,"tr^L_{TG}",swap] & G \ltimes \mathfrak{g} & 
\end{tikzcd}  
\end{center}
Hence, $T^*TG$ can be identified with $(G \ltimes \mathfrak{g}) \ltimes (\mathfrak{g}^* \ltimes \mathfrak{g}^*)$ via the map $tr^L_{T^*(G \ltimes \mathfrak{g})} \circ (T^*tr^L_{TG})^{-1}$. 
This identification is a Lie group isomorphism. Indeed, both $tr^L_{TG}$ and $tr^L_{T^*(G \ltimes \mathfrak{g})}$ are Lie group isomorphisms, and the cotangent lift of a Lie group isomorphism is again a Lie group isomorphism with respect to the canonical Lie group structure on the cotangent bundle.

Moreover, this identification is a symplectomorphism. Both $T^*tr^L_{TG}$ and $tr^L_{T^*(G \ltimes \mathfrak{g})}$ are symplectomorphisms.
\paragraph{$T^*T^*G$}
We now trivialize the double cotangent bundle of a Lie group $G$. Consider the following commutative diagram.
\begin{center}
  \begin{tikzcd}[column sep=large, row sep=large]
    T^*T^*G \arrow[d,"\pi_{T^*G}",swap] & T^*(G \ltimes \mathfrak{g}^*) \arrow[l,"T^*tr^L_{T^*G}",swap] \arrow[r,"tr^L_{T^*(G \ltimes \mathfrak{g}^*)}"] \arrow[d,"\pi_{G \ltimes \mathfrak{g}^*}"] & (G \ltimes \mathfrak{g}^*) \ltimes (\mathfrak{g}^* \ltimes \mathfrak{g}) \arrow[dl,"pr_{12}"] \\
    T^*G \arrow[r,"tr^L_{T^*G}",swap] & G \ltimes \mathfrak{g}^* & 
\end{tikzcd}  
\end{center}
Hence, $T^*T^*G$ can be identified with $(G \ltimes \mathfrak{g}^*) \ltimes (\mathfrak{g}^* \ltimes \mathfrak{g})$ via the map $tr^L_{T^*(G \ltimes \mathfrak{g}^*)} \circ (T^*tr^L_{T^*G})^{-1}$. 
This identification is a Lie group isomorphism. Indeed, both $tr^L_{T^*G}$ and $tr^L_{T^*(G \ltimes \mathfrak{g}^*)}$ are Lie group isomorphisms, and the cottangent lift of a Lie group isomorphism is again a Lie group isomorphism with respect to the canonical Lie group structure on the cotangent bundle.

Moreover, this identification is a symplectomorphism. Both $T^*tr^L_{T^*G}$ and $tr^L_{T^*(G \ltimes \mathfrak{g}^*)}$ are symplectomorphisms.
\vspace{10pt}\\
We now return to the unification of the Lagrangian and Hamiltonian formalisms. Throughout, a tilde indicates that the corresponding map has been left-trivialized according to the preceding constructions. 

The entire unification is summarized in the diagram below. The maps and expressions shown in blue correspond to the Hamiltonian picture, while those shown in red correspond to the Lagrangian picture. One again observes that the Hamilton's equations and the Euler-Lagrange equations each define Lagrangian submanifolds of $(G \ltimes \mathfrak{g}^*) \ltimes (\mathfrak{g} \ltimes \mathfrak{g}^*)$. 
\begin{center}
\begin{tikzpicture}[scale=0.7]
%
\node at (0,0) {$(G \ltimes \mathfrak{g}^*) \ltimes(\mathfrak{g} \ltimes \mathfrak{g}^*)$};
\draw[->] (2.3,0)--(2.7,0);
\node at (2.5,0.5) {\scriptsize$\Tilde{\alpha}_G$};
\draw[->] (-2.3,0)--(-2.7,0);
\node at (-2.5,0.5) {\scriptsize$\Tilde{\beta}_G$};
\draw[->] (-0.5,-0.5)--(-2,-1.5);
\node at (-1.8,-0.8) {\scriptsize$pr_{12}$};
\draw[<-,blue] (-0.3,-0.7)--(-1.8,-1.7);
\node at (-0.6,-1.4) {\textcolor{blue}{\scriptsize$\Tilde{X}_{\Tilde{\mathscr{H}}}$}};
\draw[->] (0.5,-0.5)--(2,-1.5);
\node at (1.8,-0.8) {\scriptsize$pr_{13}$};
\node at (-5,0) {$(G \ltimes \mathfrak{g}^*) \ltimes(\mathfrak{g}^* \ltimes \mathfrak{g})$};
\draw[->] (-4.5,-0.5)--(-3,-1.5);
\node at (-3.2,-0.8) {\scriptsize$pr_{12}$};
\draw[<-,blue] (-4.7,-0.7)--(-3.2,-1.7);
\node at (-4.4,-1.4) {\textcolor{blue}{\scriptsize$\Tilde{d\Tilde{\mathscr{H}}}$}};
\draw[->] (-2,-2.5)--(-0.5,-3.5);
\node at (-1.5,-3.2) {\scriptsize$pr_1$};
\draw[->] (4.5,-0.5)--(3,-1.5);
\node at (3.2,-0.8) {\scriptsize$pr_{12}$};
\draw[<-,red] (4.7,-0.7)--(3.2,-1.7);
\node at (4.4,-1.4) {\textcolor{red}{\scriptsize$\Tilde{d\Tilde{\mathscr{L}}}$}};
\draw[->] (2,-2.5)--(0.5,-3.5);
\node at (1.5,-3.2) {\scriptsize$pr_1$};
\node at (5,0) {$(G \ltimes \mathfrak{g}) \ltimes(\mathfrak{g}^* \ltimes \mathfrak{g}^*)$};
\node at (-2.5,-2) {$G \ltimes \mathfrak{g}^*$};
\draw[blue,->] (-3,-2.5)--(-4.5,-3.5);
\node at (-3.5,-3.2) {\textcolor{blue}{\scriptsize$\Tilde{\mathscr{H}}$}};
\node at (-5,-4) {$\mathbb{R}$};
\node at (2.5,-2) {$G \ltimes \mathfrak{g}$};
\node at (0,-4) {$G$};
\draw[red,->] (3,-2.5)--(4.5,-3.5);
\node at (3.5,-3.2) {\textcolor{red}{\scriptsize$\Tilde{\mathscr{L}}$}};
\node at (5,-4) {$\mathbb{R}$};
\node at (5,1.3) {\textcolor{red}{\scriptsize$(g,\xi,ad_\xi^*(\frac{\delta \Tilde{\mathscr{L}}}{\delta \xi})+T_e^*L_g(\frac{\delta \Tilde{\mathscr{L}}}{\delta g}),\frac{\delta \Tilde{\mathscr{L}}}{\delta \xi})$}};
\node at (-5,1.3) {\textcolor{blue}{\scriptsize$(g,\mu,ad_{\frac{\delta \Tilde{\mathscr{H}}}{\delta \mu}}^*(\mu)-T_e^*L_g(\frac{\delta \Tilde{\mathscr{H}}}{\delta g}),\frac{\delta \Tilde{\mathscr{H}}}{\delta \mu})$}};
\node at (0,1.3) {\textcolor{red}{\scriptsize$(g,\frac{\delta \Tilde{\mathscr{L}}}{\delta \xi},\xi,T_e^*L_g(\frac{\delta \Tilde{\mathscr{L}}}{\delta g})$}};
\node at (0,2) {\textcolor{blue}{\scriptsize$(g,\mu,-\frac{\delta \Tilde{\mathscr{H}}}{\delta \mu},T_e^*L_g(\frac{\delta \Tilde{\mathscr{H}}}{\delta g})$}};
\node at (4,-2) {\textcolor{red}{\scriptsize$(g,\xi)$}};
\node at (-4,-2) {\textcolor{blue}{\scriptsize$(g,\mu)$}};
\node at (5,-4.7) {\textcolor{red}{\scriptsize$\Tilde{\mathscr{L}}(g,\xi)$}};
\node at (-5,-4.7) {\textcolor{blue}{\scriptsize$\Tilde{\mathscr{H}}(g,\mu)$}};
\node at (0,-4.7) {\scriptsize$g$};
\end{tikzpicture}
\end{center}
In the presence of a left-translational symmetry, the diagram reduces accordingly as shown below. This 
reduced diagram highlights that the Lie-Poisson and Euler-Poincar{\'e} equations also define Lagrangian submanifolds of $(G \ltimes \mathfrak{g}^*) \ltimes (\mathfrak{g} \ltimes \mathfrak{g}^*)$. 
\begin{center}
\begin{tikzpicture}[scale=0.7]
%
\node at (0,0) {$(G \ltimes \mathfrak{g}^*) \ltimes(\mathfrak{g} \ltimes \mathfrak{g}^*)$};
\draw[->] (2.3,0)--(2.7,0);
\node at (2.5,0.5) {\scriptsize$\Tilde{\alpha}_G$};
\draw[->] (-2.3,0)--(-2.7,0);
\node at (-2.5,0.5) {\scriptsize$\Tilde{\beta}_G$};
\draw[->] (-0.5,-0.5)--(-2,-1.5);
\node at (-1.8,-0.8) {\scriptsize$pr_{12}$};
\draw[<-,blue] (-0.3,-0.7)--(-1.8,-1.7);
\node at (-0.6,-1.4) {\textcolor{blue}{\scriptsize$\Tilde{X}_h$}};
\draw[->] (0.5,-0.5)--(2,-1.5);
\node at (1.8,-0.8) {\scriptsize$pr_{13}$};
\node at (-5,0) {$(G \ltimes \mathfrak{g}^*) \ltimes(\mathfrak{g}^* \ltimes \mathfrak{g})$};
\draw[->] (-4.5,-0.5)--(-3,-1.5);
\node at (-3.2,-0.8) {\scriptsize$pr_{12}$};
\draw[<-,blue] (-4.7,-0.7)--(-3.2,-1.7);
\node at (-4.4,-1.4) {\textcolor{blue}{\scriptsize$\Tilde{dh}$}};
\draw[->] (-2,-2.5)--(-0.5,-3.5);
\node at (-1.5,-3.2) {\scriptsize$pr_1$};
\draw[->] (4.5,-0.5)--(3,-1.5);
\node at (3.2,-0.8) {\scriptsize$pr_{12}$};
\draw[<-,red] (4.7,-0.7)--(3.2,-1.7);
\node at (4.4,-1.4) {\textcolor{red}{\scriptsize$\Tilde{dl}$}};
\draw[->] (2,-2.5)--(0.5,-3.5);
\node at (1.5,-3.2) {\scriptsize$pr_1$};
\node at (5,0) {$(G \ltimes \mathfrak{g}) \ltimes(\mathfrak{g}^* \ltimes \mathfrak{g}^*)$};
\node at (-2.5,-2) {$G \ltimes \mathfrak{g}^*$};
\draw[blue,->] (-3,-2.5)--(-4.5,-3.5);
\node at (-3.5,-3.2) {\textcolor{blue}{\scriptsize$h$}};
\node at (-5,-4) {$\mathbb{R}$};
\node at (2.5,-2) {$G \ltimes \mathfrak{g}$};
\node at (0,-4) {$G$};
\draw[red,->] (3,-2.5)--(4.5,-3.5);
\node at (3.5,-3.2) {\textcolor{red}{\scriptsize$l$}};
\node at (5,-4) {$\mathbb{R}$};
\node at (5,1) {\textcolor{red}{\scriptsize$(g,\xi,ad_\xi^*(\frac{\delta l}{\delta \xi}),\frac{\delta l}{\delta \xi})$}};
\node at (-5,1) {\textcolor{blue}{\scriptsize$(g,\mu,ad_{\frac{\delta h}{\delta \mu}}^*(\mu),\frac{\delta h}{\delta \mu})$}};
\node at (5,1) {\textcolor{red}{\scriptsize$(g,\xi,ad_\xi^*(\frac{\delta l}{\delta \xi}),\frac{\delta l}{\delta \xi})$}};
\node at (0,1) {\textcolor{red}{\scriptsize$(g,\frac{\delta l}{\delta \xi},\xi,0)$}};
\node at (0,1.5) {\textcolor{blue}{\scriptsize$(g,\mu,-\frac{\delta h}{\delta \mu},0)$}};
\node at (4,-2) {\textcolor{red}{\scriptsize$(g,\xi)$}};
\node at (-4,-2) {\textcolor{blue}{\scriptsize$(g,\mu)$}};
\node at (5,-4.7) {\textcolor{red}{\scriptsize$l(\xi)$}};
\node at (-5,-4.7) {\textcolor{blue}{\scriptsize$h(\mu)$}};
\node at (0,-4.7) {\scriptsize$g$};
\end{tikzpicture}
\end{center}

\chapter{Structure-preserving Numerical Integrators}
\label{ch5}
In chapter \ref{ch1} we defined linear and nonlinear symplectic structures, in chapter \ref{ch3} we introduced a retraction maps-based framework for constructing numerical integrators for systems evolving on manifolds, and in chapter \ref{ch4} we examined the symplectic geometry underpinning mechanical systems in detail. We are now ready to introduce a retraction-maps-based framework for constructing structure-preserving numerical integrators for mechanical systems. 

Since the dynamics of mechanical systems are inherently second order, we must lift our discretization maps from the level of the configuration manifold to the level of its tangent and cotangent bundles.
\section{Tangent-Lifted Discretization Maps}
\label{sec tangent-lifted disc}
Let $\mathcal{D} : TQ \to Q \times Q$ be a discretization map on a manifold $Q$. It can be lifted to the tangent bundle via the commutative diagram below:
\begin{center}
    \begin{tikzcd}
        TTQ \arrow[r,"\mathcal{D}^T"] \arrow[d,shift right,swap,"\kappa_Q"] & TQ \times TQ \arrow[d,equal] \\
        TTQ \arrow[r,"T\mathcal{D}"] \arrow[d,swap,"\tau_{TQ}"] \arrow[u,shift right,swap,"\kappa_Q"] & T(Q \times Q) \arrow[d,"\tau_{Q \times Q}"] \\
        TQ \arrow[r,swap,"\mathcal{D}"] & Q \times Q
    \end{tikzcd}
\end{center}
and is defined by
\begin{align} \label{tangent-lifted disc}
    \mathcal{D}^T := T\mathcal{D} \circ \kappa_Q,
\end{align}
where $\kappa_Q : TTQ \to TTQ$ denotes the canonical involution introduced in the previous chapter. Here we have used the  natural identification $T(Q \times Q) \cong TQ \times TQ$.  
\section{Cotangent-Lifted Discretization Maps}
\label{sec cotangent-lifted disc}
We can similarly lift the discretization map to the cotangent bundle. Let $\mathcal{D} : TQ \to Q \times Q$ be a discretization map. Its cotangent lift is defined by the commutative diagram
\begin{center}
    \begin{tikzcd}
        TT^*Q \arrow[r,"\mathcal{D}^{T^*}"] \arrow[d,swap,"\alpha_Q"] & T^*Q \times T^*Q \arrow[d,"\Phi"] \\
        T^*TQ \arrow[d,swap,"\pi_{TQ}"] & T^*(Q \times Q) \arrow[d,"\pi_{Q \times Q}"] \arrow[l,"T^*\mathcal{D}",swap] \\
        TQ \arrow[r,swap,"\mathcal{D}"] & Q \times Q
    \end{tikzcd}
\end{center}
and is given by
\begin{align} \label{cotangent-lifted disc}
    \mathcal{D}^{T^*} := \Phi^{-1} \circ (T^*\mathcal{D})^{-1} \circ \alpha_Q.
\end{align}
Here $\alpha_Q : TT^*Q \to T^*TQ$ is the canonical symplectomorphism with $TT^*Q$ equipped with the tangent-lifted symplectic structure. The map $\Phi : T^*Q \times T^*Q \to T^*(Q \times Q)$ is the canonical vector bundle isomorphism which is a symplectomorphism when $T^*Q \times T^*Q$ is equipped with the twisted product symplectic structure
\begin{align*}
    \omega_Q^\ominus := pr_1^* \omega_Q - pr_2^* \omega_Q,
\end{align*}
and $T^*(Q \times Q)$ is equipped with its canonical symplectic form $\omega_{Q \times Q}$. Consequently, $\mathcal{D}^{T^*}$ is a symplectomorphism. This cotangent lift provides the geometric foundation for constructing variational and symplectic integrators from discretizations on $Q$.
\section{Numerical Integrators Using Lifted \\
Discretization Maps}
\label{sec lifted integrators}
Consider a manifold $Q$. Let $X : TQ \to TTQ$ be a vector field on the tangent bundle whose flow we wish to discretize. We begin by choosing a discretization map $\mathcal{D} : TQ \to Q \times Q$ on the base manifold $Q$. This map encodes how a tangent vector approximates a finite displacement on $Q$. Next, we lift this discretization to the double tangent bundle by taking its tangent lift $\mathcal{D}^T : TTQ \to TQ \times TQ$
defined in \eqref{tangent-lifted disc}. We then apply the general integrator construction introduced in \eqref{integrator 1} to obtain a discrete update rule on $TQ$, summarized by the commutative diagram
\begin{center}
    \begin{tikzcd}
        TQ \times TQ \arrow[r,"(\mathcal{D}^{T})^{-1}"] \arrow[d,swap,"(\mathcal{D}^{T})^{-1}"] & TTQ \arrow[d,"\tau_{TQ}"] \\
        TTQ & TQ \arrow[l,"X"]
    \end{tikzcd}
\end{center}
and defined implicitly by
\begin{align}
    hX(\tau_{TQ}((\mathcal{D}^{T})^{-1}(q_k,v_k;q_{k+1},v_{k+1})))=(\mathcal{D}^{T})^{-1}(q_k,v_k;q_{k+1},v_{k+1}) \label{tangent lifted integrator}
\end{align}
where $h \in \mathbb{R}$ denotes the step size. This construction is particularly well suited to the case in which $X$ is a second-order vector field or semispray. Indeed, the presence of the canonical involution in the definition of $\mathcal{D}^T$ ensures that the discrete dynamics respect the intrinsic second-order structure of the continuous system. A canonical example of such dynamics is provided by hyperregular Lagrangian systems, where the Euler–Lagrange vector field defines a semispray, as given in \eqref{Lagrangian semi-spray eq}.
\\

Now let $X: T^*Q \to TT^*Q$ be a vector field on the cotangent bundle whose flow we wish to discretize. As before, we begin by choosing a discretization map $\mathcal{D} : TQ \to Q \times Q$ on the base manifold $Q$. Next, we lift this discretization to the tangent of the cotangent bundle by taking its cotangent lift $\mathcal{D}^{T^*} : TT^*Q \to T^*Q \times T^*Q$ defined in \eqref{cotangent-lifted disc}. Applying the general integrator construction of \eqref{integrator 1}, we obtain a discrete update rule on $T^*Q$, summarized by the commutative diagram
\begin{center}
    \begin{tikzcd}
        T^*Q \times T^*Q \arrow[r,"(\mathcal{D}^{T^*})^{-1}"] \arrow[d,swap,"(\mathcal{D}^{T^*})^{-1}"] & TT^*Q \arrow[d,"\tau_{T^*Q}"] \\
        TT^*Q & T^*Q \arrow[l,"X"]
    \end{tikzcd}
\end{center}
and defined implicitly by
\begin{align}
    hX(\tau_{T^*Q}((\mathcal{D}^{T^*})^{-1}(q_k,p_k;q_{k+1},p_{k+1})))=(\mathcal{D}^{T^*})^{-1}(q_k,p_k;q_{k+1},p_{k+1}) \label{cotangent lifted integrator}
\end{align}
where $h \in \mathbb{R}$ denotes the step size. This construction is particularly well suited to the case in which $X$ is a Hamiltonian vector field $X_\mathscr{H}$ associated with a Hamiltonian function $\mathscr{H}$ as defined in \eqref{Hamilton's eq}. Since the cotangent-lifted discretization map $\mathcal{D}^{T^*}$ is a symplectomorphism, the Lagrangian submanifold of $TT^*Q$ defined by $X_\mathscr{H}$ is preserved under the discrete flow. Consequently, the resulting numerical scheme is a symplectic integrator. Let us validate this construction using a simple example.

Let $Q = \mathbb{R}^n$, so that $T^*Q \cong \mathbb{R}^n \times \mathbb{R}^n$, and consider a vector field $X \in \mathfrak{X}(T^*Q)$ defined as $X(q,p) := (q,p,f_1(q,p),f_2(q,p))$ for all $(q,p) \in T^*Q$, where $f_1, f_2 : \mathbb{R}^n \times \mathbb{R}^n \to \mathbb{R}^n$. Choose the discretization map $\mathcal{D} : TQ \to Q \times Q$ described in example \ref{eg disc R^n}. Taking its cotangent lift and applying the cotangent-lifted integrator \eqref{cotangent lifted integrator}, we obtain the following one-parameter family of integrators on $T^*Q$:
\begin{align*}
    q_{k+1} &= q_k + h f_1((1-\theta)q_k+\theta q_{k+1},\theta p_k + (1-\theta)p_{k+1}) \quad \text{and} \\
    p_{k+1} &= p_k + h f_2((1-\theta)q_k+\theta q_{k+1},\theta p_k + (1-\theta)p_{k+1})
\end{align*}
where $h \in \mathbb{R}$ is the step size and $\theta \in [0,1]$ parameterizes the family. For $\theta = 0$, this scheme reduces to symplectic Euler A, while for $\theta = 1$, it reduces to symplectic Euler B. In particular, the asymmetry in the evaluation of the position and momentum variables is a direct consequence of the twisted-product symplectic structure used in the cotangent-lifted discretization.
\section{Tangent-Lifted Trivialized Discretization \\ Maps}
\label{sec tangent-lifted triv disc}
Let $\mathcal{D}^L : G \ltimes \mathfrak{g} \to G \times G$ be left-trivialized discretization map on a Lie group $G$. It can be lifted to the tangent bundle via the commutative diagram below:
\begin{center}
    \begin{tikzcd}
        (G \ltimes \mathfrak{g}) \ltimes (\mathfrak{g} \ltimes \mathfrak{g}) \arrow[d,shift left,"\widetilde{\kappa}_G^L"]  \arrow[r,"(\mathcal{D}^{L,L})^T"] & (G \ltimes \mathfrak{g}) \times (G \ltimes \mathfrak{g}) \\
        (G \ltimes \mathfrak{g}) \ltimes (\mathfrak{g} \ltimes \mathfrak{g}) \arrow[u,shift left,"\widetilde{\kappa}_G^L"] & (G \times G) \ltimes (\mathfrak{g} \times \mathfrak{g}) \arrow[u,"sw_{23}",swap] \\
        T(G \ltimes \mathfrak{g}) \arrow[d,"\tau_{G \ltimes \mathfrak{g}}",swap] \arrow[r,"T\mathcal{D}^L"] \arrow[u,"tr^L_{T(G \ltimes \mathfrak{g})}"] & T(G \times G) \arrow[d,"\tau_{G \times G}"] \arrow[u,"tr^L_{T(G \times G)}",swap] \\
        G \ltimes \mathfrak{g} \arrow[r,"\mathcal{D}^L",swap] & G \times G
    \end{tikzcd}
\end{center}
and is defined by
\begin{align}
    (\mathcal{D}^{L,L})^T := sw_{23} \circ tr_{T(G \times G)}^L \circ T\mathcal{D}^L \circ (tr_{T(G \ltimes \mathfrak{g})}^L)^{-1} \circ \tilde{\kappa}_G^L, \label{tangent lifted triv disc}
\end{align}
where $\tilde{\kappa}_G^L$ is the left-trivialized canonical involution and $sw_{23}$ denotes the permutation exchanging the second and the third components so as to identify $(G \times G) \ltimes (\mathfrak{g} \times \mathfrak{g}) \cong (G \ltimes \mathfrak{g}) \times (G \ltimes \mathfrak{g})$. This construction mirrors the tangent lift of discretization maps on manifolds while preserving the Lie group structure through left trivialization.
\section{Cotangent-Lifted Trivialized Discretization Maps}
\label{sec cotangent-lifted triv disc}
We can similarly lift the left-trivialized discretization map $\mathcal{D}^L : G \ltimes \mathfrak{g} \to G \times G$ to the cotangent bundle via the commutative diagram below:
\begin{center}
    \begin{tikzcd}
        (G \ltimes \mathfrak{g}^*) \ltimes (\mathfrak{g} \ltimes \mathfrak{g}^*) \arrow[d,"\tilde{\alpha}_G^L",swap]  \arrow[r,"(\mathcal{D}^{L,L})^{T^*}"] & (G \ltimes \mathfrak{g}^*) \times (G \ltimes \mathfrak{g}^*) \arrow[d,"\tilde{\Phi}^L"] \\
        (G \ltimes \mathfrak{g}) \ltimes (\mathfrak{g}^* \ltimes \mathfrak{g}^*) & (G \times G) \ltimes (\mathfrak{g}^* \times \mathfrak{g}^*) \\
        T^*(G \ltimes \mathfrak{g}) \arrow[d,"\pi_{G \ltimes \mathfrak{g}}",swap] \arrow[u,"tr^L_{T^*(G \ltimes \mathfrak{g})}"] & T^*(G \times G) \arrow[d,"\pi_{G \times G}"] \arrow[l,"T^*\mathcal{D}^L",swap] \arrow[u,"tr^L_{T^*(G \times G)}",swap] \\
        G \ltimes \mathfrak{g} \arrow[r,"\mathcal{D}^L",swap] & G \times G
    \end{tikzcd}
\end{center}
and is defined by
\begin{align}
    (\mathcal{D}^{L,L})^{T^*} := (\tilde{\Phi}^L)^{-1} \circ tr_{T^*(G \times G)}^L \circ (T^*\mathcal{D}^L)^{-1} \circ (tr_{T^*(G \ltimes \mathfrak{g})}^L)^{-1} \circ \tilde{\alpha}_G^L \label{cotangent lifted triv disc}
\end{align}
where $\tilde{\alpha}_G^L$ is the left-trivialized symplectomorphism dual to the canonical involution and $\tilde{\Phi}^L$ is the left-trivialization of the symplectomorphism $\Phi$ defined in \eqref{cotangent-lifted disc}. All the maps appearing in the above composition are symplectomorphisms; consequently, the lifted discretization map $(\mathcal{D}^{L,L})^{T^*}$ is itself a symplectomorphism. This construction therefore mirrors the cotangent lift of discretization maps on manifold, while respecting the Lie group structure through left trivialization.
\section{Numerical Integrators Using Lifted \\ Trivialized 
Discretization Maps}
\label{sec lifted triv integrators}
Just as lifted discretization maps from a base manifold to its tangent or cotangent bundle were used to construct numerical integrators for mechanical systems in section \ref{sec numerical integrators using disc maps}, we now specialize this framework to the setting of Lie groups. 
\\

Let $X : TG \to TTG$ be a vector field on the tangent bundle of a Lie group $G$ whose flow we wish to discretize. We begin by left (or right) trivializing the vector field to obtain $\tilde{X}^L : G \ltimes \mathfrak{g}  \to (G \ltimes \mathfrak{g}) \ltimes (\mathfrak{g} \ltimes \mathfrak{g})$ as defined earlier. Next, we choose a left (or right) trivialized discretization map $\mathcal{D}^L : G \ltimes \mathfrak{g} \to G \times G$ on the base Lie group $G$. We then lift this discretization map to the trivialized tangent bundle to obtain $(\mathcal{D}^{L,L})^T$, as defined in \eqref{tangent lifted triv disc}. Applying the general integrator construction introduced in \eqref{tangent lifted integrator}, we arrive at a discrete update rule on $G \ltimes \mathfrak{g}$, summarized by the commutative diagram
\begin{center}
    \begin{tikzcd}[column sep=large]
        (G \ltimes \mathfrak{g}) \times (G \ltimes \mathfrak{g}) \arrow[r,"((\mathcal{D}^{L,L})^T)^{-1}"] \arrow[d,"((\mathcal{D}^{L,L})^T)^{-1}",swap] & (G \ltimes \mathfrak{g}) \ltimes (G \ltimes \mathfrak{g}) \arrow[d,"pr_{12}"] \\
        (G \ltimes \mathfrak{g}) \ltimes (G \ltimes \mathfrak{g})& G \ltimes \mathfrak{g} \arrow[l,"\tilde{X}^L"]
    \end{tikzcd}
\end{center}
and defined implicitly by
\small{
\begin{align} \label{triv integrator tangent}
    h \tilde{X}^L(pr_{12}(((&\mathcal{D}^{L,L})^T)^{-1}(g_k,\xi_k,g_{k+1},\xi_{k+1})) \nonumber\\
    = ((&\mathcal{D}^{L,L})^T)^{-1}(g_k,\xi_k,g_{k+1},\xi_{k+1}) 
\end{align}
}
where $h \in \mathbb{R}$ denotes the step size. This construction is particularly well suited for systems such as the Euler–Arnold equations, where the dynamics naturally evolve on the Lie algebra.
\\

Now let $X : T^*G \to TT^*G$ be a vector field on the cotangent bundle of a Lie group $G$ whose flow we wish to discretize. Again, we  begin by left (or right) trivializing the vector field to obtain $\tilde{X}^L : G \ltimes \mathfrak{g}^*  \to (G \ltimes \mathfrak{g}) \ltimes (\mathfrak{g}^* \ltimes \mathfrak{g}^*)$ as defined earlier. Next, we choose a left (or right) trivialized discretization map $\mathcal{D}^L : G \ltimes \mathfrak{g} \to G \times G$ on the base Lie group $G$. Finally, we lift this discretization map to the trivialized cotangent bundle to obtain $(\mathcal{D}^{L,L})^{T^*}$, as defined in \eqref{cotangent lifted triv disc}. Applying the general integrator construction introduced in \eqref{cotangent lifted integrator}, we arrive at a discrete update rule on $G \ltimes \mathfrak{g}^*$, summarized by the commutative diagram
\begin{center}
    \begin{tikzcd}[column sep=large]
        (G \ltimes \mathfrak{g}^*) \times (G \ltimes \mathfrak{g}^*) \arrow[r,"((\mathcal{D}^{L,L})^{T^*})^{-1}"] \arrow[d,"((\mathcal{D}^{L,L})^{T^*})^{-1}",swap] & (G \ltimes \mathfrak{g}^*) \ltimes (\mathfrak{g} \ltimes \mathfrak{g}^*) \arrow[d,"pr_{12}"] \\
        (G \ltimes \mathfrak{g}^*) \ltimes (\mathfrak{g} \ltimes \mathfrak{g}^*)& G \ltimes \mathfrak{g}^* \arrow[l,"\tilde{X}^L"]
    \end{tikzcd}
\end{center}
and defined implicitly by
\small{
\begin{align}
    h \tilde{X}^L(pr_{12}(((&\mathcal{D}^{L,L})^{T^*})^{-1}(g_k,\mu_k,g_{k+1},\mu_{k+1})) \nonumber \\
    = ((&\mathcal{D}^{L,L})^{T^*})^{-1}(g_k,\mu_k,g_{k+1},\mu_{k+1}) \label{triv integrator cotangent}
\end{align}}
where $h \in \mathbb{R}$ denotes the step size. This construction is particularly well suited for Lie–Poisson type systems, whose dynamics evolve on the dual of a Lie algebra and are governed by an underlying Poisson geometric structure.
\\

Note that one is free to choose either left or right trivialization for both the vector field (L or R) and the discretization map (LL, LR, RL, RR). Different combinations of these choices lead to distinct, but equally natural, classes of numerical integrators on Lie groups.

\chapter{Illustrative Examples}
\label{ch6}
Finally, we demonstrate the construction of retraction map–based numerical integrators through examples that are particularly relevant to engineering and physics. We also juxtapose the resulting schemes with well-established methods in order to highlight their advantages. The examples considered are the rotational dynamics of a rigid body, the dynamics of a heavy top and the dynamics of a quadrotor. 
\section{Rigid Body}
\label{sec rigid body}
The freely rotating rigid body is the quintessential example of an Euler-Poincar{\'e} or Lie-Poisson system, admitting a quadratic Lagrangian or Hamiltonian, respectively. We therefore begin by constructing numerical integrators for such systems on a general Lie group $G$, before specializing to the group $G = \mathbb{SO}(3)$ which serves as the configuration space for this problem. 

We begin by making two choices: (i) a trivialized retraction map, and (ii) a parameter $s \in [0,1]$. We choose the trivialized retraction map 
\begin{align*}
    \mathcal{R}^L : G \times \mathfrak{g} \to G, \quad \mathcal{R}^L(g,\xi) := g \tau(\xi),
\end{align*}
as in proposition \ref{prop 2.1}. Setting $s=0$ in proposition \ref{prop 2.4} yields the corresponding trivialized discretization map
\begin{align*}
    \mathcal{D}^L : G \times \mathfrak{g} \to G \times G, \quad \mathcal{D}^L(g,\xi) := (g, g \tau(\xi)).
\end{align*}
Finally, we make the third choice: (iii) lifting the trivialized discretization map to the cotangent bundle, either by left or right trivialization. Choosing the left lift via \eqref{cotangent lifted triv disc} as the dynamics are left-invariant, we obtain the cotangent-lifted trivialized discretization
\begin{align*}
    (\mathcal{D}^{L,L})^{T^*} : (G \times \mathfrak{g}^*) \times (\mathfrak{g} \times \mathfrak{g}^*) \to (G \times \mathfrak{g}^*) \times (G \times \mathfrak{g}^*),
\end{align*}
whose inverse is
\begin{align}
    ((\mathcal{D}^{L,L})^{T^*})^{-1} (g_k, \mu_k; g_{k+1}, \mu_{k+1}) = (&g_k, d^{L*}_{\tau^{-1}(g_k^{-1}g_{k+1})}\tau(\mu_{k+1});\tau^{-1}(g_k^{-1}g_{k+1}) \nonumber \\
    &Ad^*_{(g_k^{-1}g_{k+1})^{-1}}(\mu_{k+1})-\mu_k). \label{inverse of left triv disc rigid body}
\end{align}
Choosing instead the right lift gives
\begin{align*}
    (\mathcal{D}^{L,R})^{T^*} : (G \times \mathfrak{g}^*) \times (\mathfrak{g} \times \mathfrak{g}^*) \to (G \times \mathfrak{g}^*) \times (G \times \mathfrak{g}^*),
\end{align*}
whose inverse is 
\begin{align}
    ((\mathcal{D}^{L,R})^{T^*})^{-1} (g_k, \mu_k; g_{k+1}, \mu_{k+1}) = (&g_k, d^{L*}_{\tau^{-1}(g_k^{-1}g_{k+1})}\tau(Ad^*_{g_{k+1}}(\mu_{k+1})); \nonumber \\
    &\tau^{-1}(g_k^{-1}g_{k+1}),\mu_{k+1}-\mu_k). \label{inverse of right triv disc rigid body}
\end{align}
Let the left-reduced Hamiltonian be $h(\mu) := \tfrac{1}{2} \langle \mu, \mathbb{I}^{-1}(\mu) \rangle$ where $\mathbb{I} : \mathfrak{g} \to \mathfrak{g}^*$ is the \emph{inertia map}, see \cite{arnold1999topological, khesin2008geometry}. The corresponding left Lie–Poisson equations are 
\begin{align}
    \dot{\mu} = ad^*_{\mathbb{I}^{-1}(\mu)}(\mu) \quad \text{and} \quad \dot{g} = T_eL_g(\mathbb{I}^{-1}(\mu)). \label{LP eq rigid body}
\end{align}
The associated left-trivialized vector field is
\begin{align}
    \Tilde{X}_h^L(g,\mu) := (g,\mu, \mathbb{I}^{-1}(\mu),0). \label{triv Lie-Poisson vector field}
\end{align}
Substituting \eqref{inverse of left triv disc rigid body} and \eqref{triv Lie-Poisson vector field} into \eqref{triv integrator cotangent} produces the scheme
\begin{equation} \label{integrator eg 1 left}
    \boxed{
    \begin{aligned}
        g_{k+1} &= g_k \tau(t\xi_k) \\
        d^{L*}_{t\xi_k} \tau (\mu_{k+1}) &= \mathbb{I}(t\xi_k) \\
        \mu_{k+1} &= Ad^*_{\tau(t\xi_k)}(\mu_k)
    \end{aligned}
    }
\end{equation}
where $t \in \mathbb{R}$ is the step size. 
Similarly, inserting \eqref{inverse of right triv disc rigid body} yields
\begin{equation} \label{integrator eg 1 right}
    \boxed{
    \begin{aligned}
        g_{k+1} &= g_k \tau(t\xi_k) \\
        d^{L*}_{t\xi_k} \tau (Ad^*_{g_{k+1}}(\mu_{k+1})) &= \mathbb{I}(t\xi_k) \\
        \mu_{k+1} &= \mu_k
    \end{aligned}
    }.
\end{equation}
\begin{center}
\begin{tikzpicture}[scale = 0.6]
\shadedraw[rotate around={-45:(5,5)}] (5,5) ellipse (2 and 4);
\node at (6,8.5) {$\mathcal{I}$};
\draw[thick,red,-stealth] (4,4)--(8,4);
\draw[thick,red,-stealth] (4,4)--(4,8);
\draw[thick,red,-stealth] (4,4)--(1.5,1.5);
\node at (1,1) {\textcolor{red}{$\mathcal{S}$}};
\draw[thick,blue,-stealth] (4,4)--(7,5.5);
\draw[thick,blue,-stealth] (4,4)--(2,7.5);
\draw[thick,blue,-stealth] (4,4)--(7.5,2);
\node at (8,2) {\textcolor{blue}{$\mathcal{B}$}};
\draw[->] (1.9,7.9) arc (113:90:5);
\node at (2.7,8.5) {$R$};
\end{tikzpicture}
\end{center}
Consider the rigid body shown above, with spatial frame $\mathcal{S}$, body frame $\mathcal{B}$ and transformation $R \in \mathbb{SO}(3)$ mapping body coordinates to spatial coordinates. We begin by fixing the following identifications which can be found in \cite{holm2008geometric, holm2011geometric1, holm2011geometric2} or \cite{marsden2002introduction}.
\paragraph{Hat map} The isomorphism $\hat{\cdot}:\mathbb{R}^3 \to \mathfrak{so}(3)$ is defined by  $x=\left(\begin{smallmatrix}
    x_1 \\
    x_2 \\
    x_3
\end{smallmatrix}\right) \mapsto \hat{x} := \left(\begin{smallmatrix}
    0 & -x_3 & x_2 \\
    x_3 & 0 & -x_1 \\
    -x_2 & x_1 & 0
\end{smallmatrix}\right)$. It satisfies $\hat{x}y = x \times y$ for all $x,y \in \mathbb{R}^3$.
\paragraph{Breve map} The isomorphism $\Breve{\cdot}:\mathbb{R}^3 \to \mathfrak{so}(3)^*$ is defined as $\breve{x} := \tfrac{1}{2}\hat{x}$ via the trace pairing $\langle \breve{x}, \hat{y} \rangle := \tfrac{1}{2}\text{trace}(\breve{x}\hat{y}^T) = x \cdot y = x^Ty = -\tfrac{1}{2}\text{trace}(\breve{x}\hat{y})$ for all $x,y \in \mathbb{R}^3$.
\vspace{10pt}
\\
With these identifications, the left-reduced Hamiltonian takes the form
\begin{align*}
    h(\breve{\Pi})=\tfrac{1}{2}\langle \breve{\Pi},\mathbb{I}^{-1}(\breve{\Pi}) \rangle = \tfrac{1}{2} \Pi \cdot \mathcal{I}^{-1}\Pi,
\end{align*}
where $\Pi \in \mathbb{R}^3$ denotes the body angular momentum and $\mathcal{I} \in \mathbb{R}^{3 \times 3}$ is the moment of inertia matrix. The inertia map $\mathbb{I} : \mathfrak{so}(3) \to \mathfrak{so}(3)^*$ 
is therefore given by $\mathbb{I}(\hat{\Omega}) := \breve{(\mathcal{I}\Omega)}$ where $\Omega \in \mathbb{R}^3$ denotes the body angular velocity. Under these identifications, the left Lie-Poisson equations reduce to the classical rigid body equations 
\begin{equation}
    \boxed{
    \begin{aligned}
        \dot{R} &= R \hat{\Omega} \\
        \Pi &= \mathcal{I} \Omega \\
        \dot{\Pi} &= \Pi \times \Omega
    \end{aligned}
    }
\end{equation}
where $R \in \mathbb{SO}(3)$ describes the orientation of the rigid body. In addition to the energy (Hamiltonian), the squared magnitude of the body angular momentum $\|\Pi\|^2$ is a casimir arising from the coadjoint orbit geometry. The integrator \eqref{integrator eg 1 left} becomes
\begin{equation}
    \boxed{
    \begin{aligned}
        R_{k+1} &= R_k \tau(t \hat{\Omega}_k) \\
        d^{L*}_{t\hat{\Omega}_k} \tau (\breve{\Pi}_{k+1}) &= \mathbb{I}(t\hat{\Omega}_k) \\
        \breve{\Pi}_{k+1} &= Ad^*_{\tau(t\hat{\Omega}_k)}(\breve{\Pi}_k)
    \end{aligned} \label{integrator eg 1 left SO(3)}
    }
\end{equation}
while the right-lift variant reads
\begin{equation}
    \boxed{
    \begin{aligned}
        R_{k+1} &= R_k \tau(t \hat{\Omega}_k) \\
        d^{L*}_{t\hat{\Omega}_k} \tau (Ad^*_{R_{k+1}}(\breve{\Pi}_{k+1})) &= \mathbb{I}(t\hat{\Omega}_k) \\
        \breve{\Pi}_{k+1} &= \breve{\Pi}_k
    \end{aligned} \label{integrator eg 1 right SO(3)}
    }.
\end{equation}
where $t \in \mathbb{R}$ is the time step. This scheme preserves the coadjoint orbits and therefore the Casimir invariants.

The map $\tau : \mathfrak{so}(3) \to \mathbb{SO}(3)$ is a local diffeomorphism; the most commonly used choices in practice are the matrix exponential and the Cayley map.
\subsection{Numerical Integrator Using Exponential Map}
\label{subsec rigid body exp integrator}
We first choose the matrix exponential map $\exp : \mathfrak{so}(3) \to \mathbb{SO}(3)$ as the local diffeomorphism $\tau$. For the $\mathbb{SO}(3)$ group, the coadjoint $ad^*$ and Coadjoint $Ad^*$ actions are given by
\begin{align*}
    ad^*_{\hat{\Omega}} (\breve{\Pi}):=\breve{(\Pi \times \Omega)} \quad \text{and} \quad Ad^*_{R}(\breve{\Pi}):=\breve{(R^{-1}\Pi)}
\end{align*}
where $\Omega, \Pi \in \mathbb{R}^3$ and $R\in \mathbb{SO}(3)$. Using \eqref{integrator eg 1 left SO(3)} , the resulting integrator takes the form 
\begin{equation}
    \boxed{
    \begin{aligned}
        R_{k+1}&=R_k\exp(t\hat{\Omega}_k)\\
        d^{L*}_{t\hat{\Omega}_k}\exp(\breve{\Pi}_{k+1})&=\mathbb{I}(t\hat{\Omega}_k)\\
        \breve{\Pi}_{k+1}&=Ad^*_{\tau(t\hat{\Omega}_k)}(\breve{\Pi}_k) \label{integrator exp SO(3)}
    \end{aligned}
    }
\end{equation}
where $d^{L}\exp$ denotes the left logarithmic derivative of the exponential map. Substituting the explicit expression for the left logarithmic derivative of the exponential map (see \cite{iserles2000lie}), we obtain the fully explicit update rule
\begin{equation}
    \boxed{
    \begin{aligned}
        R_{k+1}&=R_k\exp(t\hat{\Omega}_k)\\
        \Bigl(I_3-\frac{\sin^2{\|t\Omega_k/2\|}}{2\|t\Omega_k/2\|^2}t\hat{\Omega}_k
        +\frac{\|t\Omega_k\|-\sin{\|t\Omega_k\|}}{\|t\Omega_k\|^3}t^2\hat{\Omega}_k^2\Big)\Pi_{k+1}&=t\mathcal{I}\Omega_k\\
        \Pi_{k+1}&=\exp{(t\hat{\Omega}_k)}\Pi_k \label{integrator exp SO(3) explicit}
    \end{aligned}
    }
\end{equation} 
where $I_3$ is the $3 \times 3$ identity matrix. The plots shown in figure \ref{fig_rigid_body} clearly highlight the structure-preserving properties of this integrator. 
\subsection{Numerical Integrators Using Cayley Map}
\label{subsec rigid body cay integrator}
We now choose the \emph{Cayley map} $Cay: \mathfrak{so}(3) \to \mathbb{SO}(3)$ as the local diffeomorphism $\tau$, defined by 
\begin{align*}
    Cay(\hat{x}) := \left(I_3-\hat{x}\right)^{-1}\left(I_3+\hat{x}\right)^{-1}
\end{align*}
for all $x \in \mathbb{R}^3$. Using \eqref{integrator eg 1 left SO(3)}, the corresponding integrator is given by
\begin{equation}
    \boxed{
    \begin{aligned}
        R_{k+1}&=R_k\,Cay(t\hat{\Omega}_k)\\
        d^{L*}_{t\hat{\Omega}_k}Cay(\breve{\Pi}_{k+1})&=\mathbb{I}(t\hat{\Omega}_k)\\
        \breve{\Pi}_{k+1}&=Ad^*_{Cay(t\hat{\Omega}_k)}(\breve{\Pi}_k) \label{integrator Cayley SO(3)}
    \end{aligned}
    }
\end{equation}
where $d^L \,Cay$ denotes the left logarithmic derivative of the Cayley map. Substituting the explicit expression for the left logarithmic derivative of the Cayley map (see, e.g., \cite{iserles2000lie}), we obtain the fully explicit update rule
\begin{equation}
    \boxed{
    \begin{aligned}
        R_{k+1}&=R_k(I_3-t\hat{\Omega})^{-1}(I_3+t\hat{\Omega}_k)\\
        \left(I_3-{t\hat{\Omega}_k}\right)\Pi_{k+1}&=\left(\frac{1+\|t\Omega_k\|^2}{2}\right)t\mathcal{I}\Omega_k\\
        \Pi_{k+1}&=(I_3-t\hat{\Omega})^{-1}(I_3+t\hat{\Omega}_k)\Pi_k \label{integrator Cayley SO(3) explicit}
    \end{aligned}
    }
\end{equation}
where $I_3$ is the $3 \times 3$ identity matrix. The plots shown in figure \ref{fig_rigid_body} clearly highlight the structure-preserving properties of the resulting integrators. 
\subsection{Comparison with Classical Integrators}
\begin{figure}
\centering
\begin{subfigure}{.45\textwidth}
  \centering
  \includegraphics[width=1\linewidth]{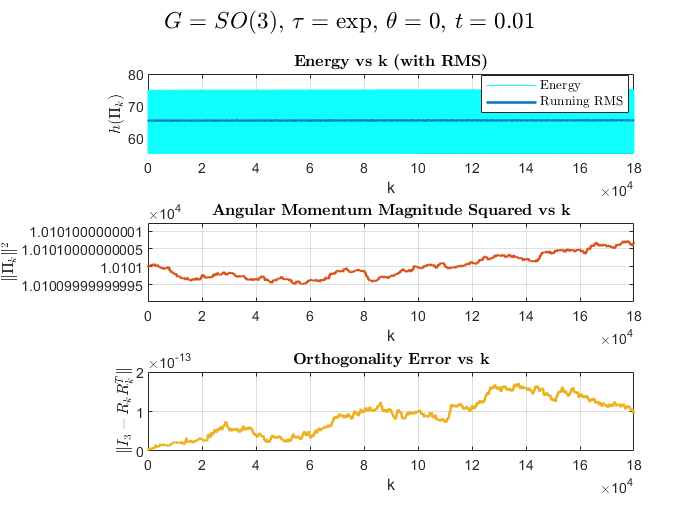}
  \caption{}
\end{subfigure}
\begin{subfigure}{.45\textwidth}
  \centering
  \includegraphics[width=1\linewidth]{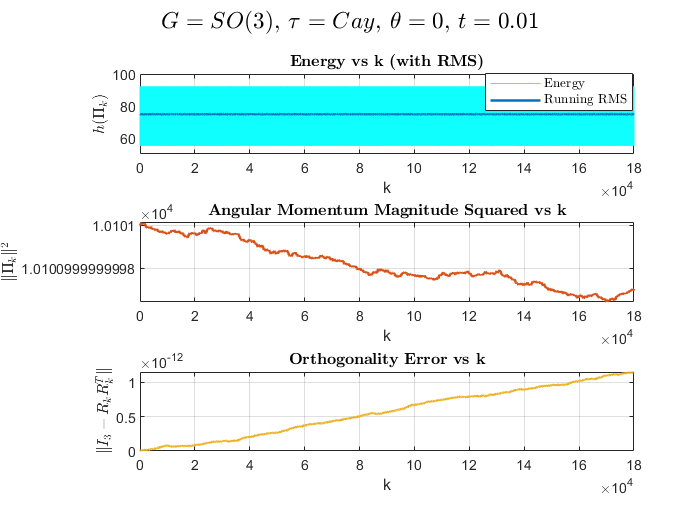}
  \caption{}
\end{subfigure}
\begin{subfigure}{.45\textwidth}
  \centering
  \includegraphics[width=1\linewidth]{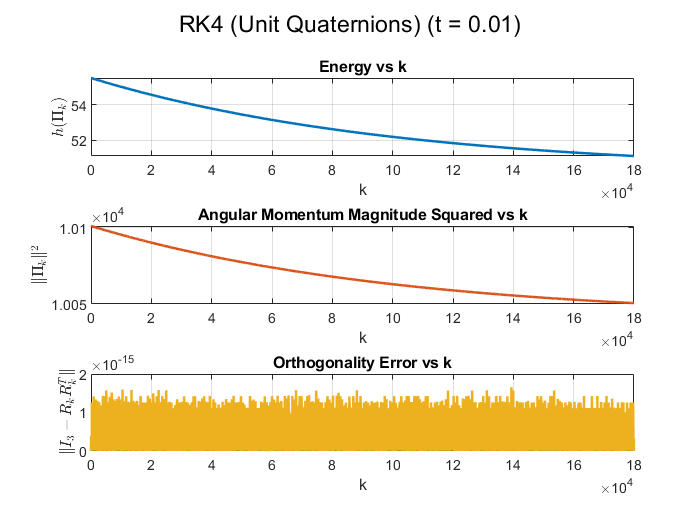}
  \caption{}
\end{subfigure}
\begin{subfigure}{.45\textwidth}
  \centering
  \includegraphics[width=1\linewidth]{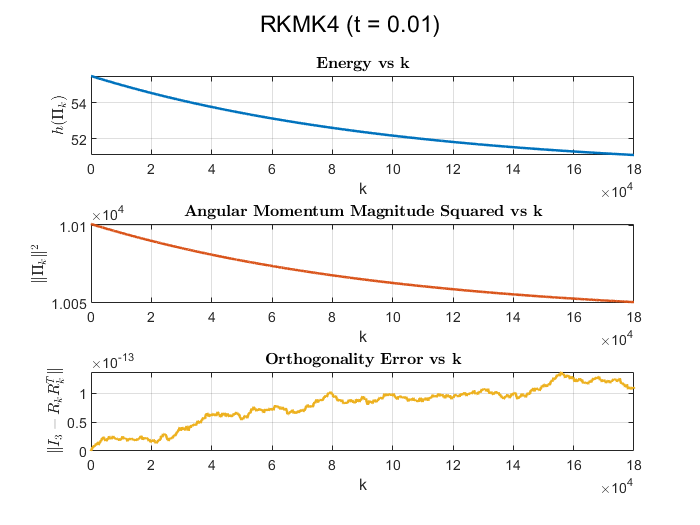}
  \caption{}
\end{subfigure}
\caption{(a) Exponential map-based integrator, (b) Cayley map-based integrator, (c) Runge-Kutta 4th order with projection-based integrator (unit quaternion), (d) Runge-Kutta-Munthe-Kaas 4th order-based integrator. All simulations were performed over a time horizon of 30 minutes for a rigid body with inertia tensor $\mathcal{I}=diag(1,10,100)$, initial conditions $\Pi_0=(1,1,1)$, $R_0=diag(1,1,1)$ and step size $h = 0.01$ seconds (in standard units).}
\label{fig_rigid_body}
\end{figure}

A widely used approach for numerically integrating rigid body dynamics is to represent the attitude using unit quaternions rather than rotation matrices, and to apply standard Runge–Kutta methods in the ambient space $\mathbb{R}^4$. After each time step, the quaternion is renormalized to enforce the unit-length constraint, thereby remaining on $\mathbb{S}^3 \subset \mathbb{R}^4$, the space of unit quaternions. From a geometric viewpoint, this corresponds to lifting the dynamics from $\mathbb{SO}(3)$ to its universal cover $Spin(3) \cong \mathbb{SU}(2)$, on which the rigid body motion is represented without singularities. Despite their simplicity and numerical robustness, such projection-based schemes do not respect the underlying Lie–Poisson structure of the rigid body equations and consequently fail to preserve key invariants, such as the Hamiltonian and the Casimir $\|\Pi\|^2$, as illustrated in figure \ref{fig_rigid_body}.
\\

Another popular class of methods consists of Lie group Runge–Kutta schemes, in which the numerical update is performed in the Lie algebra, exploiting its linear structure, while the exponential map is used to reconstruct the group element at each step. Although these methods exactly preserve the group constraint $R_k \in \mathbb{SO}(3)$, they are not, in general, Poisson or energy preserving. As a result, they also exhibit long-time drift in both the Hamiltonian and the squared magnitude of the body angular momentum, as shown in figure \ref{fig_rigid_body}. 
\section{Heavy Top}
The rotating rigid body pivoted at a point other than its center of gravity and subject to a uniform gravitational field in the vertical direction provides another important example of an Euler-Poincar{\'e} or Lie-Poisson system, now with advected parameters; see \cite{HOLM19981}. In this setting, the presence of gravity breaks the full rotational symmetry of the free rigid body, as shown in the figure below. This loss of symmetry can be systematically addressed by augmenting the configuration space to include the vertical unit vector expressed in the body frame, thereby enlarging the underlying group to the semidirect product $\mathbb{SE}(3) := \mathbb{SO}(3) \ltimes \mathbb{R}^3$. 
\begin{center}
\begin{tikzpicture}[scale = 0.6]
\shadedraw[rotate around={-45:(5,5)}] (5,5) ellipse (2 and 4);
\node at (6,8.5) {$m \, ,\,\mathcal{I}$};
\draw[thick,red,-stealth] (4,4)--(8,4);
\draw[thick,red,-stealth] (4,4)--(4,8);
\draw[thick,red,-stealth] (4,4)--(1.5,1.5);
\node at (1,1) {\textcolor{red}{$\mathcal{S}$}};
\draw[thick,blue,-stealth] (4,4)--(7,5.5);
\draw[thick,blue,-stealth] (4,4)--(2,7.5);
\draw[thick,blue,-stealth] (4,4)--(7.5,2);
\node at (8,2) {\textcolor{blue}{$\mathcal{B}$}};
\draw[fill] (5,5) circle (4pt);
\draw[ultra thick,-stealth, violet] (4,4)--(5,5);
\node at (4.5,5) {$\chi$};
\draw[ultra thick,-stealth, violet] (4,4)--(4,5);
\node at (3.7,5) {$\Gamma$};
\draw[thick,-stealth] (9,6)--(9,5);
\node at (9,4.5) {$g$};
\draw[->] (1.9,7.9) arc (113:90:5);
\node at (2.7,8.5) {$R$};
\end{tikzpicture}
\end{center}
The additional variable—the vertical unit vector in body coordinates—is referred to as an advected parameter, since it is constant in the spatial frame but evolves in the body frame according to the dynamics. Incorporating such parameters allows the equations of motion to be formulated in a reduced Euler–Poincaré or Lie–Poisson framework on the dual of a semidirect product Lie algebra. 

The kinetic energy of a heavy top rotating with angular momentum $\Pi \in \mathbb{R}^3$ expressed in the body frame is given by 
\begin{align*}
    T = \frac{1}{2}\Pi \cdot \mathcal{I}^{-1} \Pi = \frac{1}{2} \Pi^T \mathcal{I}^{-1} \Pi
\end{align*}
where $\mathcal{I}$ denotes the moment of inertia matrix. As before, the body angular velocity is defined by $\Omega := \mathcal{I}^{-1}\Pi \in \mathbb{R}^3$. The potential energy of the heavy top in a uniform gravitational field is
\begin{align*}
    V = mg \Gamma \cdot \chi = mg \Gamma^T \chi,
\end{align*}
where $\Gamma \in \mathbb{R}^3$ is the vertical unit vector expressed in the body frame, and $\chi \in \mathbb{R}^3$ denotes the vector from the pivot point to the center of gravity, also expressed in the body frame. The vector $\Gamma$ is an advected parameter, evolving dynamically in the body frame while remaining constant in the spatial frame. The total energy (Hamiltonian) of the system is therefore
\begin{align*}
    h(\Pi,\Gamma) = \frac{1}{2}\Pi \cdot \mathcal{I}^{-1} \Pi + mg \Gamma \cdot \chi.
\end{align*}
The corresponding Lie–Poisson equations with advected parameters can be written as
\begin{equation}
    \boxed{
    \begin{aligned}
        \dot{R} &= R \hat{\Omega} \\
        \Pi &= \mathcal{I} \Omega \\
        \dot{\Pi} &= \Pi \times \Omega + mg \Gamma \times \chi \\
        \dot{\Gamma} &= \Gamma \times \Omega
    \end{aligned}
    }
\end{equation}
where $R \in \mathbb{SO}(3)$ represents the orientation of the heavy top, while $\Omega \in \mathbb{R}^3$ and $\Pi \in \mathbb{R}^3$ denote the angular velocity and angular momentum in the body frame, respectively. In addition to the Hamiltonian, the vertical component of the body angular momentum $\Pi \cdot \Gamma$ is conserved, reflecting the residual symmetry about the vertical axis. Furthermore, the squared magnitude $\|\Gamma\|^2 = 1$ is also preserved. We now construct structure-preserving integrators for these equations along the lines of \eqref{integrator eg 1 left}, adapting the construction to the semidirect-product structure of the heavy top.
\subsection{Numerical Integrator Using Exponential Map}
We choose the Lie group exponential map as the local diffeomorphism $\tau$, which in this case coincides with the matrix exponential. For the $\mathbb{SE}(3)$ group, we construct the integrator using \eqref{integrator eg 1 left} as
\begin{equation}
    \boxed{
    \begin{aligned}
        (R_k^{-1}R_{k+1},R_k^{-1}(x_{k+1}-x_k)) &= \exp{(t\hat{\Omega}_k,tv_k)} \\
            d^{L*}_{(t\hat{\Omega}_k,tv_k)}\exp (\breve{\Pi_{k+1}+\Gamma_{k+1} \times x_{k+1}},\Gamma_{k+1}) &= \mathbb{I} (t\hat{\Omega}_k,tv_k) \\
            (\breve{\Pi}_{k+1},\Gamma_{k+1}) &= Ad^*_{\exp(t\hat{\Omega}_k,tv_k)}(\breve{\Pi}_k,\Gamma_k)
    \end{aligned}
    }
\end{equation}
and after substituting the explicit expressions for the left logarithmic derivative of the exponential, see \cite{kobilarov2011discrete, solà2021microlietheorystate}, we obtain the following update equations:
\begin{equation}
    \boxed{
            \begin{aligned}
                &R_{k+1} = R_k \exp{(t\widehat{\Omega}_k)} \\
                &x_{k+1} = x_k + R_k J(t\Omega_k)mg\chi \\
                &J(t\Omega_k)^TR_k^T(\Pi_k+\Gamma_k \times x_k) + Q(t\Omega_k,mg\chi)^TR_k^T\Gamma_k = t\mathcal{I}\Omega_k \\
                & \Pi_{k+1} =\exp{(t\hat{\Omega}_k)}(\Pi_k-\Gamma_k \times J(t\Omega_k)mg\chi) \\
                & \Gamma_{k+1} = \exp{(-t\hat{\Omega}_k)}\Gamma_k
            \end{aligned}
            }
\end{equation}
where $x \in \mathbb{R}^3$ is an auxiliary variable arising from the semidirect product structure. The matrix-valued functions 
\begin{align*}
    J(y) := I_3+\frac{\sin^2 \|y/2\|}{\|y/2\|^2}\hat{x} + \frac{\|y\|-\sin \|y\|}{\|y\|^3}\hat{y}^2
\end{align*}
and
\begin{align*}
    Q(y,z) := &\left(\frac{\sin^2\|y/2\|}{\|y/2\|^2}\right)\hat{z} + \left(\frac{\|y\|-\sin{\|y\|}}{\|y\|^3}\right)(\hat{y}\widehat{z}-\hat{z}\hat{y}) \\
            &+ \left(\frac{\|y\|\sin{\|y\|+2\cos{\|y\|-2}}}{\|y\|^4}\right)y^Tz\hat{y} \\
            &+ \left(\frac{3\sin{\|y\|-\|y\|\cos{\|y\|}-2\|y\|}}{\|y\|^5}\right)y^Tz\hat{y}^2
\end{align*}
for all $y,z \in \mathbb{R}^3$. The numerical results shown in figure \ref{fig_heavy_top} clearly demonstrate the structure-preserving properties of this integrator. The simulations were run for 30 minutes for a heavy top with inertia tensor $\mathcal{I} = diag(1,10,100)$ and initial conditions $\Pi_0 = (1,1,1)$, $R_0 = I_3$, $\Gamma_0 = (0,0,1)$.
\subsection{Numerical Integrators Using Cayley Map}
We now choose the Cayley map
\begin{align*}
    Cay : \mathfrak{so}(3) &\to \mathbb{SO}(3) \\
    \hat{x} &\mapsto (I_3 - \hat{x})^{-1}(I_3 + \hat{x})
\end{align*}
as the local diffeomorphism $\tau$. Using this choice, we construct a integrator via \eqref{integrator eg 1 left}, yielding
\begin{equation}
    \boxed{
    \begin{aligned}
        (R_k^{-1}R_{k+1},R_k^{-1}(x_{k+1}-x_k)) &= Cay{(t\hat{\Omega}_k,tv_k)} \\
            d^{L*}_{(t\hat{\Omega}_k,tv_k)} Cay (\breve{\Pi_{k+1}+\Gamma_{k+1} \times x_{k+1}},\Gamma_{k+1}) &= \mathbb{I} (t\hat{\Omega}_k,tv_k) \\
            (\breve{\Pi}_{k+1},\Gamma_{k+1}) &= Ad^*_{Cay(t\hat{\Omega}_k,tv_k)}(\breve{\Pi}_k,\Gamma_k)
    \end{aligned}
    }
\end{equation}
and after substituting the explicit expressions of the right logarithmic derivative of the Cayley map, see \cite{muller2021review}, we obtain the following update equations:
\begin{equation}
    \boxed{ 
            \begin{aligned}
                &R_{k+1} = R_k \, Cay{(t\hat{\Omega}_k)} \\
                &x_{k+1} = x_k + R_k^T(I_3+(I_3+t\hat{\Omega}_k)(I_3-t\hat{\Omega}_k)^{-1})mg\chi \\
                &2(I_3-t\hat{\Omega}_k)(R_k^T(\Pi_{k+1}+\Gamma_k \times x_k)+(R_k^T\hat{\Gamma}_{k+1}R_kmg\chi)=(1+\|t\Omega_k\|^2) t\mathcal{I}\Omega_k \\
                & \Pi_{k+1} =\exp{(t\hat{\Omega}_k)}(\Pi_k-\Gamma_k \times J(t\Omega_k)mg\chi) \\
                & \Gamma_{k+1} = \exp{(-t\hat{\Omega}_k)}\Gamma_k
            \end{aligned}
            }
\end{equation}
where $x \in \mathbb{R}^3$ is an auxiliary variable. The numerical results presented in figure \ref{fig_heavy_top} clearly demonstrate the structure-preserving properties of this integrator. In particular, the method nearly preserves the Casimir invariants $\Pi \cdot \Gamma$ and exactly preserves $\|\Gamma\|^2$, while exhibiting excellent long-term energy behavior.
\subsection{Comparison with Classical Integrators}
\begin{figure}
\centering
\begin{subfigure}{.45\textwidth}
  \centering
  \includegraphics[width=1\linewidth]{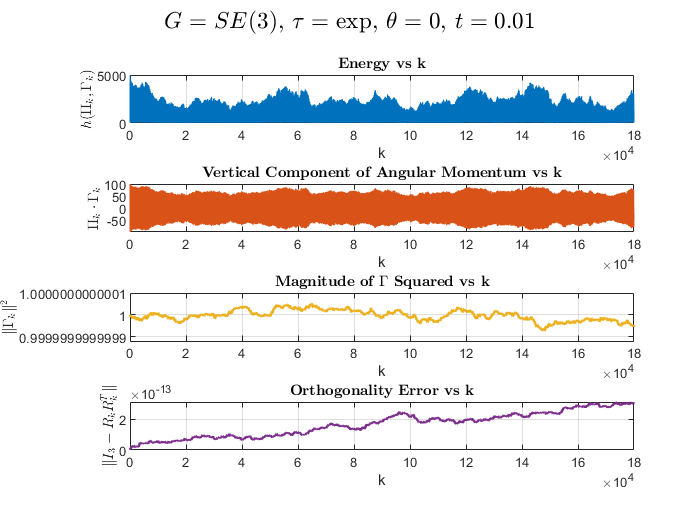}
  \caption{}
\end{subfigure}
\begin{subfigure}{.45\textwidth}
  \centering
  \includegraphics[width=1\linewidth]{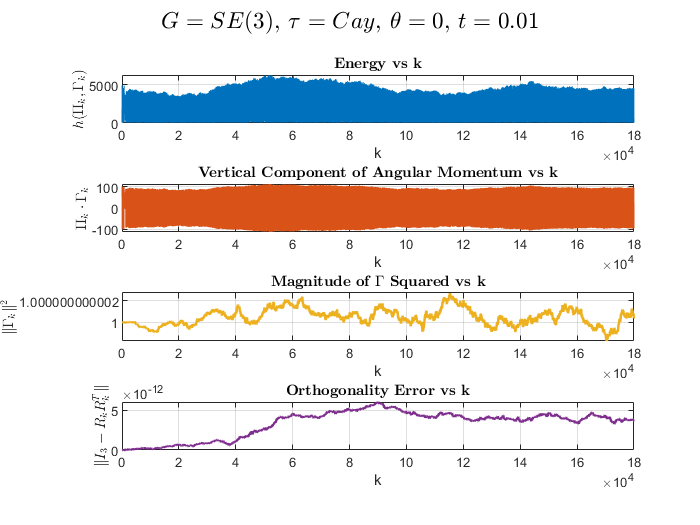}
  \caption{}
\end{subfigure}
\begin{subfigure}{.45\textwidth}
  \centering
  \includegraphics[width=1\linewidth]{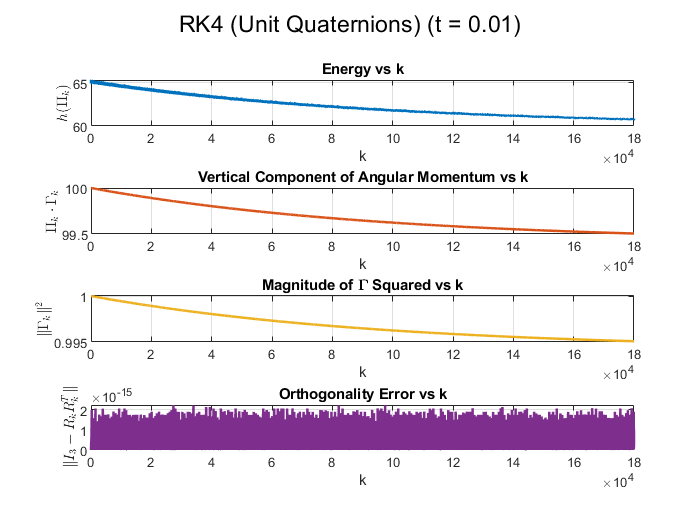}
  \caption{}
\end{subfigure}
\begin{subfigure}{.45\textwidth}
  \centering
  \includegraphics[width=1\linewidth]{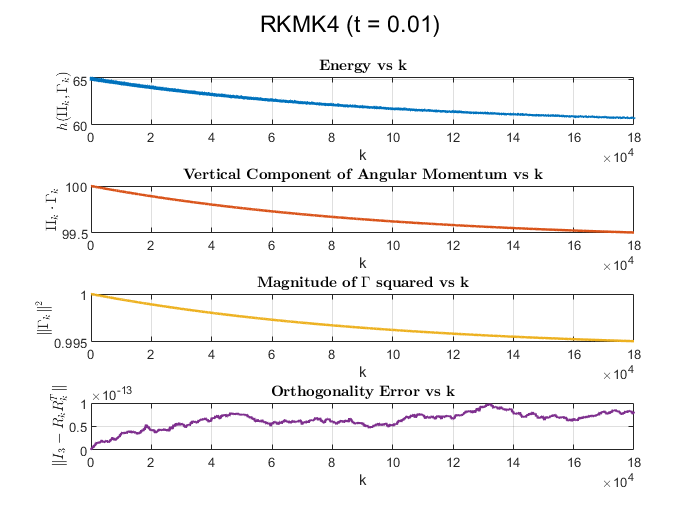}
  \caption{}
\end{subfigure}
\caption{(a) Exponential map-based integrator, (b) Cayley map-based integrator, (c) Runge-Kutta 4th order with projection-based integrator (unit quaternion), (d) Runge-Kutta-Munthe-Kaas 4th order-based integrator. All simulations were performed over a time horizon of 30 minutes for a heavy top with inertia tensor $\mathcal{I}=diag(1,10,100)$, initial conditions $\Pi_0=(1,1,1)$, $R_0=diag(1,1,1)$, $\Gamma_0=(0,0,1)$ and step size $h = 0.01$ seconds (in standard units).}
\label{fig_heavy_top}
\end{figure}

Once again, we compare the performance of the proposed integrator with the standard fourth-order Runge–Kutta (RK4) method implemented on unit quaternions along with the fourth-order Runge-Kutta-Munthe-Kaas (RKMK4) method. The plots presented in figure \ref{fig_heavy_top} clearly demonstrate that both the RK4 and RKMK4 schemes fail to preserve any of the geometric invariants of the system, including the Hamiltonian and the Casimir quantities.
\section{Quadrotor}
Quadrotors constitute a canonical example of an underactuated mechanical system evolving on a nonlinear configuration space, combining rigid body attitude dynamics with translational motion in Euclidean space. Their configuration space may be identified with the Lie group–manifold product $\mathbb{SO}(3)\times\mathbb{R}^3$, where the rotational dynamics evolve on $\mathbb{SO}(3)$ while the center-of-mass motion evolves in $\mathbb{R}^3$ as illustrated in the figure below. 
\begin{center}
\begin{tikzpicture}
\draw[gray] (3,2)--(5,3);
\draw[gray] (3.2,3.2)--(4.8,1.8);
\fill[fill opacity=0.2] (3,2) ellipse (15pt and 10pt);
\fill[fill opacity=0.2] (5,3) ellipse (15pt and 10pt);
\fill[fill opacity=0.2] (3.2,3.2) ellipse (15pt and 10pt);
\fill[fill opacity=0.2] (4.8,1.8) ellipse (15pt and 10pt);
\draw[thick,-stealth,red] (1.8,1)--(2.8,1);
\draw[thick,-stealth,red] (1.8,1)--(1.8,2);
\draw[thick,-stealth,red] (1.8,1)--(1.2,0.4);
\node at (2.5,0.8) {\footnotesize $\textcolor{red}{\mathcal{S}}$};
\draw[thick,-stealth,blue] (4,2.5)--(4.5,2.75);
\draw[thick,-stealth,blue] (4,2.5)--(4.5,2.065);
\draw[thick,-stealth,blue] (4,2.5)--(4.3,3.3);
\node at (4.1,2.1) {\footnotesize $\textcolor{blue}{\mathcal{B}}$};
\draw[thick,dotted,-stealth,red] (4,2.5)--(4.7,2.5);
\draw[thick,dotted,-stealth,red] (4,2.5)--(4,3.2);
\draw[thick,dotted,-stealth,red] (4,2.5)--(3.5,2);
\draw[->] (1.8,1)--(4,2.5);
\node at (2.7,1.4) {\scriptsize $q$};
\draw[->] (4,3.4) arc (90:45:10pt);
\node at (4.2,3.6) {\footnotesize $R$};
\draw[-stealth] (2,3)--(2,2.2);
\node at (1.8,2.6) {\footnotesize $g$};
\end{tikzpicture}
\end{center}
Unlike classical free rigid bodies or heavy tops, quadrotors are subject to continuous external forcing through thrust and control moments generated by the rotors. This attitude-dependent forcing breaks the full $\mathbb{SE}(3)$ symmetry of the system, preventing a formulation as a purely Euler–Poincar{\'e} or Lie–Poisson system. Nevertheless, the intrinsic geometric structure of the rotational dynamics remains that of a forced rigid body on $\mathbb{SO}(3)$, making quadrotors a natural testbed for Lie group and structure-preserving integrators. In this section, we exploit this decomposition to integrate the rotational and translational dynamics using geometric methods that respect the underlying manifold structure while accommodating external forces and moments.

The dynamics of a quadrotor are given by
\begin{equation}
    \boxed{
    \begin{aligned}
        \dot{R} &= R \widehat{\mathcal{I}^{-1}\Pi} \\
        \dot{\Pi} &= \Pi \times \mathcal{I}^{-1}\Pi + M \\        
        \dot{q} &= \frac{p}{m} \\
        \dot{p} &= -mg e_3 + FRe_3
    \end{aligned}
    }
\end{equation}
where the top two equations describe the rotational dynamics and the bottom two equations describe the translational dynamics. Here $R \in \mathbb{SO}(3)$ denotes the attitude of the quadrotor, $\Pi \in \mathbb{R}^3$ is the body angular momentum, and $\mathcal{I} \in \mathbb{R}^{3 \times 3}$ is the inertia tensor. The variable $q \in \mathbb{R}^3$ denotes the position of the centre of mass, $m \in \mathbb{R}$ is the mass, and $p \in \mathbb{R}^3$ is the linear momentum expressed in the inertial frame. Furthermore, $M \in \mathbb{R}^3$ denotes the net body moment generated by the rotors and $F \in \mathbb{R}$ is the total thrust. The dependence of $M$ and $F$ on rotor speeds and geometry can be found in standard quadrotor models, see for example \cite{lee2010geometric}. 

The presence of an attitude-dependent force - namely the thrust term - breaks the left-invariance of the dynamics on $\mathbb{SE}(3)$. Consequently, the system cannot be formulated as an Euler-Poincar{\'e} or Lie-Poisson system on $\mathbb{SE}(3)$. Nevertheless, the rotational dynamics are decoupled from the translational dynamics and correspond to those of a rigid body subject to an external torque. These equations may therefore be integrated using the Lie group methods developed in section \ref{sec rigid body}. The translational dynamics evolve on $T^*\mathbb{R}^3 \cong \mathbb{R}^3 \times \mathbb{R}^3$ and maybe integrated independently using the symplectic Euler method described in paragraph \ref{para symplectic Euler}. The resulting discrete-time integrator is given by
\begin{equation}
    \boxed{
    \begin{aligned}
        R_{k+1} &= R_k \tau(t \hat{\Omega}_k) \\
        d^{L*}_{t\hat{\Omega}_k} \tau (\breve{\Pi}_{k+1}) &= \mathbb{I}(t\hat{\Omega}_k) \\
        \breve{\Pi}_{k+1} &= Ad^*_{\tau(t\hat{\Omega}_k)}(\breve{\Pi}_k) + \breve{M}_k \\
        q_{k+1} &= q_k + t \frac{p_{k+1}}{m} \\
        p_{k+1} &= -p_k + t mg e_3 -tF_kR_ke_3
    \end{aligned} \label{quadrotor integrator}
    }
\end{equation}
where $t \in \mathbb{R}$ denotes the time step. Due to the presence of external forces and moments, the resulting integrator does not admit conserved quantities such as energy or momentum. Nevertheless, it respects the underlying configuration manifold structure and provides a geometrically consistent discretization of the quadrotor dynamics. While the translational dynamics are integrated here using the first-order symplectic Euler method, one could alternatively employ higher-order symplectic integrators, such as St{\"o}rmer–Verlet or higher-order splitting methods, to improve accuracy while still preserving the underlying symplectic structure of the translational phase space.
\subsection{Numerical Integrator Using Exponential Map}
From subsection \ref{subsec rigid body exp integrator}, the resulting update equations are
\begin{equation}
    \boxed{
    \begin{aligned}
        &R_{k+1}=R_k\exp(t\hat{\Omega}_k)\\
        &\Bigl(I_3-\frac{\sin^2{\|t\Omega_k/2\|}}{2\|t\Omega_k/2\|^2}t\hat{\Omega}_k
        +\frac{\|t\Omega_k\|-\sin{\|t\Omega_k\|}}{\|t\Omega_k\|^3}t^2\hat{\Omega}_k^2\Big)\Pi_{k+1}=t\mathcal{I}\Omega_k\\
        &\Pi_{k+1}=\exp{(t\hat{\Omega}_k)}\Pi_k + M_k\\
        &q_{k+1} = q_k + t \frac{p_{k+1}}{m} \\
        &p_{k+1} = -p_k + t mg e_3 -tF_kR_ke_3
    \end{aligned}
    }.
\end{equation}
In the special case of hovering flight, the performance of the integrator is shown in figure \ref{fig_quadrotor}. Although energy is not preserved due to external forcing, the method preserves angular momentum and respects the underlying manifold structure.
\subsection{Numerical Integrator Using Cayley Map}
\begin{figure}
\centering
\begin{subfigure}{.45\textwidth}
  \centering
  \includegraphics[width=1\linewidth]{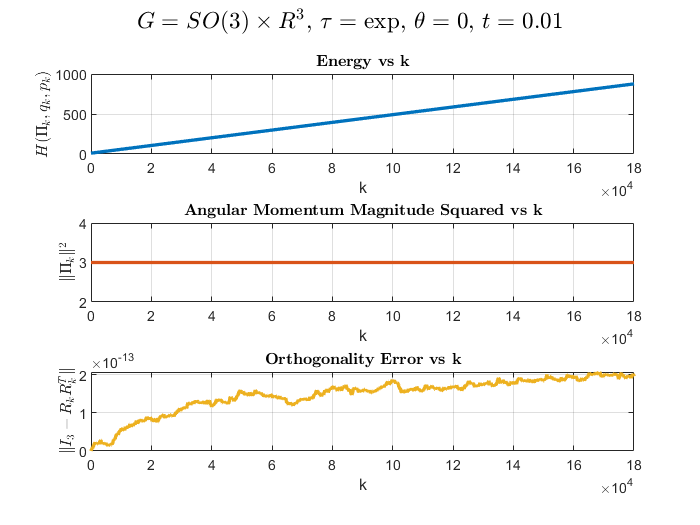}
  \caption{}
\end{subfigure}
\begin{subfigure}{.45\textwidth}
  \centering
  \includegraphics[width=1\linewidth]{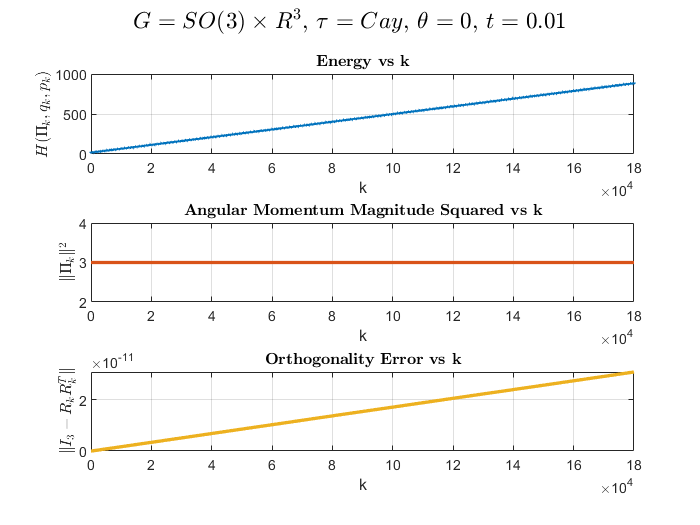}
  \caption{}
\end{subfigure}
\caption{(a) Exponential map-based integrator, (b) Cayley map-based integrator. All simulations were performed over a time horizon of 30 minutes for a quadrotor with inertia tensor $\mathcal{I}=diag(1,10,100)$, initial conditions $\Pi_0=(1,1,1)$, $R_0=diag(1,1,1)$, $p_0=(0,0,0)$, $q_0=(0,0,1)$ and step size $h = 0.01$ seconds (in standard units).}
\label{fig_quadrotor}
\end{figure}

From subsection \ref{subsec rigid body cay integrator}, the resulting update equations are
\begin{equation}
    \boxed{
    \begin{aligned}
        &R_{k+1}=R_k(I_3-t\hat{\Omega})^{-1}(I_3+t\hat{\Omega})\\
        &\left(I_3-t\hat{\Omega}\right)\Pi_{k+1}=\left(\frac{1+\|t\Omega_k\|^2}{2}\right)t\mathcal{I}\Omega_k\\
        &\Pi_{k+1}=(I_3-t\hat{\Omega})^{-1}(I_3+t\hat{\Omega})\Pi_k + M_k \\
        &q_{k+1} = q_k + t \frac{p_{k+1}}{m} \\
        &p_{k+1} = -p_k + t mg e_3 -tF_kR_ke_3
    \end{aligned}
    }.
\end{equation}
Again in the special case of hovering flight, the performance of the integrator is shown in figure \ref{fig_quadrotor}. Although energy is not preserved due to external forcing, the method preserves angular momentum and respects the underlying manifold structure just like before.

\bibliographystyle{plain}

\end{document}